\documentclass[times,12pt,preprint,nopreprintline]{elsarticle}

\usepackage{amssymb}

\usepackage[caption=false]{subfig} 
\usepackage{algorithm}
\usepackage{algorithmic}
\usepackage{amsfonts}
\usepackage{amsmath}
\usepackage{array}
\usepackage[T1]{fontenc}%
\usepackage{hyperref}
\hypersetup{colorlinks,linkcolor=black,urlcolor=black,citecolor=black}
\usepackage{bm}
\newtheorem{remark}{Remark}
\makeatletter
\def\tnotemark[#1]{\@for\mytmark:=#1\do{%
       \@ifundefined{X@\mytmark}{\expandafter\def\csname X@\mytmark\endcsname{0}}{}%
       \expandafter\ifcase\csname X@\mytmark\endcsname\or$^{\star}$\or
       $^{,\star\star}$\fi
    }%
}
\newenvironment{breakablealgorithm}
{%
		\begin{center}
			\refstepcounter{algorithm}%
			\hrule height.8pt depth0pt \kern2pt%
			\renewcommand{\caption}[2][\relax]{%
				{\raggedright\textbf{\ALG@name~\thealgorithm} ##2\par}%
				\ifx\relax##1\relax %
				\addcontentsline{loa}{algorithm}{\protect\numberline{\thealgorithm}##2}%
				\else %
				\addcontentsline{loa}{algorithm}{\protect\numberline{\thealgorithm}##1}%
				\fi
				\kern2pt\hrule\kern2pt
			}
		}{%
		\kern2pt\hrule\relax%
	\end{center}
}
\makeatother
\usepackage{placeins}
\usepackage{graphicx,epstopdf}

\begin{document}
	
	\begin{frontmatter}
		
		\title{A model order reduction based adaptive parareal method for time-dependent partial differential equations \tnoteref{tnote1}}%
		
		\tnotetext[tnote1]{This work was supported by  the National Natural Science Foundation of China under grant 92270206  and the Strategic Priority Research Program of the Chinese Academy of Sciences under grant XDB0640000.}
		
		\author[1]{Xiaoying Dai\corref{cor1}}
		\cortext[cor1]{Corresponding author:
			email: daixy@lsec.cc.ac.cn;}
		\author[2]{Miao Hu} 
		\ead{humiao@cpu.edu.cn}
		\author[1]{Shuwei Shen}
		\ead{shenshuwei@lsec.cc.ac.cn}
		
		\address[1]{State Key Laboratory of Mathematical Sciences, Academy of Mathematics and Systems Science, Chinese Academy of Sciences, Beijing 100190, China; and School of Mathematical Sciences,
			University of Chinese Academy of Sciences, Beijing 100049, China}
		\address[2]{China Pharmaceutical University, Nanjing 211198, China}

		\begin{abstract}
			In this paper, we propose a model order reduction based adaptive parareal method for time-dependent  partial differential equations. By using the data obtained by the fine propagator from  neighboring time subintervals and previous neighboring parareal iterations together with some model order reduction technique, we construct some POD subspaces locally in time and use them to construct the coarse propagator for each parareal iteration. We apply this new method to solve some 3D time-dependent advection-diffusion equations with the Kolmogorov flow and the ABC flow. Numerical results show that our method can achieve high accuracy with significant  speedup for long-term evolution problems.
		\end{abstract}
	\end{frontmatter}
	\section{Introduction}
	\label{sec:Introduction}
	
	With the increasing number and efficiency of cores in supercomputer platforms, there is a growing demand for scalability and parallelism of algorithms. However, when dealing with time-dependent  problems, the inherent sequential nature of the time direction becomes a bottleneck for scalability. To overcome this obstacle, researchers attempt to extend parallelism into the time direction. We refer to \cite {gander201550, gander2025time} and references therein for an overview of the development of parallel in time methods. One of the most widely used methods is the parareal method. 
	
	The parareal method was first proposed in \cite {lions2001resolution}. The main idea is to break the problem defined over the entire time interval into separate  subproblems on time subintervals. After that, some analysis on the convergence and the stability for the parareal method is given in \cite {bal2005convergence, gander2008nonlinear, gander2007analysis, staff2005stability}. Over the past few years, the parareal method has been widely used in many fields (see  \cite {maday2010parareal} for a review), e.g., fluid and structure problems \cite {farhat2003time},  the Navier-Stokes equations  \cite {fischer2005parareal},  reservoir simulation  \cite {garrido2005convergent},  the simulation of molecular dynamics  \cite {legoll2022adaptive}, time-periodic problems \cite {gander2013analysis} and so on. Furthermore, some modified parareal methods have been proposed to deal with the long time simulation of  hyperbolic problems and Hamiltonian systems \cite {bal2008symplectic,  dai2013symmetric, dai2013stable, farhat2003time}, which are shown to have better  convergence and stability.
	
	The parareal method involves two propagators: a fine propagator with high accuracy but high computational cost that computes in parallel on time subintervals, and a coarse propagator that is inexpensive to compute over the entire time interval. The choice of the fine propagator depends on the specific problem and the expected accuracy for the approximate solution. The fine propagator usually applies a classical time integrator (e.g., Euler or Runge-Kutta methods) with a fine step size, combined with a classical spatial discretization method that has fine spatial resolution. The classical spatial discretization methods include the finite element method (FEM) \cite{brenner2008mathematical}, the finite difference method (FDM) \cite{strikwerda2004finite} and so on. These methods often result in semi-discretized systems with very large dimensionality, particularly for problems with high spatial dimension, leading to computationally intensive tasks involving millions or billions of degrees of freedom.
	
	Meanwhile, the choice of the coarse propagator is more flexible and also has an important impact on the performance of the parareal method. Numerical analysis \cite{bal2005convergence, gander2008nonlinear, gander2007analysis} demonstrates the effect of the accuracy of the coarse propagator on the convergence speed of the parareal method: the higher the accuracy of the coarse propagator, the faster the parareal method converges. However, a coarse propagator with high accuracy often results in expensive computational cost. Therefore, developing a coarse propagator with adequate accuracy and low cost is crucial for improving the performance of the parareal method. In other words, a good coarse propagator usually needs to satisfy the following two requirements: sufficient accuracy and low computational cost. The most commonly used coarse propagator is the same scheme as that used in the fine propagator but with a coarser step size \cite{lions2001resolution}. Additionally, the coarse propagator can employ an explicit time discretization scheme \cite{nielsen2012feasibility} or a scheme with lower order accuracy \cite{blouza2011parallel}.  
	 In addition, there are some works on using neural networks \cite{ibrahim2023parareal,ibrahim2024space,pamela2024neural}, simplified physics modules \cite{baffico2002parallel,legoll2022adaptive,maday2003parallel} or model order reduction \cite{caldas2021modified, carlberg2019data, chen2014use, he-2010, iizuka2018influence, maday-2007} to construct the coarse propagators.  
	
	For the plain parareal method, the accuracy of both the coarse and fine propagator is usually fixed throughout all iterations. To further accelerate the convergence and improve the efficiency, several adaptive parareal methods have been proposed in recent years. In \cite{maday2020adaptive}, Maday and   Mula introduced an adaptive parareal method by improving the accuracy of the fine propagator during iterations incrementally, and also analyzed the convergence and parallel efficiency. In \cite{legoll2022adaptive}, Legoll et al. proposed an adaptive variant by dividing the entire time interval into smaller time-slabs adaptively, and applied it to the simulation of molecular dynamics trajectories.
	
	In this article, we aim to design an adaptive parareal method for time-dependent partial differential equations by updating the coarse propagator adaptively during the parareal iterations.  As mentioned earlier, the computational cost of the coarse propagator is determined by both temporal and spatial discretization. For temporal discretization, we employ the same time integrator used in the fine propagator but with a coarse step size. For spatial discretization, we note that in each iteration of the parareal method, we need to run the fine propagator on each time subinterval, which will provide a large amount of data. Our basic idea is to make use of the data to construct a subspace, which will serve as the discrete space for the coarse propagator. As the parareal iteration proceeds, the iterative solution becomes progressively  closer to the solution obtained by applying the fine propagator starting from the exact initial value. The subspace constructed from the data may thus capture the behavior of the fine-scale solution increasingly well. Therefore, we have reason to believe that as the iteration proceeds, the spatial discretization accuracy of the coarse propagator will be progressively closer to that of the fine propagator. Furthermore, to reduce the cost of the coarse propagator, we hope the dimension of the subspace constructed from the data can be as small as possible. To achieve this, we employ some model order reduction technique \cite{benner2015survey, chinesta2017model, maday2006reduced} to construct a subspace from the data. Specifically, we employ the proper orthogonal decomposition (POD) \cite {pinnau2008model} to perform the model order reduction in this article.
	
	The main steps for our adaptive parareal method for time-dependent PDEs are as follows: at first, using the data obtained from the fine propagator on the early time interval and the POD method, we initialize the coarse propagator by constructing its POD subspace for discretization. The initial values on each time subinterval are then obtained by sequentially propagating the coarse propagator over the entire time interval. During the $k$-th iteration with $k\geq1$, we first solve the subproblem by the fine propagator on each time subinterval and collect snapshots in parallel, and then update the POD subspace on each time subinterval in parallel using these snapshots, from which we can refine the coarse propagator iteratively. Besides, we augment these subspaces by adding the initial values of the subproblems on time subintervals. This step can improve the accuracy of the coarse propagator by removing the error introduced by Galerkin projection. Finally, we obtain more accurate numerical solutions by using this coarse propagator with higher accuracy. This process establishes a positive feedback loop, wherein the coarse propagator improves progressively, leading to increasingly accurate approximate solutions obtained by the parareal method.

	  The idea of building a coarse propagator by using  information from the corresponding fine propagator at previous iterations in parareal method  is not new.   The Krylov subspace parareal method \cite{farhat-2006, gander2008analysis} builds a dynamically enriched coarse propagator by using
	  all numerical solutions at the end of each time subinterval computed by the corresponding fine propagator during all  previous iterations.   The idea was then applied to a linear acoustic-advection system  and an explicit implementation suitable for hyperbolic flow problems was constructed in \cite{ruprecht2012explicit}. Inspired by the Krylov subspace parareal approach, in \cite{chen2014use}, Chen et al. replaced the full accumulated Krylov subspace with a subspace spanned by the reduced basis constructed during the parareal iterations, and  proposed a reduced basis parareal method. Besides, they extended the idea to nonlinear systems using empirical interpolation methods.    
	  In \cite{iizuka2018influence},  the authors  combined the time evolution basis method with the reduced basis parareal method  and showed that reducing the number of bases causes the saturation of the decrease in the error in the linear advection-diffusion equation. 	  
	  In \cite{caldas2021modified},  
	   the authors  proposed a modified parareal method based on the proper orthogonal decomposition-discrete empirical interpolation method (POD-DEIM) for solving  nonlinear shallow-water equations.  	
	   More recently, Gander et al. \cite{gander-2026} proposed a new parareal algorithm  with low-rank coarse solvers  based on localized approximations of the transfer operator that maps initial values for a  given time interval to the solution at the end of the interval.
	  
	  The key difference between our method and the existing works lies in the construction of the coarse propagator.  
	  In all these existing works, the subspace used for the coarse propagator is built by collecting information from all previous iterations over the entire time interval, and the subspace remains the same across all time subintervals. In contrast, our method uses snapshots from neighboring time subintervals and previous neighboring iterations to construct the subspace for each time subinterval, which makes  the subspaces used  for different time subintervals  different. Furthermore, our method augments each time local reduced subspace with the current parareal interface state to eliminate the corresponding projection error. This produces a coarse propagator whose spatial approximation progressively improves while retaining a low  computational cost.

	The remainder of this paper is organized as follows. We review some fundamental concepts of the parareal method in section \ref{sec:parareal}. In section \ref{sec: coarse}, we introduce the details on how to construct our adaptive parareal method, especially those on how to use data obtained from propagating the fine propagator in each time subinterval together with the  POD method to construct the POD subspace which will be used as the spatial discretization of the coarse propagator, and the analysis of the speedup and computational complexity for our adaptive parareal method. In section \ref{sec: examples}, we apply our method to some typical time-dependent partial differential equations, i.e., the advection-diffusion equations with three dimensional velocity fields, including both the Kolmogorov flow and the ABC flow. The numerical results indicate that our method performs well in simulating the long-term evolution of this kind of problem, especially for the advection-dominated cases. Finally, we present our concluding remarks in section \ref{sec: Conclusion}.
	\section{Preliminaries}
	\label{sec:parareal}
	In this paper, we use bold letters to denote matrices and vectors, and regular letters for functions and scalars. Some standard notation for Sobolev spaces is also used, see e.g., \cite{verfurth2013posteriori}. Let $V$ be a Banach space with norm $||*||_V$ and $V^*$ be the dual space of $V$. We first introduce the Banach space
	$$
	W^2\left(0,T; V, V^* \right) = \left\{ u \in L^2(0,T; V) \ \middle| \ \frac{\partial u}{\partial t} \in L^2(0,T; V^*) \right\},
	$$
	where $
	L^p(0,T;V) = \left\{ u: (0,T) \to V \,\middle|\, u \text{ is measurable and } \int_0^T \|u(t)\|_V^p \, dt < \infty \right\}
	$. $\frac{\partial u}{\partial t}$ is the partial derivative with respect to time in the distributional sense. We call $u\in C(0,T;V)$ if  
	$$
	\lim_{t\to t_0} \|u(t) - u(t_0)\|_V = 0, \quad \forall t_0 \in [0,T]
	$$
	with $\| u \|_{C(0,T;V)} = \max\limits_{t \in [0,T]} \| u(t) \|_V$. 
	
	We consider the following time-dependent partial differential equations with periodic boundary conditions
	\begin{equation}\label{eq0}
		\begin{aligned}
			\frac{\partial u\left(x_1,\ldots ,x_d,t\right)}{\partial t}+\mathcal{A}\left(t, u\left(x_1,\ldots ,x_d,t\right)\right)&=f\left(x_1,\ldots ,x_d\right),\quad \text{in}\ \Omega\times(0,T]\\
			\quad u\left(x_1,\ldots ,x_d,0\right)&=u_0\left(x_1,\ldots ,x_d\right),\quad \text{in}\ \Omega \\ 
			u(x_1 + L_1, x_2, \dots, x_d, t) &= u(x_1, x_2, \dots, x_d, t),\quad \text{on}\ \partial\Omega\times(0,T]\\
			\dots\\
			u(x_1, x_2, \dots, x_d + L_d, t)  &= u(x_1, x_2, \dots, x_d, t),\quad \text{on}\ \partial\Omega\times(0,T]\\
		\end{aligned}
	\end{equation}
	in a bounded space-time cylinder $\Omega\times[0,T]\subset \mathbb{R}^d\times[0,T],d\geq 2$, where $\Omega:=[0,L_1]\times \ldots\times[0,L_d]$. Let 
	$$
	V = \left\{ v \in H^1\left(\Omega\right) \ \middle| \ v \text{ satisfies the boundary conditions in \eqref{eq0}} \right\},
	$$  
	$u_0\in V$ and $f\in L^2(0,T;V)$. The final time $T$ is arbitrary but fixed in what follows.  $\mathcal{A}$ is a partial differential operator, possibly nonlinear, from $[0,T]\times H^1(\Omega)$ into its dual space $H^{-1}(\Omega)$. 
	
	Then we assume the problem is well-posed with sufficient regularity, and denote the solution of \eqref{eq0} as $u \in W^2\left(0,T; V,V^* \right)$. All methods introduced later are not limited to any specific boundary conditions, with only minor adjustments required in the definition of the bilinear form for different boundaries.
	
	We begin with a brief review of the parareal method, which was first proposed in \cite {lions2001resolution}. The exact propagator $\mathcal{E}: [0, T] \times [0, T] \times {V} \rightarrow {V}$ is defined as follows:
	$$
	u\left(t\right)=\mathcal{E}\left(s,  t,  u\left(s\right)\right), \quad s \leq t, 
	$$
	where $u(s)$ represents the solution at time $s$. Let 
	$$0=t_0<t_1<t_2<\ldots<t_n<t_{n+1}<\ldots<t_{N}=T
	$$
	be a uniform partition of the time domain $\left[0,  T\right]$,  where $t_n=n \cdot \Delta T$, $\Delta T=\frac{ T}{N}$. Then we have
	$$
	u\left(t_{n+1}\right)=\mathcal{E}\left(t_n,  t_{n+1},  u\left(t_n\right)\right),  \quad n=0,  \ldots,  N-1.
	$$
	
	Let the fine propagator $\mathcal{F}$ be an approximation of $\mathcal{E}$ with high resolution. Both the exact and fine propagators are essentially sequential.   Let $\mathcal{G}$ be a coarse propagator, which has lower computational cost and lower accuracy compared with the fine propagator $\mathcal{F}$. 
	
	Let $k$ denote the index of iteration. The parareal method \cite {lions2001resolution} builds a sequence of $N$-tuples $\{U_n^k\}$ iteratively by the following formula: 
	\begin{equation}\label{formula1}
		\begin{aligned}
			 \left\{
			\begin{aligned}
				U_{n+1}^0 &= \mathcal{G}\left(t_n,  t_{n+1}, U_n^0\right),\quad  U_0^0=u_0, \quad n=0,1 \ldots,N-1, \\
				U_{n+1}^{k} &= \mathcal{G}\left(t_n,  t_{n+1},  {U}^{k}_n\right) +\mathcal{F}\left(t_n,  t_{n+1},  U_n^{k-1}\right) \\
				&\quad -{\mathcal{G}}\left(t_n,  t_{n+1}, {U}^{k-1}_n\right),\quad U_0^{k}=u_0, \quad n=0,1 \ldots,N-1,\quad  k=1,2 \ldots.
			\end{aligned}\right. \\
		\end{aligned}
	\end{equation}
	
	This method breaks the original problem defined over the entire time interval $[0,T]$ into $N$ separate subproblems on time subintervals $[t_n,t_{n+1}]$. During the $k$-th iteration, we first obtain $\mathcal{F}\left(t_n,  t_{n+1},  U_n^{k-1}\right)$ by propagating the fine propagator starting from $t_n$ with initial value  $U_n^{k-1}$ to $t_{n+1}$ in each time subinterval $[t_n, t_{n+1}]$. Since $U_n^{k-1}$ are obtained from the last iteration, the calculation of $\mathcal{F}\left(t_n,  t_{n+1},  U_n^{k-1}\right)$ can be carried out in parallel. After this step, since $\mathcal{G}\left(t_n, t_{n+1}, {U}^{k-1}_n\right)$ can be obtained from the previous iteration, we only need to calculate $\mathcal{G}\left(t_n, t_{n+1}, {U}^{k}_n\right)$ and update $U_{n+1}^{k}$ sequentially.
	
	It is obvious that within at most $N$ iterations, the parareal method will converge to the solution obtained sequentially by the fine propagator in the whole time domain \cite{baffico2002parallel,Bal2002}. However, it is desirable for the method to be able to achieve the convergence in significantly fewer than $N$ iterations. Otherwise, the accumulated computational cost of coarse propagators will exceed the computational cost of the traditional sequential computation, resulting in no parallel speedup (the rigorous definition and analysis for the speedup will be detailed in Section \ref{sec: Complex}).
	\section{Adaptive parareal method based on model order reduction}
	\label{sec: coarse}
	
	In this section, we propose our adaptive parareal method. Previous studies \cite{bal2005convergence, gander2008nonlinear, gander2007analysis} show that the convergence rate of the parareal method depends critically on the accuracy of the coarse propagator. Higher accuracy of the coarse propagator yields faster convergence of the parareal method. However, high accuracy also means high computational cost. This trade-off highlights the importance of developing computationally efficient coarse propagators that maintain sufficient accuracy for improving the performance of the parareal method.
	
	Usually, the resolution of the coarse propagator remains fixed throughout all iterations in the existing works. To accelerate the convergence, we propose a method that can improve the accuracy of the coarse propagator adaptively while maintaining low computational cost. This adaptive parareal method can be expressed as the following general formula:
	\begin{equation}\label{formula 2}
		\begin{aligned}
			\left\{
			\begin{aligned}
				{U}_{n+1}^0 &= \mathcal{G}_{0}\left({t}_n,  {t}_{n+1},  {U}_n^0\right),\quad  {U}_0^0 = u_0,\quad n=0,1 \ldots,N-1, \\
				{U}_{n+1}^{k} &= {\mathcal{G}}_{k}\left({t}_n,  {t}_{n+1}, {U}^{k}_n\right) +\mathcal{F}\left({t}_n,  {t}_{n+1},  {U}_n^{k-1}\right) \\
				&\quad -{\mathcal{G}}_{k}\left({t}_n,  {t}_{n+1},{U}^{k-1}_n\right) , \quad {U}_0^{k} = {u}_0,\quad n=0,1 \ldots,N-1, \quad  k=1,2 \ldots.
			\end{aligned}\right. 
		\end{aligned}
	\end{equation}
	
	Unlike the plain parareal method \eqref{formula1}, where the coarse propagator $\mathcal{G}$ is not changed during all iterations, here, at different iterations, the coarse propagator may be different. Most importantly, the accuracy of the coarse propagator is expected to increase adaptively as the iteration proceeds.  We notice that during each parareal iteration, the propagation of the fine propagator on each time subinterval produces a large amount of data, which is in fact an approximation to the solution obtained by applying the fine propagator starting from the exact initial value, with increasing accuracy over iterations. Therefore, we consider to use the data to construct the spatial discretization of the coarse propagator. We believe that as the iteration proceeds, the spatial discretization accuracy of the coarse propagator can approach that of the fine propagator's spatial discretization progressively. Moreover, to minimize the computational cost of the coarse propagator, we aim to construct the lowest-dimensional subspace from the data. Thanks to the model order reduction methods \cite{benner2015survey, chinesta2017model, maday2006reduced}, we can achieve this successfully. Specifically, we employ the proper orthogonal decomposition (POD) \cite {pinnau2008model} to perform the model order reduction in this article.
	
	Roughly speaking, the main process for constructing the coarse propagator is as follows. In the $0$-th iteration, we first apply the fine propagator to the initial part of the time domain to generate the data, which is then used to construct the POD subspace for the coarse propagator. We then apply the coarse propagator over the entire time interval to provide the initial values for the subproblems on  time subintervals. In subsequent iterations, we parallelize the fine propagator across all time subintervals and collect snapshots. These snapshots are then used to update the POD subspace in parallel, leading to a refined coarse propagator. Furthermore, to improve accuracy,  we augment the updated subspaces with the initial values from each time subinterval, thereby eliminating the Galerkin projection errors.
	
	In the  remainder of this section, we will introduce the details on how to apply the adaptive parareal method \eqref{formula 2} to deal with problem \eqref{eq0}. 
	
	We take an inner product with $v\in V$ in \eqref{eq0} and define $\left(\cdot,\cdot\right)$ as the inner product in $L^2\left(\Omega\right)$. We then obtain the variational form of \eqref{eq0} as follows: find $u\in W^2\left(0,T; V,V^* \right)$ such that $u(0)=u_0 \in V^*$ and for almost every $t\in (0,T)$,
	\begin{equation}\label{eq2}
		\left(\frac{\partial u}{\partial t},  v\right)+\left(\mathcal{A}\left(t,u\right),  v\right)=(f,v), \quad \forall v \in {V}. 
	\end{equation}
	
	We first introduce the temporal discretization of \eqref{eq2}. In this article, we choose the implicit Euler scheme \cite{acary2010implicit} as the temporal discretization for both the fine and coarse propagators, but with time step sizes $\delta t$ and $dT$ respectively. We take the fine propagator as an example to introduce the details of the temporal discretization. Let $u_n\in V$ be the approximation of $u(n \delta t)$. Applying the implicit Euler scheme to  \eqref{eq2}, then we can get the following semi-discretization scheme:
	\begin{equation}\label{semi-dis-eq}
		\left(\frac{u_n-u_{n-1}}{\delta t},  v\right)+\left(\mathcal{A}\left(n \delta t, u_n\right),  v\right)=(f_n,v), \quad \forall v \in {V}. 
	\end{equation}
	The temporal discretization of \eqref{eq2} by coarse propagator is the same as \eqref{semi-dis-eq} but with time step $\delta t$ being replaced by $dT$. 
	
	In the following part of this section, we will introduce the spatial discretization of the fine propagator and coarse propagator respectively.
	\subsection{Spatial discretization of the fine propagator}
	\label{sec: construction}
	For simplicity, the iteration index $k$ will be omitted in this subsection. For the spatial discretization of the fine propagator, we choose the standard finite element discretization \cite{brenner2008mathematical}. Let $\mathcal{T}_h$ be a family of nested, shape-regular conforming meshes over $\Omega$ with mesh size $h$, i.e., there exists a constant $\gamma^*$ such that
	$$
	\frac{h_\tau}{\rho_\tau} \leq \gamma^*, \quad \forall \tau \in \mathcal{T}_h,
	$$
	where $h_\tau$ is the diameter of $\tau$ for each $\tau \in \mathcal{T}_h$, $\rho_\tau$ is the diameter of the biggest ball contained in $\tau$, and $h=\max \left\{h_\tau: \tau \in \mathcal{T}_h\right\}$. Denote the finite element space on $\Omega$ corresponding to mesh $\mathcal{T}_h$ as $V_h$,
	$$
	V_h=\left\{v_h:\left.v_h\right|_\tau \in \mathbb{P}_\tau, \forall \tau \in \mathcal{T}_h \text { and } v_h \in C^0\left( \bar{\Omega}\right)\right\} \text {, }
	$$
	where $\mathbb{P}_\tau$ is a set of polynomials on element $\tau$. Let $\left\{\phi_{h,  1},  \phi_{h,  2}, \ldots, \phi_{h,  N_g}\right\}$ be a basis of $V_h$, where $N_g$ is the dimension of $V_h$. Let
	\begin{equation}\label{basis}
		\mathbf{\Phi}_h: =\left(\phi_{h,  1},  \phi_{h,  2},  \ldots,  \phi_{h,  N_g}\right),
	\end{equation}
	and $u_{h, n}\in V_h$ be the approximation of $u(n \delta t)$.
	We then obtain the following full-discretized problem: find $\{u_{h, n}\}_{n=1,2,\ldots,N}$ with each $u_{h,n} \in V_h$, such that 
	\begin{equation}\label{eq5}
		\left(\frac{u_{h, n}-u_{h, n-1}}{\delta t},  v_h\right)+\left(\mathcal{A}(n\delta t, u_{h, n}), v_h\right)=(f_n,v_h), \quad v_h \in V_h.
	\end{equation}
	Note that each $u_{h, n}$ can be expressed as
	\begin{equation}\label{eq6}
		u_{h, n}=\sum_{i=1}^{N_g} a_{n, i} \phi_{h,  i}, 
	\end{equation}
	where $a_{ n, i}\in \mathbb{R}$. Inserting \eqref{eq6} into \eqref{eq5}, and choosing $v_h=\phi_{h,  j}$,  $j=1, 2,  \ldots,  N_g$, we obtain 
	\begin{equation}\label{eq7}
		\left(\sum_{i=1}^{N_g} a_{n, i} \phi_{h,  i}-\sum_{i=1}^{N_g} a_{n-1, i} \phi_{h,  i},  \phi_{h,  j}\right)+\delta t\cdot \left(\mathcal{A}(n\delta t,\sum_{i=1}^{N_g} a_{n, i} \phi_{h,  i}), \phi_{h,  j}\right)=\delta t\cdot \left(f_n,\phi_{h,  j}\right).
	\end{equation}
	Denoting 
	$$
	\begin{aligned}
		& \mathbf{u}_{h,n}=\left(a_{n,  1},  a_{n,  2},  \ldots,  a_{n,  N_g}\right)^{\top},\quad \mathbf{A}_{h,n,ij}= \left(\phi_{h,  j},\phi_{h,  i}\right)+\delta t\cdot\left(\mathcal{A}\left(n \delta t, \phi_{h,  j}\right),  \phi_{h,  i}\right),
	\end{aligned}
	$$
	$$
	\begin{aligned}
		& \mathbf{b}_{h,n,i}= \delta t\cdot\left(f_n, \phi_{h,  i}\right),
	\end{aligned}
	$$
	we then obtain the following algebraic form
	\begin{equation}\label{algform}
		\mathbf{A}_{h, n} \mathbf{u}_{h, n}=\mathbf{b}_{h, n}+\mathbf{M} \mathbf{u}_{h,n-1}, 
	\end{equation}
	where $\mathbf{A}_{h,n}=\left(\mathbf{A}_{h,n,ij}\right)_{N_g\times N_g}$, $\mathbf{b}_{h, n}=\left(\mathbf{b}_{h,n,i}\right)_{N_g}$ , and $\mathbf{M}=\left(\phi_{h,  j},\phi_{h,  i}\right)_{N_g\times N_g}$ is the mass matrix.
	
	\subsection{Adaptive spatial discretization of the coarse propagator}
	\label{sec: Adaptive}
	It is well-known that to make the parareal method efficient, the coarse propagator should be cheap enough. Therefore,  we do not employ the spatial discretization used in the fine propagator, i.e., the finite element method, for the coarse propagator. Instead, following \cite{dai2024augmented, dai2020two}, we employ the POD method for the spatial discretization of the coarse propagator. That is, during the propagation using the coarse propagator, we will discretize the problem in a subspace spanned by some POD modes, whose dimension is usually much smaller than that of the finite element space. As stated in the Introduction, to improve the performance of the parareal method, we hope the accuracy of the coarse propagator can be as high as possible. In the following of this subsection, we will introduce the details on how we update the POD subspace to improve the accuracy of the spatial discretization of the coarse propagator adaptively. As mentioned at the beginning of this section, the basic idea for constructing the spatial discretization for the coarse propagator is to make use of the data obtained by the fine propagator in each parareal iteration, together with some model order reduction technique, to update the subspace for the spatial discretization. Our introductions are separated into two cases: spatial discretization of the coarse propagator at the beginning of the parareal iteration($k=0$) and the adaptive update of the subspace for spatial discretization of the coarse propagator for each $k$-th parareal iteration with $k\geq1$. 
	
	{\bf Case I: Spatial discretization of the coarse propagator for $k$-th parareal iteration with $k=0$}
	
	For the $0$-th iteration, we use the POD method to do the spatial discretization of the coarse propagator. That is, we first propagate \eqref{eq0} on time interval $[0,T_0]$ using the fine propagator, then use the snapshots obtained to construct some POD subspace, and use this POD subspace as the spatial discretization of the coarse propagator $\mathcal{G}_0$. After that, on time interval $[T_0, T]$, we propagate \eqref{eq0} by the coarse propagator $\mathcal{G}_0$. More details on how to use POD method to solve the time-dependent partial differential equations can be found in \cite{dai2024augmented, dai2020two}.
	
	Now we introduce the details for each step. At first, we discretize \eqref{semi-dis-eq} on $\left[0,  T_0\right]$ by the fine propagator, and collect the numerical solutions at different times $t=0,\delta M\delta t,\ldots,\left(n_s-1\right)\delta M\delta t$ as a snapshot matrix $\mathbf{W}_h\in \mathbb{R}^{N_g\times n_s}$. Here, $\delta M$ is an integer parameter and $n_s=\lfloor\frac{T_0}{\delta t \cdot \delta M} \rfloor+1 $, where $\lfloor \cdot \rfloor$ denotes rounding down. The approximation of the solution at $T_0$ is denoted as $\hat{u}_{0}$. Then we use $\mathbf{W}_h$ to construct the POD modes.  For details on how to construct the POD modes from existing data, please refer to \cite{volkwein2011model}. We first solve the eigenvalue problem of the symmetric positive definite matrix \((\mathbf{W}_h)^{\top}\mathbf{M}\mathbf{W}_h\) and obtain the eigenvalues \(\lambda_1 \geq \lambda_2 \geq \ldots \geq \lambda_r > \lambda_{r+1} = \ldots = \lambda_{n_s} = 0\), together with the eigenvectors corresponding to the nonzero eigenvalues, denoted by \(\mathbf{Y} = [\mathbf{y}_1, \mathbf{y}_2, \ldots, \mathbf{y}_r] \in \mathbb{R}^{n_s \times r}\). Then we obtain the POD modes $\mathbf{R} = [\mathbf{r}_1, \mathbf{r}_2, \ldots, \mathbf{r}_r] \in \mathbb{R}^{N_g \times r}$ by $\mathbf{r}_i= \frac{1}{\sqrt{\lambda_i}}\mathbf{W}_h\mathbf{y}_i$.  Setting $\gamma_1\in(0, 1]$ as a given parameter, the dimension of the POD subspace is
	$$
	m=\min \left\{k \mid \sum_{i=1}^k \sqrt{\lambda_i} \geq\gamma_1 \cdot \sum_{i=1}^{r} \sqrt{\lambda_i}\right\}. 
	$$
	Denote $\widetilde{\mathbf{R}}=\mathbf{R}[: ,  1:  m]$ as the first $m$ columns of $\mathbf{R}$. Then we obtain the POD modes
	$$
	\mathbf{\Psi}_h: =\left(\psi_{ h, 1},  \psi_{ h, 2},  \ldots,  \psi_{h,m}\right)=\mathbf{\Phi}_h\widetilde{\mathbf{R}},
	$$
	and its corresponding POD subspace ${V}^0_{\mathrm{POD}}  = \text{span}\{\psi_{ h, 1},  \psi_{ h, 2},  \ldots,  \psi_{h,m}\}$, here $\mathbf{\Phi}_h$ is defined in \eqref{basis}. Then, we finish the construction of the POD subspace ${V}^0_{\mathrm{POD}}$ for the spatial discretization of the coarse propagator $\mathcal{G}_{0}$. We summarize the process for constructing POD modes as a routine \texttt{POD\_Mode} $\left(\mathbf{W}_h, \gamma, {\bf M}, \Phi_h,  m,  \mathbf{\Psi}_h,  \widetilde{\mathbf{R}}\right)$ in Algorithm \ref{PODMode}, see also in \cite{dai2024augmented,dai2020two}. 
	
	\begin{algorithm}
		\flushleft \caption{\texttt{POD\_Mode}($\mathbf{W}_h, \gamma, {\bf M}, \mathbf{\Phi}_h,  m,  \mathbf{\Psi}_h,  \widetilde{\mathbf{R}}$)}
		\label{PODMode}
		\hspace*{0.02in} {\bf Input: }$\mathbf{W}_h,  \gamma, {\bf M}, \mathbf{\Phi}_h=\left(\phi_{h,  1},  \phi_{h,  2},  \ldots,  \phi_{h,  N_g}\right)$.\\
		\hspace*{0.02in} {\bf Output: }$ m$,  POD modes $\mathbf{\Psi}_h=\left(\psi_{h,  1},  \psi_{h,  2},  \ldots,  \psi_{h,  m}\right), \widetilde{\mathbf{R}}$.
		\begin{algorithmic}[1]
			\STATE Solve the eigenvalue problem of \((\mathbf{W}_h)^{\top}\mathbf{M}\mathbf{W}_h\) and obtain $\mathbf{R}$.
			\STATE Let $m=\min \left\{k \mid \sum_{i=1}^k \sqrt{\lambda_i} \geq\gamma_1 \cdot \sum_{i=1}^{r} \sqrt{\lambda_i}\right\}$.
			\STATE $\widetilde{\mathbf{R}}=\mathbf{R}[: ,  1:  m]$.
			\STATE $\mathbf{\Psi}_h=\left(\psi_{h,  1},  \psi_{h,  2},  \ldots,  \psi_{h,  m}\right)=\mathbf{\Phi}_h \widetilde{\mathbf{R}}$.
		\end{algorithmic}
	\end{algorithm}
	
	Next, we show how to discretize \eqref{semi-dis-eq} in $V^0_{\mathrm{POD}}$ for $t\geq T_0$. We first set $u_{ l,\mathrm{POD}}=\sum_{j=1}^m {\alpha}_{l,  j} \psi_{h,  j}$, as the approximation of $u(T_0+l dT)$, and $\mathbf{u}_{l,\mathrm{POD}}=\left({\alpha}_{ l, 1}, \ldots, {\alpha}_{l, m} \right)^{\top}$, ${\alpha}_{l,j} \in \mathbb{R}$. Defining
	\begin{equation*} 
	\begin{aligned}
		& \mathbf{A}_{l,\mathrm{POD},ij}= \left(\psi_{h,  j},\psi_{h,  i}\right)+dT\cdot\left(\mathcal{A}\left(T_0+l dT, \psi_{h,  j}\right),  \psi_{h,  i}\right),\\
		& \mathbf{b}_{l,\mathrm{POD},i}= dT\cdot\left(f(T_0+l dT), \psi_{h,  i}\right),
	\end{aligned}
	\end{equation*} 
	
	we arrive at the reduced algebraic form 
	\begin{equation}\label{eq8}
		\mathbf{A}_{l,\mathrm{POD}} \mathbf{u}_{l,\mathrm{POD}}=\mathbf{b}_{l,\mathrm{POD}}+{\mathbf{M}}_{\mathrm{POD}}\mathbf{u}_{l-1,\mathrm{POD}}, 
	\end{equation}
	where $ \mathbf{A}_{l,\mathrm{POD}} =\left(\mathbf{A}_{l,\mathrm{POD},ij} \right)_{m\times m}$, $\mathbf{b}_{l,\mathrm{POD}}=\left(\mathbf{b}_{l,\mathrm{POD},i} \right)_{m}$ and ${\mathbf{M}}_{\mathrm{POD}}=\left(\psi_{h,  j},\psi_{h,  i}\right)_{m\times m}$. Considering that the POD modes can be expressed in $\mathbf{\Psi}_h=\mathbf{\Phi}_h\widetilde{\mathbf{R}}$, we can also write the equation \eqref{eq8} as  
	\begin{equation}\label{eq00}
		\widetilde{\mathbf{R}}^{\top}\mathbf{A}_{h, l}\widetilde{\mathbf{R}} \mathbf{u}_{l,\mathrm{POD}}=\widetilde{\mathbf{R}}^{\top}\mathbf{b}_{h, l}+\widetilde{\mathbf{R}}^{\top}\mathbf{M}\widetilde{\mathbf{R}}\mathbf{u}_{l-1,\mathrm{POD}}.
	\end{equation}
	
	{\bf Case II: Adaptive update of the spatial discretization of the coarse propagator for $k\geq1$}
	
	We now introduce how to update the POD subspace adaptively which is used for the spatial discretization of the coarse propagator during each parareal iteration. The construction of the POD subspace used in the coarse propagator $\mathcal{G}_{0}$ only uses the snapshots obtained on $[0,T_0]$. The accuracy of the coarse propagator with this POD subspace for the spatial discretization is usually rather limited. As mentioned previously, we notice that during the parareal iterations, the propagation of the fine propagator on each time subinterval produces a large amount of data, which is in fact an  approximation to the solution obtained by applying the fine propagator starting  from the exact initial value  with increasing accuracy during the iteration. Therefore, we consider to use the data to update the POD subspace. We believe that as the iteration proceeds, the accuracy of the coarse propagator using this updated POD subspace can  approach that of the fine propagator's spatial discretization progressively.  
	
	Let
	$$
	T_0=\hat{t}_0<\hat{t}_1<\hat{t}_2<\ldots<\hat{t}_n<\hat{t}_{n+1}<\ldots<\hat{t}_N=T
	$$
	be a partition of the time domain $\left[T_0,T\right]$.  The solution obtained by the fine propagator at $T_0$ is denoted as $\hat{u}_0$. Meanwhile, the solution obtained by applying the fine propagator starting from the exact initial value is denoted as $U_n$. We denote $\mathcal{G}_{k,n}$ the coarse propagator on the time subinterval $\left[\hat{t}_n, \hat{t}_{n+1}\right]$ during the $k$-th iteration, whose spatial discretization is carried out in the subspace $V^k_{n,\mathrm{POD}}$. Our goal is to construct subspaces  $\{V^k_{n,\mathrm{POD}}\}_{k}$ such that the spatial discretization of the coarse propagator on these subspaces yields approximations with increasingly high accuracy as $k$ increases. From the data obtained by calculating $\mathcal{F}\left(\hat{t}_n, \hat{t}_{n+1}, U_n^{k-1}\right)$, we collect the snapshots on each time subinterval in parallel. These snapshots can be used to construct new POD subspaces. In this way, we design an adaptive strategy, which fully utilizes data from previous iterations and adjacent time subintervals, to progressively refine the spatial discretization of the coarse propagator. 
	
		Now, we introduce the details. We first collect the numerical solutions per $\delta M$ steps from calculating $\mathcal{F}\left(\hat{t}_n, \hat{t}_{n+1},U_n^{k-1}\right)$ to form the snapshot matrix $\mathbf{W}^k_{n,0}$.  Then, choosing parameter $\gamma_2\in(0, 1]$ and performing model reduction on $\mathbf{W}^k_{n,0}$ by \texttt{POD\_Mode} $\left(\mathbf{W}^k_{n,0}, \gamma_2, {\bf M}, \mathbf{\Phi}_h,  m^k_{n,0},  \mathbf{\Psi}^k_{n,0},  {\bf\widetilde{R}}^k_{n, 0}\right)$ as in Algorithm \ref{PODMode}, we obtain new POD modes $\mathbf{\Psi}^k_{n,0}$ on each time subinterval. For the time-dependent problems, the solution at a given time depends on that at earlier times. We consider to utilize the data ${\bf\widetilde{R}}^k_{n, 0}$ from the previous $p$ iterations and 
		 the neighboring time subintervals on the left to construct the POD subspace $\widetilde{V}^k_{n,{\mathrm{POD}}}$ on the time subinterval $\left[\hat{t}_n, \hat{t}_{n+1}\right]$ and obtain
	\begin{equation}\label{W^k_n}
		\mathbf{W}^k_{n}=\left[{\bf\widetilde{R}}^{k}_{n-m_l, 0}, {\bf\widetilde{R}}^{k}_{n-m_l+1, 0}, \ldots, {\bf\widetilde{R}}^{k}_{n, 0},{\bf\widetilde{R}}^{k-1}_{n, 0}, \ldots, {\bf\widetilde{R}}^{k-p}_{n, 0}\right],
	\end{equation}
	where $m_l$ denotes the number of time subintervals selected on the left side. Finally, choosing parameter $\gamma_3\in(0, 1]$ and performing model reduction on $\mathbf{W}^k_{n}$ by \texttt{POD\_Mode} $\left(\mathbf{W}^k_{n}, \gamma_3, {\bf M}, \mathbf{\Phi}_h,  m^k_{n},  \mathbf{\Psi}^k_{n},  {\bf\widetilde{R}}^k_{n}\right)$ as in Algorithm \ref{PODMode}, we obtain the POD subspace $\widetilde{V}^k_{n,{\mathrm{POD}}}$ spanned by POD modes $\mathbf{\Psi}^k_{n}$ on each subinterval.  
	
	We note that $U^{k}_n \notin \widetilde{V}^k_{n, \mathrm{POD}}$, which will introduce projection error when projecting the initial value onto the POD subspace during the sequential update of $U_n^k$. This part of the error cannot be removed by updating the POD subspaces. To remove the errors introduced by the Galerkin projection, we improve our method by adding $U^k_n$ to $\widetilde{V}^k_{n,\mathrm{POD}}$. Following is the details for how to add $U^k_n$ to $\widetilde{V}^k_{n,\mathrm{POD}}$.
	 
	The POD modes for the $n$-th time subinterval in the $k$-th iteration are 
	$$
	\mathbf{\Psi}^k_n=\left(\psi^k_{n, 1},  \psi^k_{n, 2},  \ldots,  \psi^k_{n, m^k_n}\right)=\mathbf{\Phi}_h \widetilde{\mathbf{R}}^k_{n}.
	$$
	Denote $\mathbb{P}_{\widetilde{V}^k_{n, \mathrm{POD}}}: V_h\to \widetilde{V}^k_{n, \mathrm{POD}}$  the Galerkin projection onto the POD subspace $\widetilde{V}^k_{n, \mathrm{POD}}$. Considering $U^{k}_n=\mathbf{\Phi}_h\mathbf{U}^{k}_{n}$, and denoting $\mathbf{\widetilde{d}}^{k}_{n}$  the coefficient column vector of $U^{k}_n-\mathbb{P}_{\widetilde{V}^k_{n, \mathrm{POD}}}\left(U^{k}_n\right)$ in the basis $\mathbf{\Phi}_h$, we obtain  
	\begin{equation}\label{eqaug}
		\mathbf{\widetilde{d}}_{n}^{k}=\mathbf{U}^{k}_{n}-\widetilde{\mathbf{R}}^k_{n}{\left(\widetilde{\mathbf{R}}^k_{n}\right)}^{\top}\mathbf{M} \mathbf{U}^{k}_{n}.
	\end{equation}
	Let $\mathbf{d}^{k}_{n}=\frac{\mathbf{\widetilde{d}}^{k}_{n}}{\sqrt{  \left( \mathbf{\widetilde{d}}^{k}_{n}\right)^{\top}  \mathbf{M} \mathbf{\widetilde{d}}^{k}_{n}}}$,  the new POD modes are
	\begin{equation}\label{eq12}
		\left(\psi^k_{n, 1},  \psi^k_{n, 2},  \ldots,  \psi^k_{n, {m^k_n+1}}\right)=\mathbf{\Phi}_h \left(\widetilde{\mathbf{R}}^k_{n}, \mathbf{d}^{k}_n\right).
	\end{equation}
	
	After adding $U_n^{k}$ to $\widetilde{V}^k_{n,\mathrm{POD}}$, we obtain the final subspace
	\begin{equation*}
		{V}^k_{n, \mathrm{POD}}  := \text{span}\{\psi^k_{n, 1},  \psi^k_{n, 2},  \ldots,  \psi^k_{n, {m^k_n+1}} \}.
	\end{equation*}
	Set $\widetilde{u}^k_{n,l}$  the POD approximation at $t=\hat{t}_n+l\cdot dT$. We have $\widetilde{u}^k_{n,l}=\sum_{i=1}^{m^k_n +1} {\beta}^k_{n,l, i} \psi^k_{n,  i}$, and define
	\begin{equation*}
		\widetilde{\mathbf{u}}^k_{n,l}=\left({\beta}^k_{ n,l, 1}, \ldots, {\beta}^k_{n,l, m^k_n+1} \right)^{\top}.
	\end{equation*}
	
	Denote $\mathbf{A}_{n,l}$ and $\mathbf{b}_{n,l}$  the corresponding stiff matrix and vector at $\hat{t}_n+l\cdot dT$, respectively. By setting $V_h$ in \eqref{eq5} to $V_{n, \mathrm{POD}}^{k}$, we obtain 
	\begin{equation*}
		\begin{split}
			\left(\sum_{i=1}^{m^k_n+1} \beta^k_{n,l, i} \psi^k_{n,i}-\sum_{i=1}^{m^k_n +1} \beta^k_{n,l-1, i} \psi^k_{n,i},  \psi^k_{n,j}\right) &+ d T\cdot \left(\mathcal{A}(\hat{t}_n+l\cdot dT, \sum_{i=1}^{m^k_n+1} \beta^k_{n,l, i} \psi^k_{n,i}), \psi^k_{n,j}\right)\\
			&=dT\cdot(f(\hat t_n+l\cdot dT), \psi^k_{n,j}).
		\end{split}
	\end{equation*}
	Combining with \eqref{eq8} and \eqref{eq12}, we can rewrite the equations above as 
	\begin{equation}\label{eq13}
		\begin{split}
			\left[\begin{array}{cc}
				{\left(\widetilde{\mathbf{R}}^k_{n}\right)}^{\top} \mathbf{A}_{n,l} \widetilde{\mathbf{R}}^k_{n} & {\left(\widetilde{\mathbf{R}}^k_{n}\right)}^{\top} \mathbf{A}_{n,l} \mathbf{d}^{k}_n \\
				\left(\mathbf{d}^{k}_n\right)^{\top} \mathbf{A}_{n,l} {\left(\widetilde{\mathbf{R}}^k_{n}\right)}^{\top} & \left(\mathbf{d}^{k}_n\right)^{\top} \mathbf{A}_{n,l} \mathbf{d}^{k}_n
			\end{array}\right] \widetilde{\mathbf{u}}^k_{n,l} &= \left[\begin{array}{c}
				{\left(\widetilde{\mathbf{R}}^k_{n}\right)}^{\top} \mathbf{b}_{n,l} \\
				\left(\mathbf{d}^{k}_n\right)^{\top} \mathbf{b}_{n,l}
			\end{array}\right] \\
			&\quad + \left[\begin{array}{cc}
				{\left(\widetilde{\mathbf{R}}^k_{n}\right)}^{\top} \mathbf{M} \widetilde{\mathbf{R}}^k_{n} & {\left(\widetilde{\mathbf{R}}^k_{n}\right)}^{\top} \mathbf{M} \mathbf{d}^{k}_n \\
				\left(\mathbf{d}^{k}_n\right)^{\top} \mathbf{M} {\left(\widetilde{\mathbf{R}}^k_{n}\right)}^{\top} & \left(\mathbf{d}^{k}_n\right)^{\top} \mathbf{M} \mathbf{d}^{k}_n
			\end{array}\right]  \widetilde{\mathbf{u}}^k_{n,{l-1}}.
		\end{split}
	\end{equation}
	
	We then obtain our adaptive parareal method as the following formula \eqref{formula4}:
	\begin{equation}\label{formula4}
		\begin{aligned}
			& \left\{
			\begin{aligned}
				U_{n+1}^0 &= \mathcal{G}_{0}\left(\hat{t}_n, \hat{t}_{n+1},  \mathbb{P}_{{V}^0_{\mathrm{POD}}}\left( U_n^0\right)\right),\quad  U_0^0 = \hat{u}_0,  \\
				U_{n+1}^{k} &= {\mathcal{G}}_{k, n}\left(\hat{t}_n, \hat{t}_{n+1}, \mathbb{P}_{{V}^k_{n, \mathrm{POD}}}\left( U^{k}_n\right)\right) +\mathcal{F}\left(\hat{t}_n, \hat{t}_{n+1},  U_n^{k-1}\right) \\
				&\quad -{\mathcal{G}}_{k, n}\left(\hat{t}_n, \hat{t}_{n+1},  \mathbb{P}_{{V}^k_{n, \mathrm{POD}}}\left( U^{k-1}_n\right)\right) , \quad U_0^{k} = \hat{u}_0, \quad k=1,2 \ldots.
			\end{aligned}\right. \\
		\end{aligned}
	\end{equation}
	where $\mathcal{G}_{k,n}$ is  the propagator whose spatial discretization is   introduced above. 
	
	We summarize the general framework of the adaptive parareal method based on model order reduction as in Algorithm \ref{FEM-APODalg}.
	\\
	\begin{breakablealgorithm}
		\flushleft	\caption{General framework of the adaptive parareal method based on model order reduction}\label{FEM-APODalg}
		\begin{algorithmic}[1]
			\STATE Given the mesh $\mathcal{T}_h$,  time steps $\delta t$,  $dT$,  parameters $T$, $T_0$, $\Delta T$, $\delta M$, $\gamma_1$, $\gamma_2$, $\gamma_3$, $p, m_l$;
			\STATE Construct $\mathcal{F}$; 
			\STATE Discretize \eqref{semi-dis-eq} on $[0, T_0]$ by $\mathcal{F}$ and obtain the snapshot matrix ${\bf W}_{h}$;
			\STATE Construct POD modes $\mathbf{\Psi}_h$ by \texttt{POD\_Mode} $\left(\mathbf{W}_h, \gamma_1, {\bf M}, \mathbf{\Phi}_h,  m,  \mathbf{\Psi}_h,  \widetilde{\mathbf{R}}\right)$ and obtain $\mathcal{G}_0$;
			\STATE Discretize \eqref{semi-dis-eq} on $[T_0, T]$ by $\mathcal{G}_0$ to obtain initial values $U_{n}^0$ for each time subinterval 
			$$
			U_0^0 = \hat{u}_0,  \quad U_{n+1}^0 = \mathcal{G}_{0}\left(\hat{t}_n, \hat{t}_{n+1},  \mathbb{P}_{{V}^0_{\mathrm{POD}}}\left(U_n^0\right)\right),  n=0,  \ldots,  N-1;
			$$
			\STATE $k=1$;
			\WHILE {not converged}
			\STATE Discretize \eqref{semi-dis-eq} on $\left[\hat{t}_n, \hat{t}_{n+1}\right]$ by $\mathcal{F}$ in parallel to obtain $\mathcal{F}\left(\hat{t}_n, \hat{t}_{n+1},  U_n^{k-1}\right)$ and the snapshot matrix $\mathbf{W}^k_{n,0}$, $n=k-1,  \ldots,  N-1;$ 
			
			\STATE Construct ${\bf\widetilde{R}}^k_{n, 0}$ by \texttt{POD\_Mode} $\left(\mathbf{W}^k_{n,0}, \gamma_2, {\bf M}, \mathbf{\Phi}_h,  m^k_{n,0},  \mathbf{\Psi}^k_{n,0},  {\bf\widetilde{R}}^k_{n, 0}\right)$ and obtain
			$$
			\mathbf{W}^k_{n}=\left[{\bf\widetilde{R}}^{k}_{n-m_l, 0}, {\bf\widetilde{R}}^{k}_{n-m_l+1, 0}, \ldots, {\bf\widetilde{R}}^{k}_{n, 0},{\bf\widetilde{R}}^{k-1}_{n, 0}, \ldots, {\bf\widetilde{R}}^{k-p}_{n, 0}\right],
			$$
			in parallel, $n=0,\ldots,N-1$;
			
			\STATE Construct ${\bf\widetilde{R}}^k_{n}$ by \texttt{POD\_Mode} $\left(\mathbf{W}^k_{n}, \gamma_3, {\bf M}, \mathbf{\Phi}_h,  m^k_{n},  \mathbf{\Psi}^k_{n},  {\bf\widetilde{R}}^k_{n}\right)$ in parallel, $n=k-1,  \ldots,  N-1;$ 
			
			\STATE $n=k$;
			\WHILE {$n<N$} 
			\STATE Add $U_{n}^{k}$ to the POD modes $\mathbf{\Psi}^k_{n}$ to get the final subspace ${V}^k_{n, \mathrm{POD}}$ and then obtain $\mathcal{G}_{k, n}$;
			
			\STATE Discretize \eqref{semi-dis-eq} on $\left[\hat{t}_n, \hat{t}_{n+1}\right]$ by $\mathcal{G}_{k, n}$ to obtain ${\mathcal{G}}_{k, n}\left(\hat{t}_n, \hat{t}_{n+1}, \mathbb{P}_{{V}^k_{n, \mathrm{POD}}}\left(U^{k-1}_n\right)\right)$ and  ${\mathcal{G}}_{k, n}\left(\hat{t}_n, \hat{t}_{n+1}, \mathbb{P}_{{V}^k_{n, \mathrm{POD}}}\left( U^{k}_n\right)\right)$ in parallel;
			\STATE Update $U_{n+1}^{k}$ by 
			\begin{equation*}
				\begin{split}
					U_{n+1}^{k}&={\mathcal{G}}_{k, n}\left(\hat{t}_n, \hat{t}_{n+1}, \mathbb{P}_{{V}^k_{n, \mathrm{POD}}}\left( U^{k}_n\right)\right) +\mathcal{F}\left(\hat{t}_n, \hat{t}_{n+1},  U_n^{k-1}\right)\\
					&\quad -{\mathcal{G}}_{k, n}\left(\hat{t}_n, \hat{t}_{n+1},  \mathbb{P}_{{V}^k_{n, \mathrm{POD}}}\left(U^{k-1}_n\right)\right); 
				\end{split}
			\end{equation*}
			\STATE $n=n+1$;
			\ENDWHILE
			\STATE $k=k+1$.
			\ENDWHILE
		\end{algorithmic}
	\end{breakablealgorithm}
	\quad\\\par
	
	\begin{remark}
		Compared with the plain parareal method, where $\mathcal{G}\left(\hat{t}_n, \hat{t}_{n+1}, U^{k-1}_n\right)$ is obtained from the previous iteration, the addition of $U_n^{k}$ to $\widetilde{V}^k_{n,\mathrm{POD}}$ and the calculation of ${\mathcal{G}}_{k, n}\left(\hat{t}_n, \hat{t}_{n+1},  \mathbb{P}_{{V}^k_{n, \mathrm{POD}}}\left( U^{k-1}_n\right)\right)$ in our method must be done sequentially. Nonetheless, considering that the number of POD modes is significantly smaller than the dimension of the finite element space, and only one element is added, the additional computational cost remains little compared to that of the fine propagator.
	\end{remark}

	\subsection{Speedup}
	\label{sec: Complex}
	
	It is known that the goal of a parareal method is to obtain high speedup. In this subsection,  we analyze the speedup of our adaptive parareal method Algorithm \ref{FEM-APODalg}.
	
	We use the notation ${\operatorname {cost}}_{\operatorname {par}}\left([0,T]\right)$ to denote the time cost   of our adaptive parareal method Algorithm \ref{FEM-APODalg}, and use  notation ${\operatorname {cost}}_{\operatorname {seq}}\left([0,T]\right)$ to denote the time cost for propagating the equation \eqref{eq0}  by the fine propagator sequentially.

	The speedup is then defined as follows \cite {fischer2005parareal}
	\begin{align*}
		{\operatorname {speedup}}\left([0,T]\right)=\frac{{\operatorname {cost}}_{\operatorname {seq}}\left([0,T]\right)}{ {\operatorname {cost}}_{\operatorname {par}}\left([0,T]\right)}.
	\end{align*}

	We denote $\tau_{\mathcal{F}}$ and $\tau_{\mathcal{G}}$ the time cost for propagating the equation \eqref{eq0}  one  step  by using the fine and coarse propagators respectively. Then we can deduce that the time cost ${\operatorname {cost}}_{\operatorname {seq}}([0,T])$ to carry the fine propagator serially is $\left(N\frac{\Delta T}{\delta t}+\frac{T_0}{\delta t}\right)\tau_{\mathcal{F}}$. That is, 
	\begin{align*}
		{\operatorname {cost}}_{\operatorname {seq}}\left([0,T]\right) = \left(N\frac{\Delta T}{\delta t}+\frac{T_0}{\delta t}\right)\tau_{\mathcal{F}}.
	\end{align*}

	We now analyze the cost for our adaptive parareal method Algorithm \ref{FEM-APODalg}. 
	The time cost for carrying the coarse propagator at the $0$-th iteration is $N\frac{\Delta T}{dT}\tau_{\mathcal{G}}$.   For each parareal iteration, the time cost for calculating $\mathcal{F}(\hat{t}_n,  \hat{t}_{n+1}, U_n^{k-1})$ is $\frac{\Delta T}{\delta t}\tau_{\mathcal{F}}$. For the calculation of ${\mathcal{G}}_{k, n}\left(\hat{t}_n,  \hat{t}_{n+1},  \mathbb{P}_{{V}^k_{n, \mathrm{POD}}}\left( U^{k-1}_n\right)\right)$ and ${\mathcal{G}}_{k, n}\left(\hat{t}_n,  \hat{t}_{n+1},  \mathbb{P}_{{V}^k_{n, \mathrm{POD}}}\left( U^{k}_n\right)\right)$, since they are independent of each other, they can be carried out in parallel essentially. As we mentioned before, according to the finite termination property of the parareal method \cite{baffico2002parallel,Bal2002}, the approximation $U_n^k$ converges exactly to the sequential fine solution $U_n$ for all time steps up to $n \le k$. This indicates that $U^{k}_k=U_k$. Justified by this property, the coarse propagator does not have to execute from $T_0$ to $T$, but only needs to cover $[T_0+k\Delta T,T]$. Therefore, the time cost is $(N-k)\frac{\Delta T}{dT}\tau_{\mathcal{G}}$. The cost of obtaining or updating the coarse propagator, which is incurred only once pre iteration, consists of constructing POD subspace (denoted $T_C$) and forming the reduced system \eqref{eq8} (denoted $T_U$).  The time cost for adding $U^k_n$ to the original POD subspace $\widetilde{V}^k_{n,\mathrm{POD}}$, which has to be carried in each time subinterval sequentially, is denoted as $T_A$. Denote $k_{\operatorname{max}}$ the number of iterations required for convergence. 
	
	Therefore, by the analysis above, we have that 
	\begin{equation*}
		\begin{aligned}
			{\operatorname {cost}}_{\operatorname {par}}\left([0,T]\right) &= \frac{T_0}{\delta t}\tau_{\mathcal{F}}+T_C  + T_U+ N\frac{\Delta T}{dT}\tau_{\mathcal{G}}\\
			&\quad +\sum_{k=1}^{k_{\operatorname{max}}} \left(\frac{\Delta T}{\delta t}\tau_{\mathcal{F}}+T_C  + T_U+\left(N-k\right) \left(T_A+ \frac{\Delta T}{dT}\tau_{\mathcal{G}}\right) \right) \\
			&= \frac{T_0}{\delta t}\tau_{\mathcal{F}}+T_C +T_U+ N\frac{\Delta T}{dT}\tau_{\mathcal{G}}\\
			&\quad +k_{\operatorname{max}}\left(\frac{\Delta T}{\delta t}\tau_{\mathcal{F}}+T_C +T_U+\left(N- \frac{k_{\operatorname{max}}+1}{2}\right) \left(T_A+ \frac{\Delta T}{dT}\tau_{\mathcal{G}}\right) \right). 
		\end{aligned}
	\end{equation*}

	In further, we analyze the cost for each part. 
	
	We first analyze $\tau_{\mathcal{F}}$.  Due to the sparsity of the matrices in the linear system \eqref{algform}, the computational complexity for  constructing \eqref{algform} is $C_{f_1} N_g$ \cite{dai2024augmented}, with the constant factor $C_{f_1}$ depending on the type of finite element, which is usually of magnitude tens or hundreds. For solving the linear system \eqref{algform},   although there exist some linear solvers which can achieve approximately linear complexity, their constant factors are highly problem-dependent and typically very large, often being of magnitude hundreds or thousands. 	
	Therefore, $\tau_{\mathcal{F}} \approx C_{f_1} N_g + C_{f_2} N_g$, with $C_{f_1} $ being of magnitude tens or hundreds and $C_{f_2}$ being  of magnitude hundreds or thousands.
	
	Furthermore, for the cases where $\mathcal{A}\left(t,u\right)$ can be separable in time and space,  we only need to construct the system \eqref{algform} at the first time step. Therefore, in this case, $\tau_{\mathcal{F}} \approx C_{f_2} N_g$.  
	
	For simplicity, we write $\tau_{\mathcal{F}} \approx C_{f} N_g$ with  $C_{f} = C_{f_1} + C_{f_2}$ or $C_{f} =C_{f_2}$, a constant being of magnitude hundreds or thousands. 
	
	We then analyze $\tau_{\mathcal{G}}$. 
	Denote $m_{\operatorname{max}}$ the maximum dimension of the POD subspaces across all time subintervals. For solving the reduced system \eqref{eq8} or \eqref{formula4} ($k\geq 1$), whose coefficient matrix is dense, we know that the computational complexity   is  $ C_{p_1} m_{\operatorname{max}}^3$, where $C_{p_1}$ can be less than $1$.   Hence, $\tau_{\mathcal{G}} \approx  C_{p_1} m_{\operatorname{max}}^3$.

	For the time cost $T_C$, we distinguish two cases. For the $0$-th iteration, the computational complexity of performing model order reduction operation on $\mathbf{W}_h$ is about $C_{p_2} n_s^2 N_g$ \cite{golub2013matrix}, where $n_s = \left\lfloor\frac{\Delta T}{\delta t \cdot \delta M}\right\rfloor+1$, and the constant factor $C_{p_2}$ can be less than $10$. So $T_C \approx C_{p_2} n_s^2 N_g $.  For the $k$-th iteration with $k\geq 1$, the operation of model order reduction has to be carried out twice. The first is performed on $\mathbf{W}^k_{n,0}$, and the second is on $\mathbf{W}^{k}_{n}$. Let   $n_{\operatorname{max}}$ be the maximum column dimension of the matrices involved in the POD across all time subintervals. Then, we have $T_C \approx  C_{p_2} \left(n_s^2+  n_{\operatorname{max}}^2 \right) N_g$ .

    For the time cost $T_U$, once we obtain the POD subspace, the time cost for obtaining  the reduced system \eqref{eq8} is $C_{p_3}  m^2_{\operatorname{max}} N_g$,  with constant factor  $C_{p_3}$ usually being  of magnitude tens.  So $T_U \approx C_{p_3}  m^2_{\operatorname{max}} N_g$.

	Finally, we consider  the time cost $T_A$.  
	 Similar to the cost for obtaining the system \eqref{eq8}, we can see that the cost for calculating \eqref{eqaug} is  about  $C_{p_3}  m_{\operatorname{max}}  N_g$. Therefore, 
	$T_A \approx C_{p_3} m_{\operatorname{max}}  N_g$.

	Therefore,  for the case   where  $\mathcal{A}\left(t,u\right)$ can be separable in time and space, we have 
	\begin{equation*}
		\begin{aligned}
			{\text{cost}}_{\text{par}}\left([0,T]\right) &= \frac{T_0}{\delta t}\tau_{\mathcal{F}}+T_C +T_U+ N\frac{\Delta T}{dT}\tau_{\mathcal{G}}\\
			&\quad +k_{\text{max}}\left(\frac{\Delta T}{\delta t}\tau_{\mathcal{F}}+T_C +T_U+\left(N- \frac{k_{\text{max}}+1}{2}\right) \left(T_A+ \frac{\Delta T}{dT}\tau_{\mathcal{G}}\right) \right) \\
			&= \frac{T_0}{\delta t} C_f N_g + C_{p_2} n_s^2 N_g + C_{p_3} m^2_{\operatorname{max}} N_g + C_{p_1} N\frac{\Delta T}{dT} m_{\text{max}}^3 \\
			&\quad +k_{\text{max}}\left(\frac{\Delta T}{\delta t} C_f N_g + C_{p_2} \left(n_s^2+ n_{\text{max}}^2 \right) N_g + C_{p_3} m^2_{\text{max}} N_g\right) \\
			&\quad +k_{\text{max}}\left(N- \frac{k_{\text{max}}+1}{2}\right) \left(C_{p_3} m_{\text{max}} N_g + \frac{\Delta T}{dT} C_{p_1} m_{\text{max}}^3\right) \\
			&= \Bigl( \frac{T_0}{\delta t} C_f + C_{p_2} n_s^2 + C_{p_3} m^2_{\text{max}} \\
			&\qquad + k_{\text{max}}\Bigl(\frac{\Delta T}{\delta t} C_f + C_{p_2}\left(n_s^2 + n_{\text{max}}^2\right) + C_{p_3} m^2_{\text{max}} \\
			&\qquad\qquad + C_{p_3} m_{\text{max}} \left(N- \frac{k_{\text{max}}+1}{2}\right)\Bigr)\Bigr) N_g \\
			&\quad + C_{p_1} \frac{\Delta T}{dT} \left( N + k_{\text{max}} \left(N- \frac{k_{\text{max}}+1}{2}\right) \right) m_{\text{max}}^3,
		\end{aligned}
	\end{equation*}
	and 
	\begin{equation*}
		\begin{aligned}
		{\text{cost}}_{\text{seq}}\left([0,T]\right) = 	\left(N\frac{\Delta T}{\delta t}+\frac{T_0}{\delta t}\right)\tau_{\mathcal{F}} = \left(N\frac{\Delta T}{\delta t}+\frac{T_0}{\delta t}\right)C_f N_g.  
      \end{aligned}
\end{equation*}
	
	From the former analysis, we know that $C_f$ is usually of magnitude hundreds or thousands, while $C_{p_1}$ is usually less than $1$, $C_{p_2}$ and $C_{p_3}$ are usually of magnitude tens. Besides, from our numerical tests, we know that  $m_{\text{max}}$, $n_s$, and $n_{\text{max}}$   generally do not exceed an order of magnitude of a hundred.  Note that $\Delta T \gg \delta t$ and $N_g \gg m^3_{\text{max}}$,  we have that $C_{p_2}\left(n_s^2+n_{\text{max}}^2\right) + C_{p_3} m_{\text{max}}^2 \le \frac{\Delta T}{\delta t}C_f$. For the convenience of our analysis for the speedup, we only  keep the dominant term in the ${\text{cost}}_{\text{par}}([0,T])$. That is, we have 
		\begin{equation*}
		\begin{aligned}
		{\text{cost}}_{\text{par}}([0,T])	&\approx
		\left( \frac{T_0}{\delta t} C_f +  k_{\text{max}} \left(\frac{\Delta T}{\delta t} C_f  + C_{p_3} m_{\text{max}} \left(N-  \frac{k_{\text{max}}+1}{2}\right) \right) \right)   N_g \\
		&\quad + C_{p_1} \frac{\Delta T}{dT}  k_{\text{max}} m_{\text{max}}^3N .
		\end{aligned}
		\end{equation*}
	
	Assume that we have $N$ groups of processors, with each group being assigned an equal number of processors per time subinterval to perform the spatial parallelization. If we ignore  the communication time and only keep the dominant term in the ${\text{cost}}_{\text{par}}([0,T])$, we have that
	\begin{equation}\label{speedup0}
		\begin{aligned}
		&{\text{speedup}}\left([0,T]\right) = \frac{{\text{cost}}_{\text{seq}}([0,T])}{{\operatorname {cost}}_{\text{par}}([0,T])} \\
		\approx~& 
		\frac{\left(N\frac{\Delta T}{\delta t}+\frac{T_0}{\delta t}\right)C_f N_g}{\left( \frac{T_0}{\delta t} C_f +  k_{\text{max}} \left(\frac{\Delta T}{\delta t} C_f  + C_{p_3} m_{\text{max}} \left(N-  \frac{k_{\text{max}}+1}{2}\right) \right) \right)   N_g +C_{p_1} \frac{\Delta T}{dT}  k_{\text{max}} m_{\text{max}}^3N}\\
		=~& 
		\frac{\left(N\frac{\Delta T}{\delta t}+\frac{T_0}{\delta t}\right)}{\left( \frac{T_0}{\delta t}  +  k_{\text{max}} \left(\frac{\Delta T}{\delta t}    + \frac{C_{p_3} m_{\text{max}}}{C_f} \left(N-  \frac{k_{\text{max}}+1}{2}\right) \right) \right)   + \frac{C_{p_1} \frac{\Delta T}{dT}   m_{\text{max}}^3}{C_f} \frac{k_{\text{max}} N}{N_g}}\\
		\approx~& 
		\frac{ N\frac{\Delta T}{\delta t}}{  k_{\text{max}} \left(\frac{\Delta T}{\delta t}    + \frac{C_{p_3} m_{\text{max}}}{C_f}  N \right)    + \frac{C_{p_1} m_{\text{max}}^3}{C_f} \frac{k_{\text{max}} N}{N_g}  \frac{\Delta T}{dT} }\\
		=~& 
		\frac{ N}{k_{\text{max}} + \frac{C_{p_3} m_{\text{max}}}{C_f} \frac{\delta t}{\Delta T} k_{\text{max}} N  + \frac{C_{p_1} m_{\text{max}}^3}{C_f} \frac{k_{\text{max}} N}{N_g} \frac{\delta t}{dT} }.\\
      \end{aligned}
\end{equation}

From our previous analysis, we have that $\frac{C_{p_1} m_{\text{max}}^3}{C_f}$ is controllable. For the case of $N_g$ being far larger than $k_{\text{max}} N$, the last term in \eqref{speedup0} is negligible compared with the other two terms. Then, we have 
	\begin{equation}\label{speedup}
\begin{aligned}
	&~ {\text{speedup}}\left([0,T]\right) 	\approx 
	\frac{ N}{k_{\text{max}} + \frac{C_{p_3} m_{\text{max}}}{C_f} \frac{\delta t}{\Delta T} k_{\text{max}} N} = \frac{ N}{\left(1+\frac{C_{p_3} m_{\text{max}}}{C_f} \frac{\delta t}{\Delta T} N \right)k_{\text{max}}}. 
  \end{aligned}
\end{equation}
 
\begin{remark}
  From \eqref{speedup}, we can see that  when high accuracy approximations are required, the performance of our method will be better. In fact, high accuracy approximation means large $N_g$ and small $\delta t$. Note that the constant $C_f$ is highly problem-dependent. Large $N_g$ means solving large scale linear systems, which means large $C_f$. Therefore, $\frac{C_{p_3} m_{\text{max}}}{C_f }$ will be far smaller than $1$. Small $\delta t$ means $\frac{\delta t}{\Delta T} \ll 1$ for given $\Delta T$. Therefore, for given $N$, if the number of parareal iterations $k_{\text{max}}$ required for convergence over the $N$ time subintervals does not increase quickly as the accuracy required increases, we will expect better speedup when high accuracy approximations are required. 
 
  Besides, from \eqref{speedup}, we can also see that, for fixed $\Delta T$ and $\delta t$ with $\Delta T \gg \delta t$,  if the  number of parareal iterations  $k_{\text{max}}$ required for convergence over $N$ time subintervals grows slower than linearly with $N$,  then the speedup increases with $N$, until some maximum speedup is obtained. For fixed $\Delta T$, a larger $N$ means a larger $T$. Therefore,   if the number of parareal iterations required for convergence by our method grows slower than linearly with $N$,  we can expect very large speedup for long time simulation, if we have enough number of processors. 
\end{remark}

	\section{Numerical experiments}
	\label{sec: examples}
	We consider the following time-dependent  partial differential equations with periodic boundary conditions
	\begin{equation}\label{eq1}
	\left\{\begin{array}{l}
	u_t-\epsilon \Delta u+\mathbf{B}(x,  y,  z,  t) \cdot \nabla u+c(x,  y,  z,  t) u=f(x,  y,  z,  t),  \quad \text{in}\ \Omega \times(0,  T],  \\
	u(x,  y,  z,  0)=h(x,  y,  z), \quad \text{in}\ \Omega \\
	u(x+L,  y,  z,  t)=u(x,  y+L,  z,  t)=u(x,  y,  z+L,  t) \\
	=u(x,  y,  z,  t),\quad  \text{on}\ \partial \Omega\times(0,  T], 
	\end{array}\right.
	\end{equation}
	where $\Omega=[0,  L]^3$, $f \in L^2\left(0,  T ; L^2(\Omega)\right)$, and $c \in C\left( 0,T; L^{\infty}\left(\Omega\right)\right)$. Here $\mathbf{B} \in C\left(0,  T;W^{1,  \infty}(\Omega)^3\right)$, $h\in L^2(\Omega) $, and $\epsilon$ is a positive constant. Generally, the smaller the $\epsilon$, the more challenging the simulation of this model becomes \cite{bal2005convergence,gander2007analysis}.
	
	In our numerical simulation, we choose two typical fluid advection fields in three dimensional space to show the performance of our methods. The relative error of approximation $U^k_n$ at $t=\hat{t}_n$ obtained by our parareal method is calculated by
	
	\begin{equation}\label{error}
	\mathrm{Error}=\frac{||U^k_n-U_n||_{L^2}}{||U_n||_{L^2}},
	\end{equation}
	where the reference solution  $U_n$ is that obtained by propagating \eqref{eq1}  using the fine propagator with exact initial values over $[0, T]$. 
	
	In the following discussions, we denote our adaptive parareal method Algorithm \ref{FEM-APODalg} as "AdapParareal". For the method "AdapParareal", we test the following $4$ cases respectively: $(m_l, p)=(0,0), (1,0), (0,1), (1,1)$. For the sake of simplicity, we denote ${\mathcal{G}}_{k, n}\left(\hat{t}_n, \hat{t}_{n+1}, \mathbb{P}_{{V}^k_{n, \mathrm{POD}}}\left(U^{k}_n\right)\right)$ as $\mathcal{G}_k(U^{k}_n)$ in this section.
	
	 For comparison, we will also show the performance of the plain parareal method. For fairness of the comparison, we set the parameters in the 	plain parareal 	method to be the same as those used in the method "AdapParareal" except for the spatial discretization of the coarse propagator. Here, we consider two types of the coarse spatial discretization. One is using the POD method to do the coarse spatial discretization with  the POD subspace being constructed from snapshots obtained from the simulation of the fine propagator over a time window $[0, T_0]$  in the $0$-th iteration. We denote the corresponding parareal method as "Parareal-I". The other one is using the finite element method with coarser mesh size than that used for fine propagator to do the coarse spatial discretization. We denote the corresponding parareal method as "Parareal-II".

	Our numerical experiments are carried out on the high performance computers of State Key Laboratory of Mathematical Sciences, Chinese Academy of Sciences, and our code is based on the toolbox PHG \cite{PHG}. In all numerical experiments, we use  $144$ processors to do the spatial parallelism. It should be noted that, due to the limitation of computational resources, the steps $8-10$ and the step $14$ in Algorithm \ref{FEM-APODalg} which should be executed in essentially in parallel across time subintervals are instead carried out sequentially, one time subinterval at a time. 
	To get a rough estimate for the speedup, in our  numerical experiments, we will record the wall time for each part of the algorithm, from which we can obtain an estimation for the speedup based on our previous analysis for the speedup in Section~\ref{sec: Complex}.

		Denote $T_{\mathcal{F}}$ and $T_{\mathcal{G}}$   the computational time of the fine propagator and the coarse propagator within each time subinterval per iteration respectively.	 Recall the notation in Section~\ref{sec: Complex}, we have 
		\begin{align*}
		T_{\mathcal{F}}=\frac{\Delta T}{\delta t}\tau_{\mathcal{F}},\quad T_{\mathcal{G}}=\frac{\Delta T}{dT}\tau_{\mathcal{G}}.
		\end{align*}
		The total time cost for the method per iteration is denoted as $T_{\text{total}}$.
	
	\subsection{Kolmogorov flow}\label{Kolmogorov flow}
	We consider the advection-diffusion equation \eqref{eq1}, where the advection field is defined by the Kolmogorov flow \cite{borue1996numerical, obukhov1983kolmogorov}. The components of the equation are given as follows:      
	\begin{equation}\label{koleq}
	\begin{aligned}
	\mathbf{B}(x,  y,  z,  t)&=(\cos (y),  \cos (z),  \cos (x))+(\sin (z),  \sin (x),  \sin (y)) \cos (t),  \\
	f(x,  y,  z,  t) &= -\cos (y)-\sin (z)\cos (t),  \\
	c(x,  y,  z,  t) &= 0, 
	~h(x,  y,  z) =0, ~L= 2 \pi,~T=1005.0.
	\end{aligned}
	\end{equation}

	In this example, $\mathbf{B}$ and $f$ are separable in time and space, so the computation for constructing the discretized system can be simplified. We refer to \cite{dai2024augmented} for more details. We use piecewise linear functions as the basis functions for the finite element method. We first divide $\Omega$ into $6$ tetrahedrons to create the initial grid and refine the initial mesh $21$ times uniformly using bisection to obtain the final mesh. The number of degrees of freedom is $2097152$. For the choice of $T_0$, we have done some tests and found  that our method "AdapParareal" is insensitive to it. In our tests, we simply set $T_0 = 5$. We fix the parameters $\delta t=1\times 10^{-2}$, $dT=0.5$, $\Delta T=5.0$, $\delta M=5$, $\gamma_1=\gamma_2=1.0-5.0\times10^{-6}$ and $\gamma_3=1.0-2.0\times10^{-8}$. The number of time subintervals is therefore $200$. For all figures in this article, the $x$-axis represents time $t$ and $k$ represents the index of iteration, and the $y$-axis represents the relative error of $U^k_n$ defined by \eqref{error}.

\begin{figure}[H]
		\centering
		\subfloat[$\epsilon=0.5.$]{\includegraphics[width=0.5\textwidth]{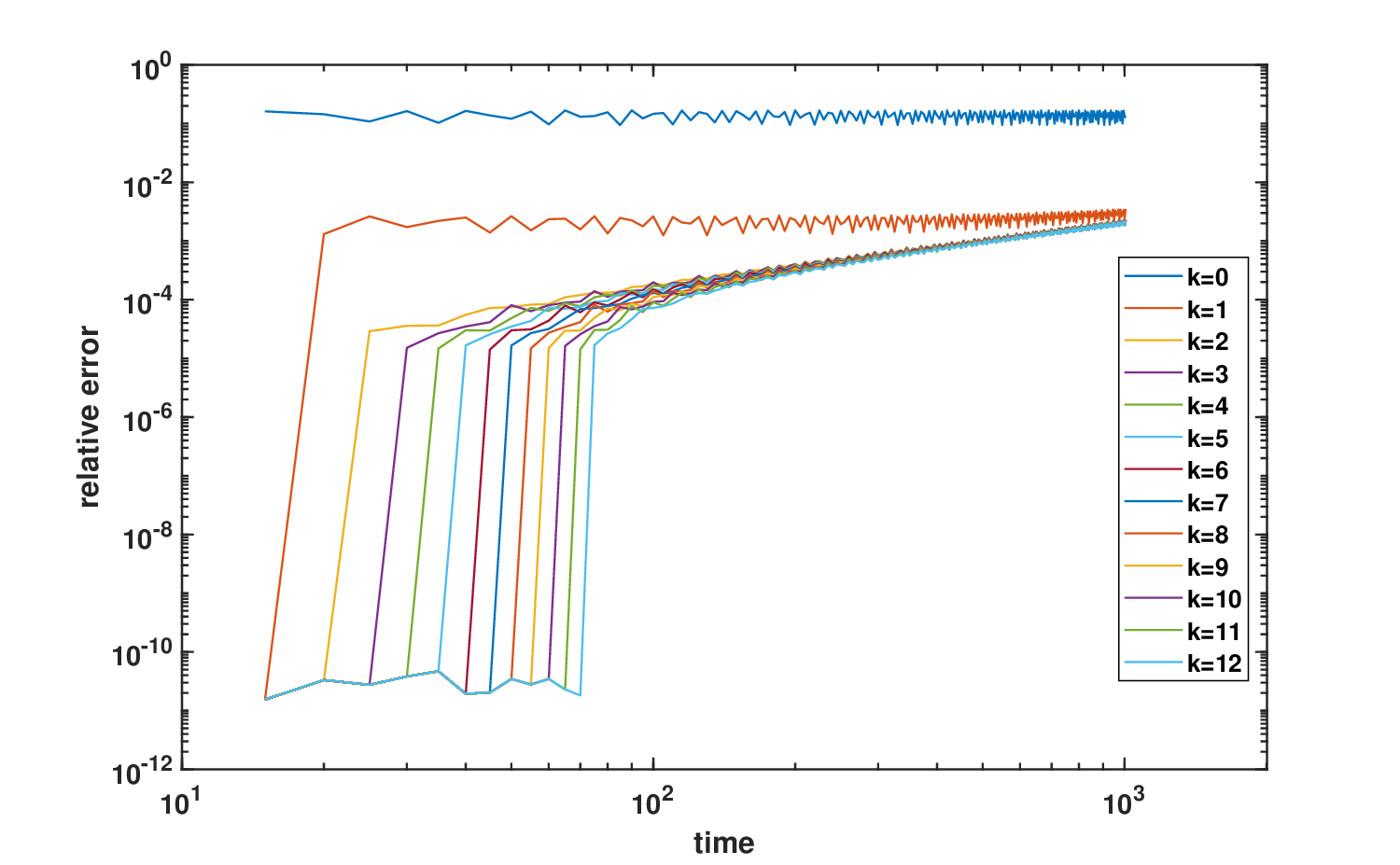}}
		\subfloat[$\epsilon=0.1.$]{\includegraphics[width=0.5\textwidth]{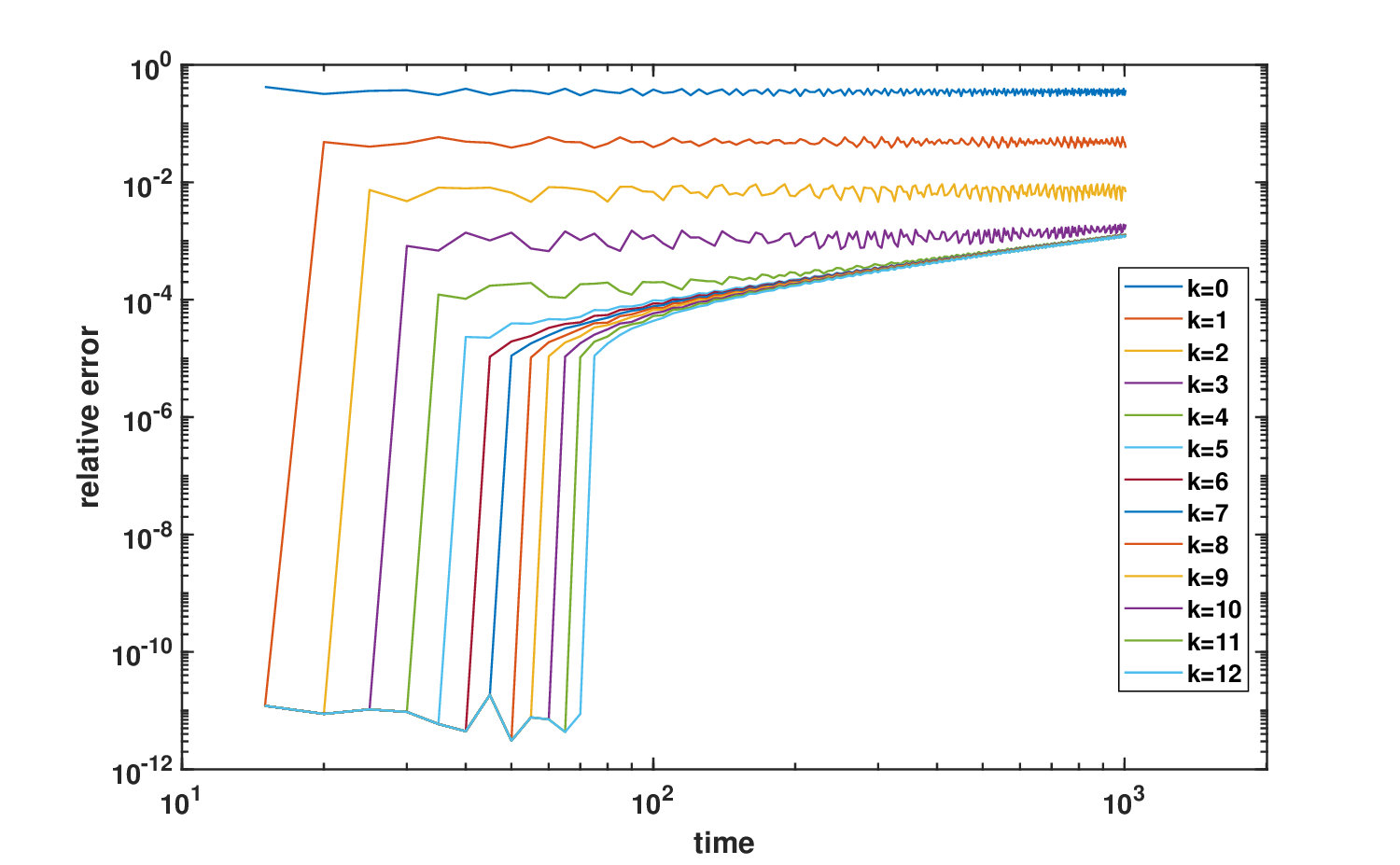}}
		\caption{The evolution curves of the relative error of $U^k_n$ obtained by the method "Parareal-I" in each parareal iteration for solving the Kolmogorov flow with $\epsilon=0.5,0.1$}
		\label{KOLFEM-POD}
	\end{figure}

	For comparison, we first show the performance of  the method "Parareal-I".  We  test  two cases with $\epsilon = 0.5$ and $\epsilon = 0.1$. Considering that the plain parareal method is sensitive to the accuracy of the coarse propagator, for method "Parareal-I", we use a relatively large $T_0$ to obtain the POD subspace. Here, we  set $T_0=10.0$. For fairness of the comparison, we set the other parameters to be the same as those used in the method "AdapParareal". The numerical results  are shown in Fig. \ref{KOLFEM-POD}. The dimensions of POD subspaces for the coarse propagator in the $0-$th iteration for these two cases are $24$ and $38$, respectively. For problem \eqref{eq1}-\eqref{koleq}, we know that the smaller the $\epsilon$, the more rapidly the solution varies. This is the reason for the difference between number of POD modes obtained.  From the two sub-figures in Fig. \ref{KOLFEM-POD}, we can see that the relative error can reach $O(10^{-3})$ after several iterations. It shows the effectiveness of method "Parareal-I". However, we also see that the accuracy cannot be improved any more as the number of iterations increases. Considering that the model will be more difficult to simulate as $\epsilon$ decreases, we cannot expect better results for smaller $\epsilon$ and hence skip the showing of the numerical results for cases of smaller $\epsilon$. 
	
	\FloatBarrier
	
\begin{figure}[H]
		\centering
		\subfloat[$\epsilon=0.5.$]{\includegraphics[width=0.45\textwidth]{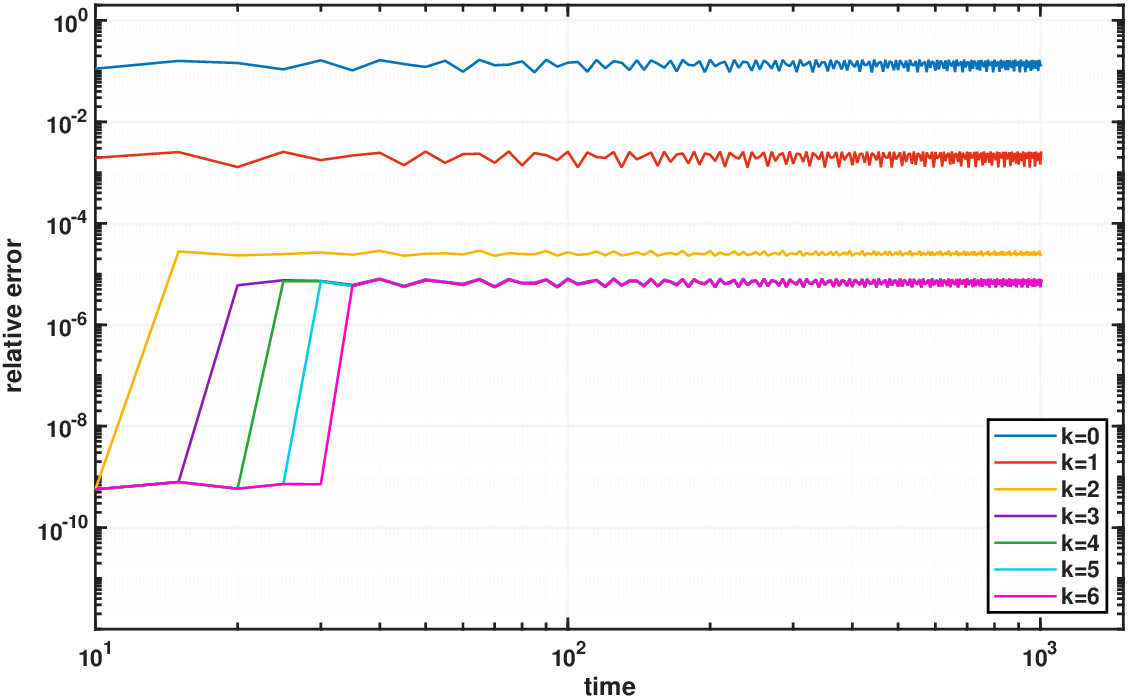}}\hspace{1em}
		\subfloat[$\epsilon=0.1.$]{\includegraphics[width=0.45\textwidth]{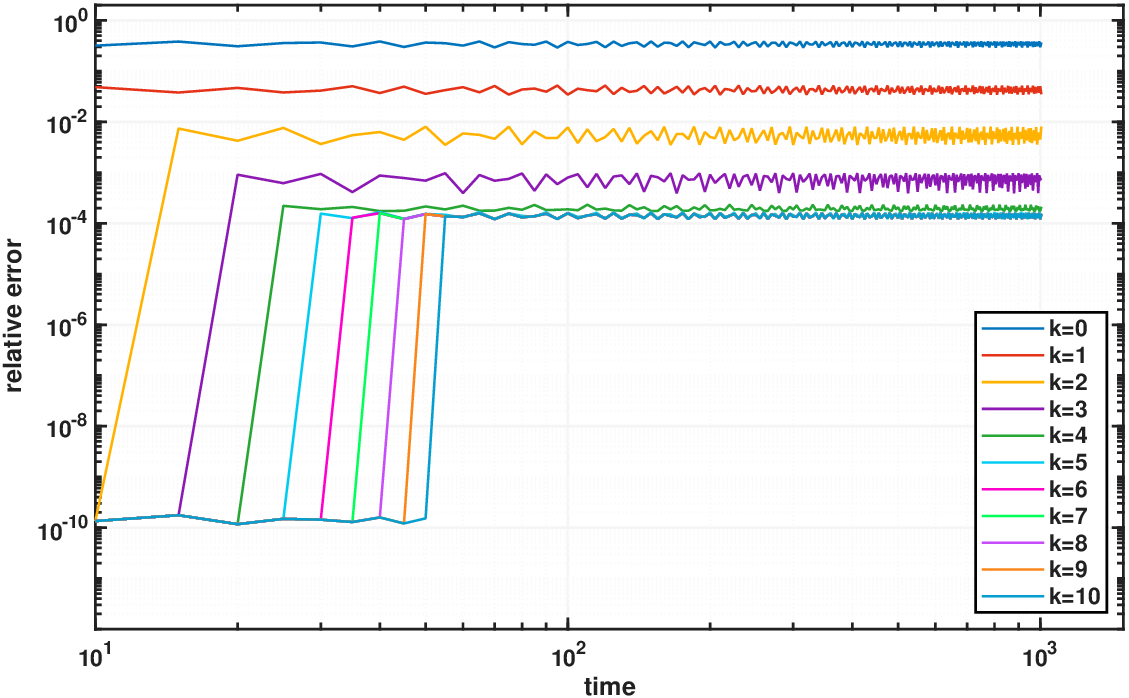}}\hspace{1em}
		\caption{The evolution curves of the relative error of $U^k_n$ obtained by the method "Parareal-II" in each parareal iteration for solving the Kolmogorov flow with different $\epsilon$.}
		\label{KOLFEM-FEM}
	\end{figure}

	We then show the performance of the method "Parareal-II", where finite element method with coarser mesh size is used for the spatial discretization in the coarse propagator.   We refine the initial mesh $18$ times uniformly using bisection to obtain the final coarse mesh and the number of degrees of freedom is $262144$. We test the  cases of $\epsilon = 0.5$ and $\epsilon = 0.1$.  The numerical results are shown in Fig. \ref{KOLFEM-FEM}. Comparing them with Fig. \ref{KOLFEM-POD}, we can see that the method "Parareal-II" performs better than  the method "Parareal-I". However,  we also see that the accuracy cannot be improved any more once  it reaches a certain level.

	For this method, we show its computational cost in  Table~\ref{tableKolmoFEM_time} (unit:s), including  the computational cost of fine propagator $T_{\mathcal{F}}$ and coarse propagator $T_{\mathcal{G}}$ in each time subinterval  per iteration, the total computational cost per iteration $T_{\text{total}}$.  For this method,  
		\begin{align*}
		T_{\text{total}}=T_{\mathcal{F}}+N\times T_{\mathcal{G}},
		\end{align*}
		and $N$ is the total number of time subintervals. In our tests, $N=200$. We see from Table~\ref{tableKolmoFEM_time} that although the computational cost of the coarse propagator for one time subinterval is much smaller than that of the fine propagator,  the total cost over $200$ time subintervals is far greater than that of the fine propagator due to the fact that the coarse propagator runs sequentially over the entire time interval, which  significantly limits the speedup of the parareal method.
	\begin{table}[!htb]
		\caption{ The computational cost of the method "Parareal-II" per iteration for solving Kolmogorov flow (unit:s).}\centering
		\label{tableKolmoFEM_time}
		\begin{tabular}{|cccc|}
			\hline
			$\epsilon$ & $T_{\mathcal{F}}$   & $T_{\mathcal{G}}$ & $T_{\text{total}}$\\
			\hline
			0.5         & 43.87& 0.77     &197.87 \\
			\hline
			0.1         & 43.27   & 0.71          &185.27 \\
			\hline
		\end{tabular}
	\end{table}
	\FloatBarrier

\begin{figure}[H]
		\centering
		\subfloat[$\epsilon=0.5.$]{\includegraphics[width=0.5\textwidth]{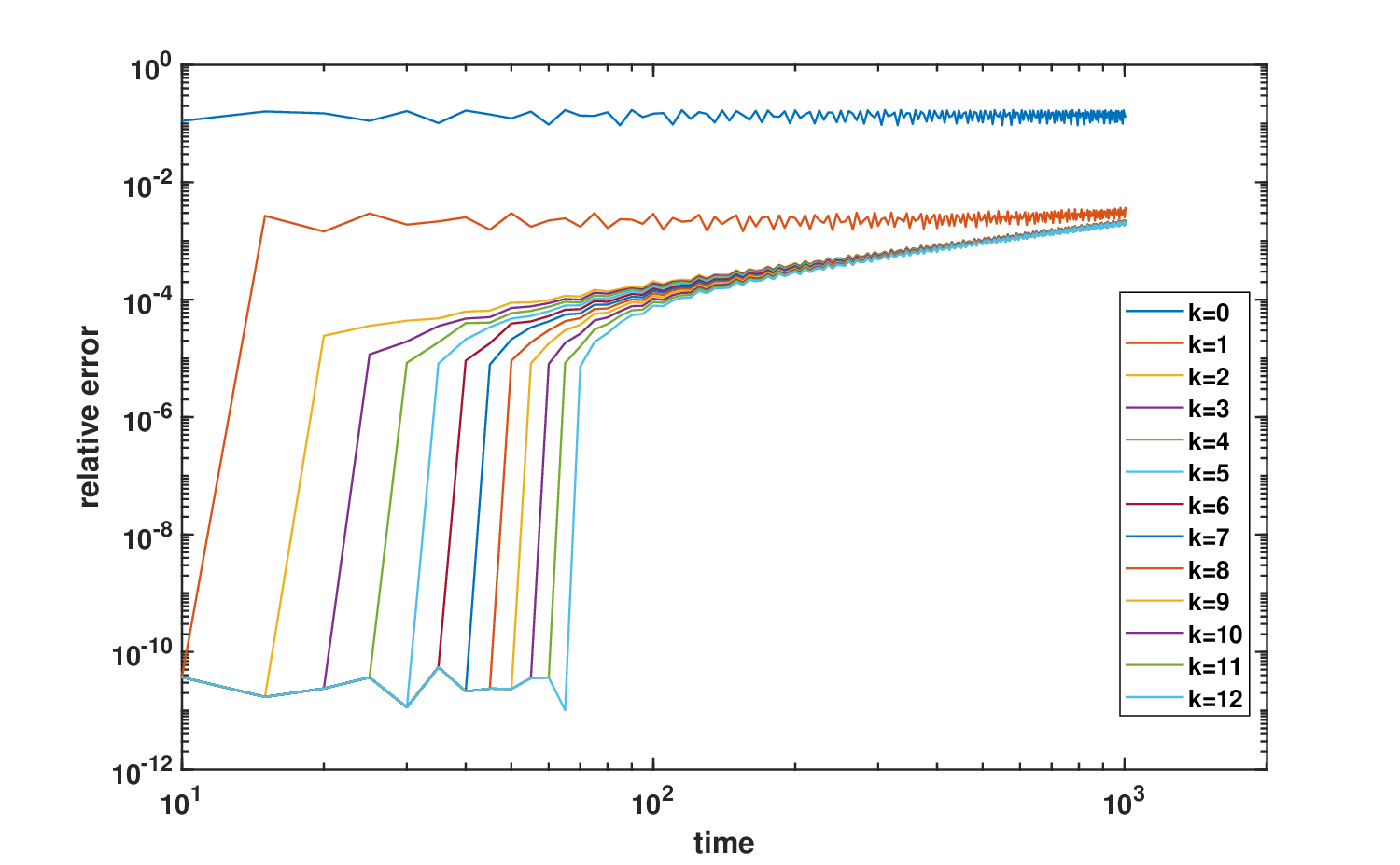}}
		\subfloat[$\epsilon=0.1.$]{\includegraphics[width=0.5\textwidth]{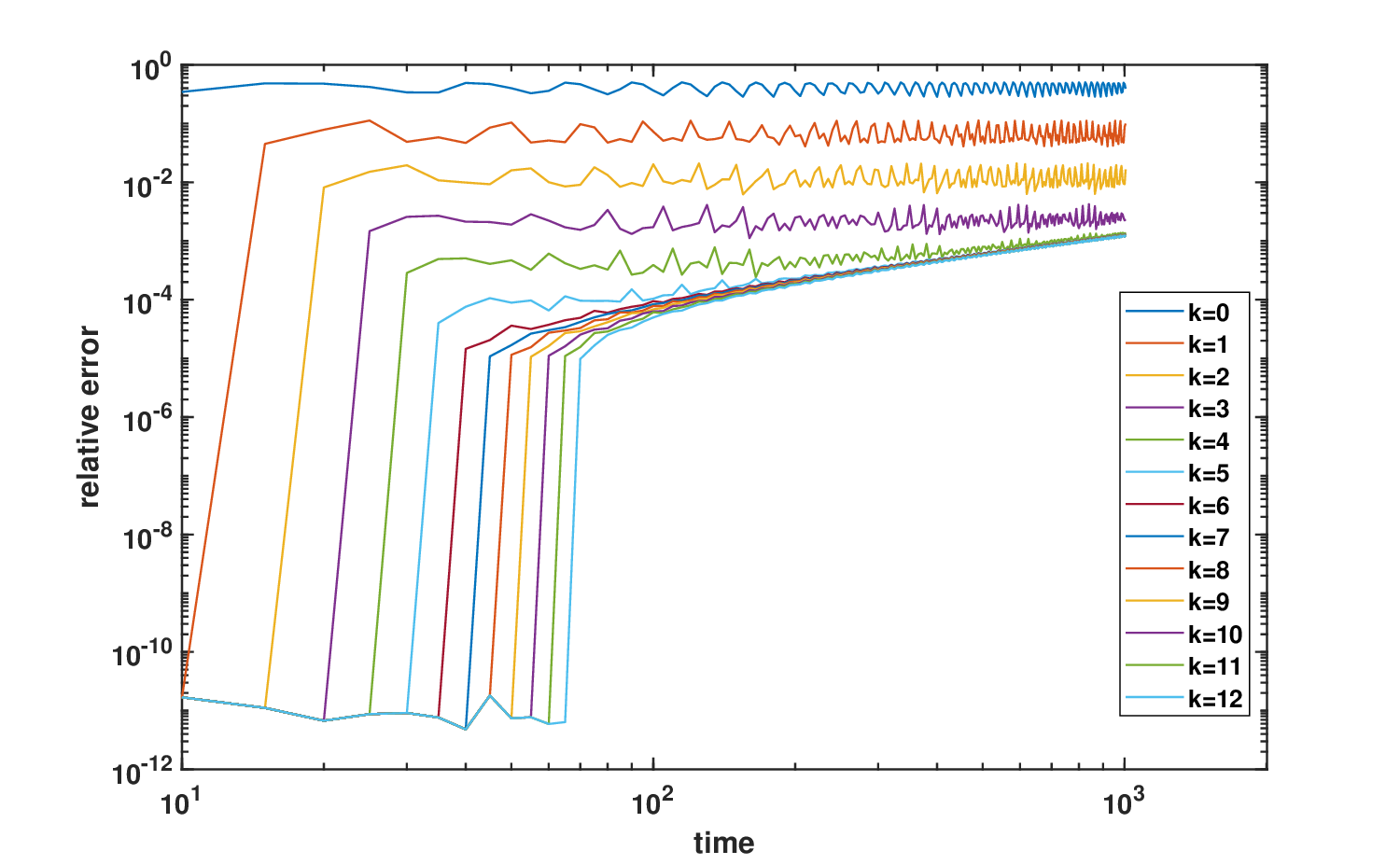}}
		
		\caption{The evolution curves of the relative error of $U^k_n$ obtained by the method "AdapParareal" with $(m_l,p)=(0,0)$  in each parareal iteration for solving the Kolmogorov flow with $\epsilon=0.5,0.1$}
		\label{KOLFEM-PODONE}
	\end{figure}

\begin{figure}[!tbhp]
		\centering
		\subfloat[$\epsilon=0.5.$]{\includegraphics[width=0.5\textwidth]{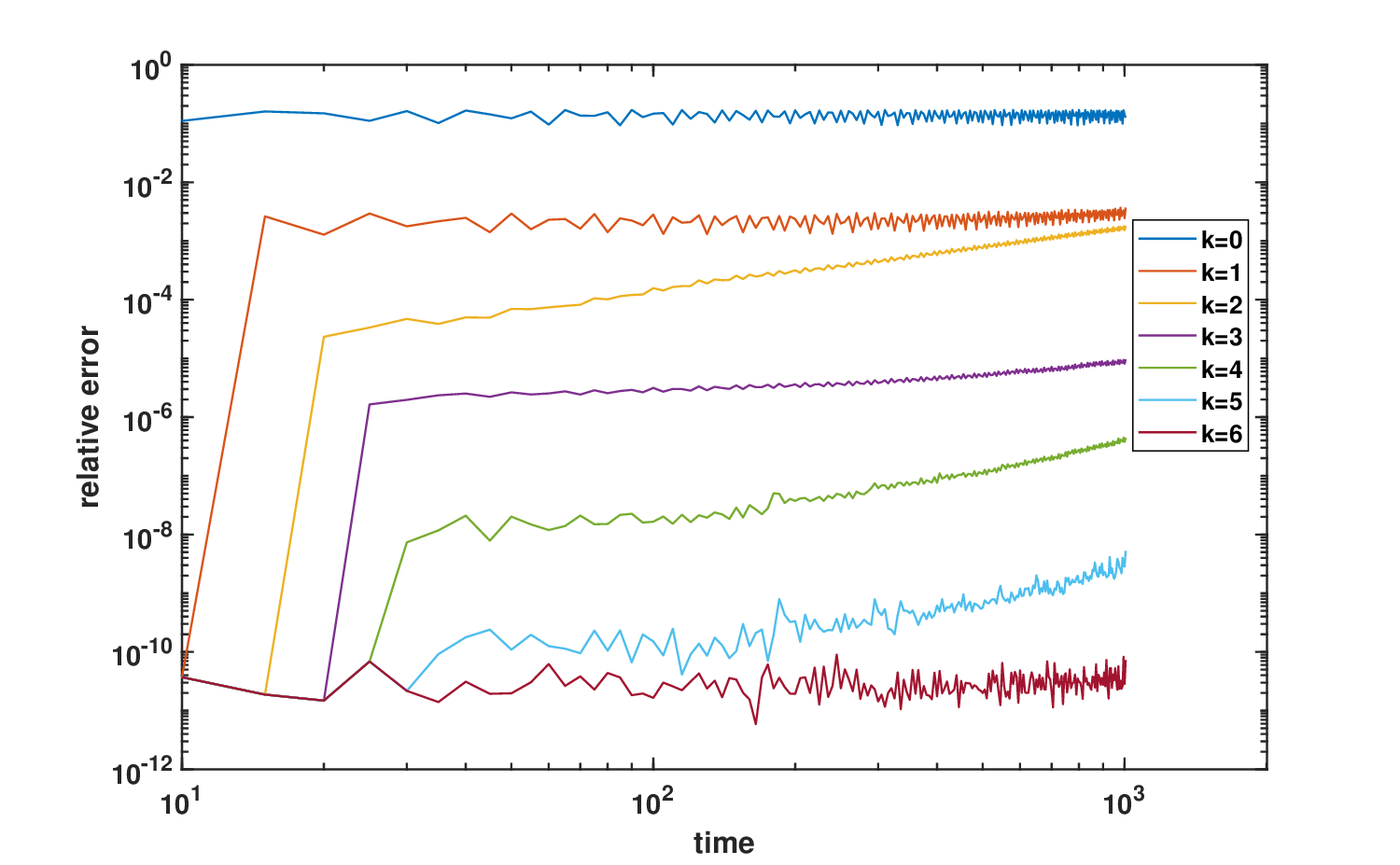}}
		\subfloat[$\epsilon=0.1.$]{\includegraphics[width=0.5\textwidth]{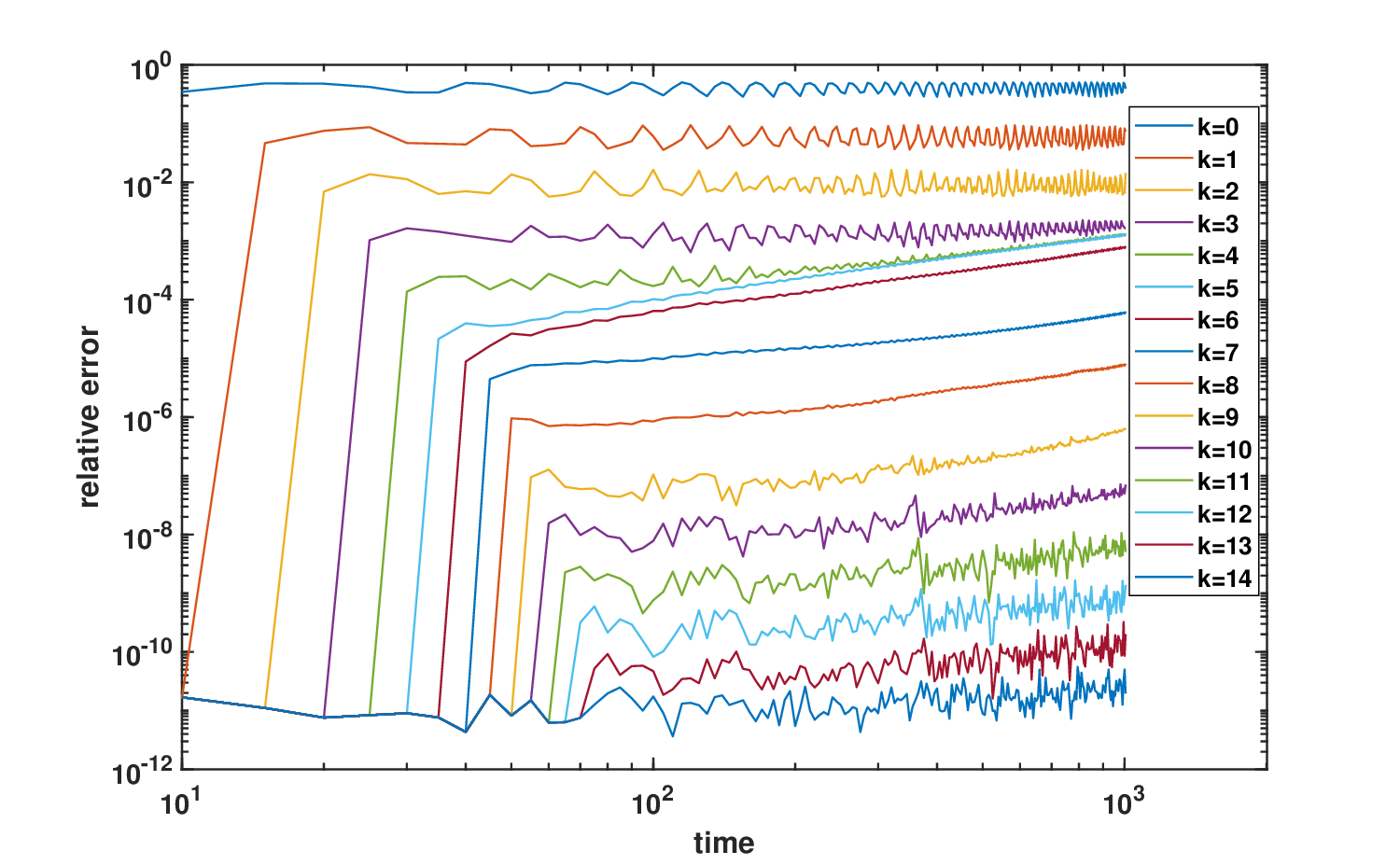}}
		\caption{The evolution curves of the relative error of $U^k_n$ obtained by the method "AdapParareal" with $(m_l,p)=(1,0)$  in each parareal iteration for solving the Kolmogorov flow with  $\epsilon=0.5,0.1$}
		\label{KOLFEM-APODlast}
	\end{figure}

\begin{figure}[!tbhp]
		\centering
		\subfloat[$\epsilon=0.5.$]{\includegraphics[width=0.5\textwidth]{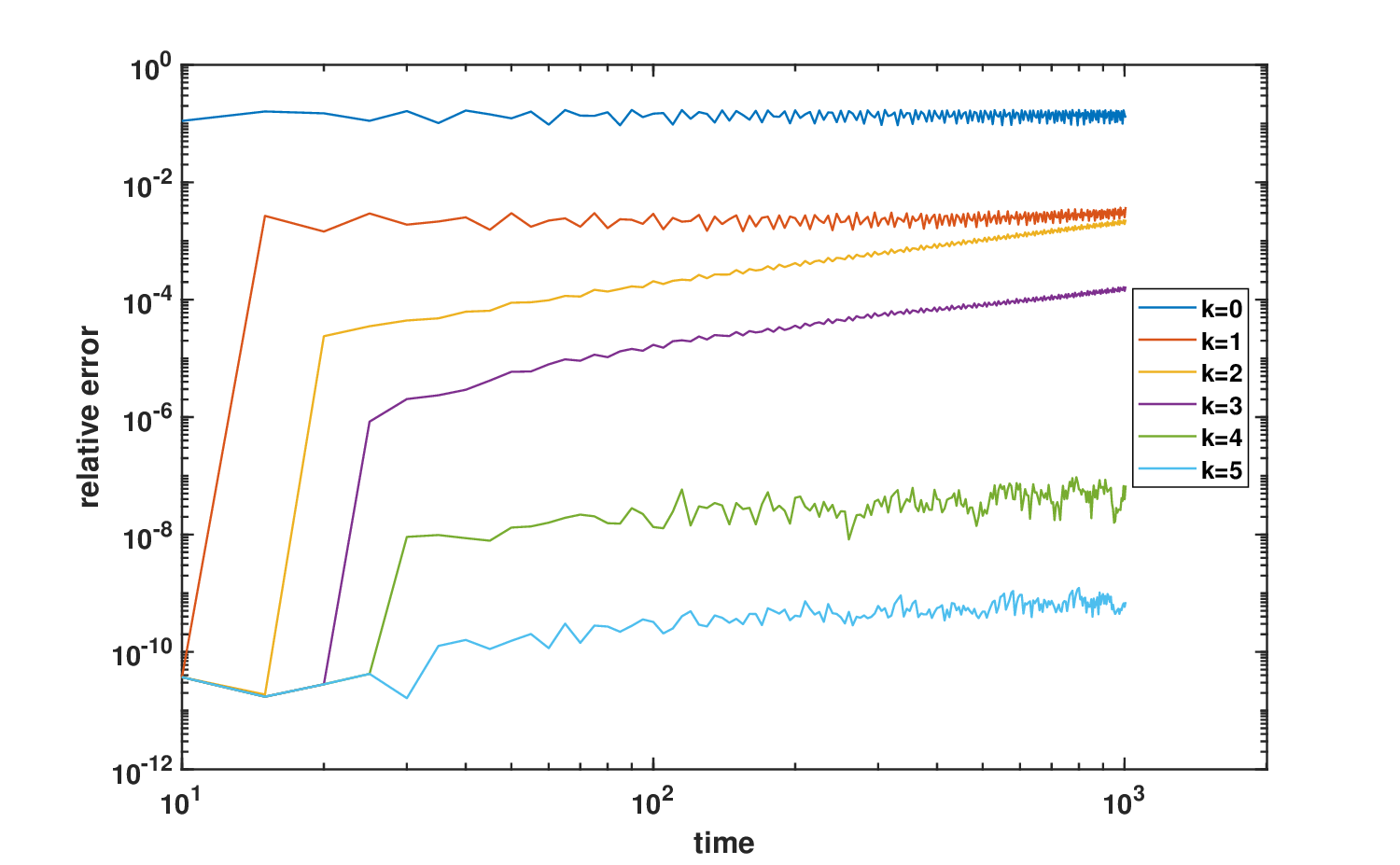}}
		\subfloat[$\epsilon=0.1.$]{\includegraphics[width=0.5\textwidth]{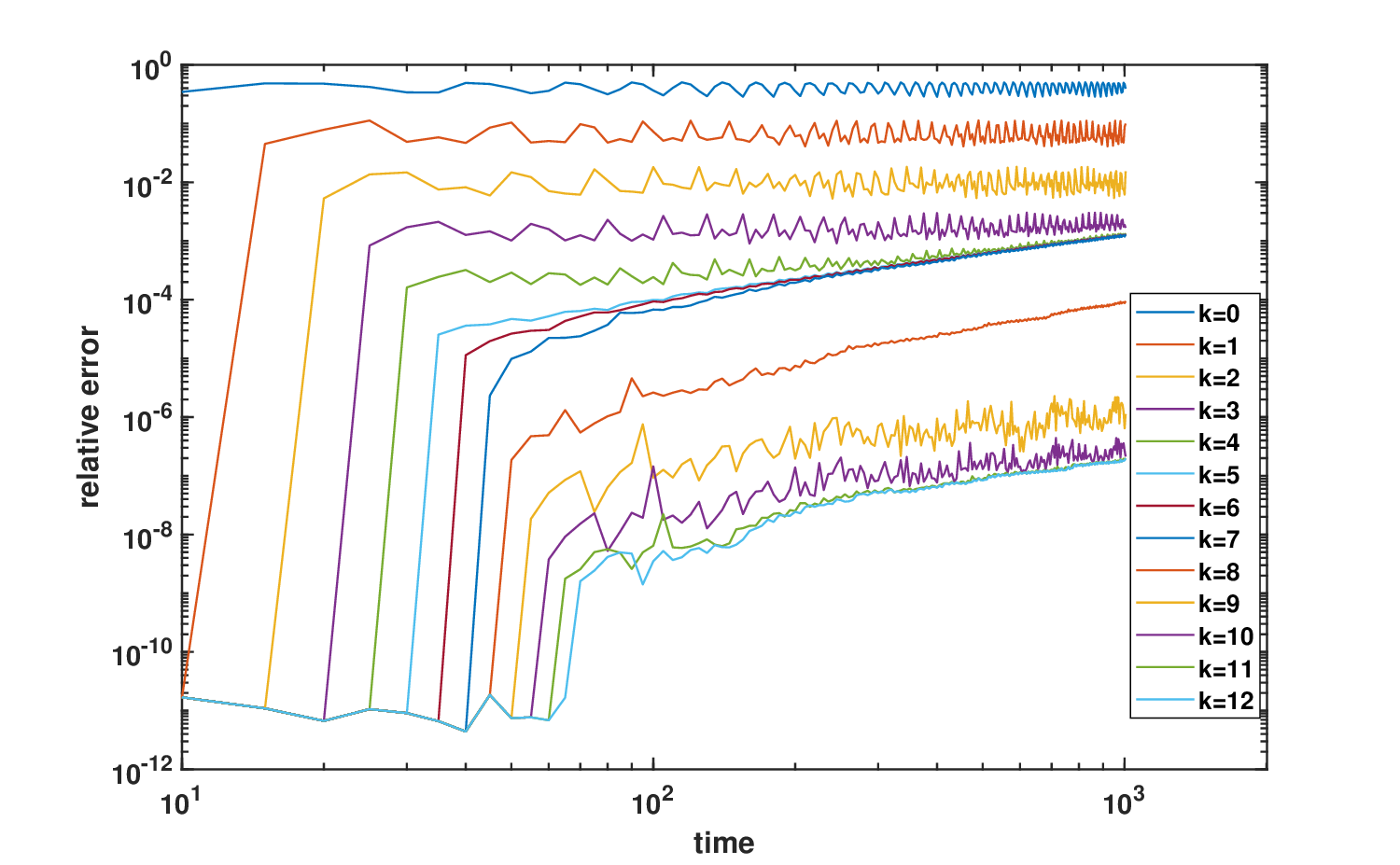}}
		\caption{The evolution curves of the relative error of $U^k_n$ obtained by the method "AdapParareal"  with $(m_l,p)=(0,1)$ in each parareal iteration for solving the Kolmogorov flow with $\epsilon=0.5,0.1$}
		\label{KOLFEM-APODtwo}
	\end{figure}

\begin{figure}[!tbhp]
		\centering
		\subfloat[$\epsilon=0.5.$]{\includegraphics[width=0.5\textwidth]{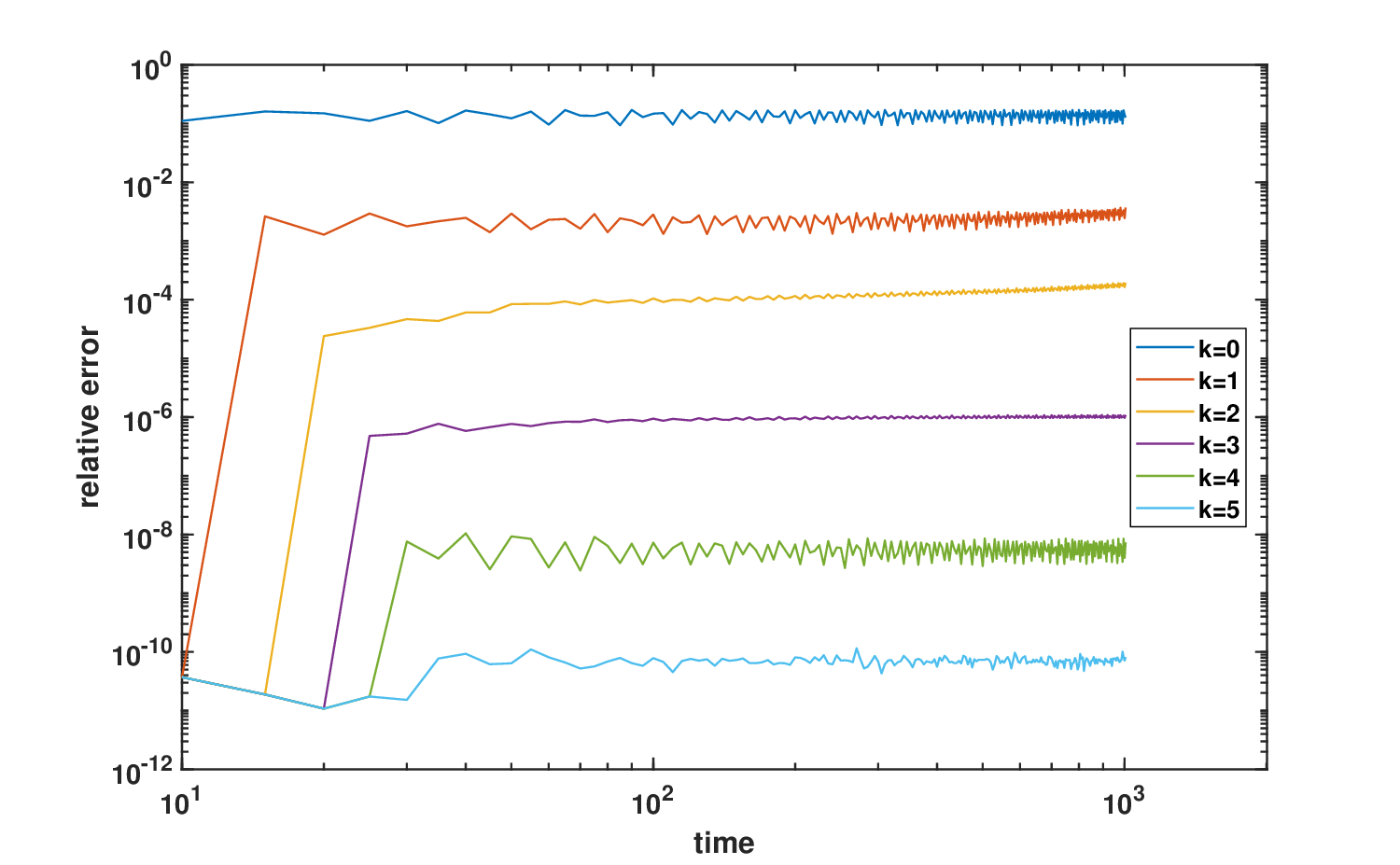}}
		\subfloat[$\epsilon=0.1.$]{\includegraphics[width=0.5\textwidth]{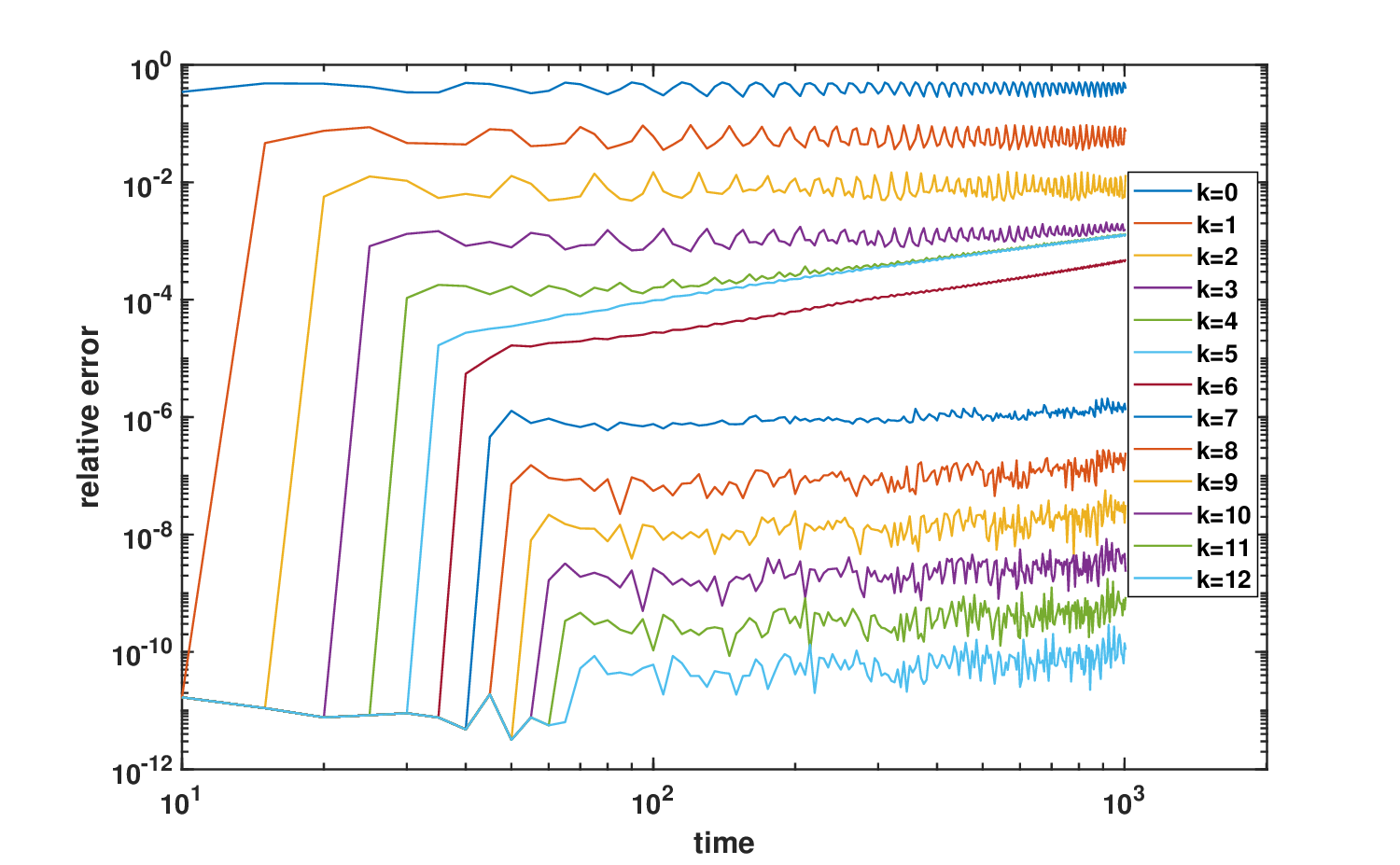}}
		\caption{The evolution curves of the relative error of $U^k_n$ obtained by the method "AdapParareal"  with $(m_l,p)=(1,1)$ in each parareal iteration for solving the Kolmogorov flow with $\epsilon=0.5,0.1$}
		\label{KOLFEM-APODTHREELONG}
	\end{figure}

	We then use these two cases with $\epsilon=0.5$ and $\epsilon=0.1$ to show the performance of our new method "AdapParareal" with different $(m_l,p)$. The dimensions of POD subspaces for the coarse propagator in the $0-$th iteration are $15$ and $20$ respectively. The corresponding results are presented in Fig. \ref{KOLFEM-PODONE}-\ref{KOLFEM-APODTHREELONG}. First, from these figures, it is easy to see that the relative error can reach at least $O(10^{-3})$ after several iterations. This indicates that our adaptive parareal method is effective. For the cases of $(m_l,p)=(1,0), (1,1)$, the accuracy can reach even $O(10^{-10})$ after several iterations. By comparing them with Fig. \ref{KOLFEM-POD}-\ref{KOLFEM-FEM}, we find that our adaptive parareal method has better performance than that of the method "Parareal-I" and "Parareal-II" in achieving high accuracy. Furthermore, the relative error does not increase rapidly over time, which indicates that our adaptive parareal method can perform well for simulating long-term evolution.

	\FloatBarrier
	
	\FloatBarrier
	
	\FloatBarrier
	
	\FloatBarrier 
	
	By comparing the numerical results using different $(m_l,p)$ in Fig.\ref{KOLFEM-PODONE}-\ref{KOLFEM-APODTHREELONG}, we observe that as the number of time subintervals used in constructing the coarse propagator increases, the convergence becomes faster. 
		We have done further tests using larger $m_l$ or $p$, and have not seen obvious improvement in the accuracy. Note that selecting more time subintervals means more computational cost for getting the POD subspace, and will result in more POD modes, which will increase the computational cost in further. Therefore, considering the accuracy and the cost, we recommend to set $(m_l, p) = (1, 1)$.
	 
\begin{figure}
		\centering
		\subfloat[$\epsilon=0.05.$]{\includegraphics[width=0.5\textwidth]{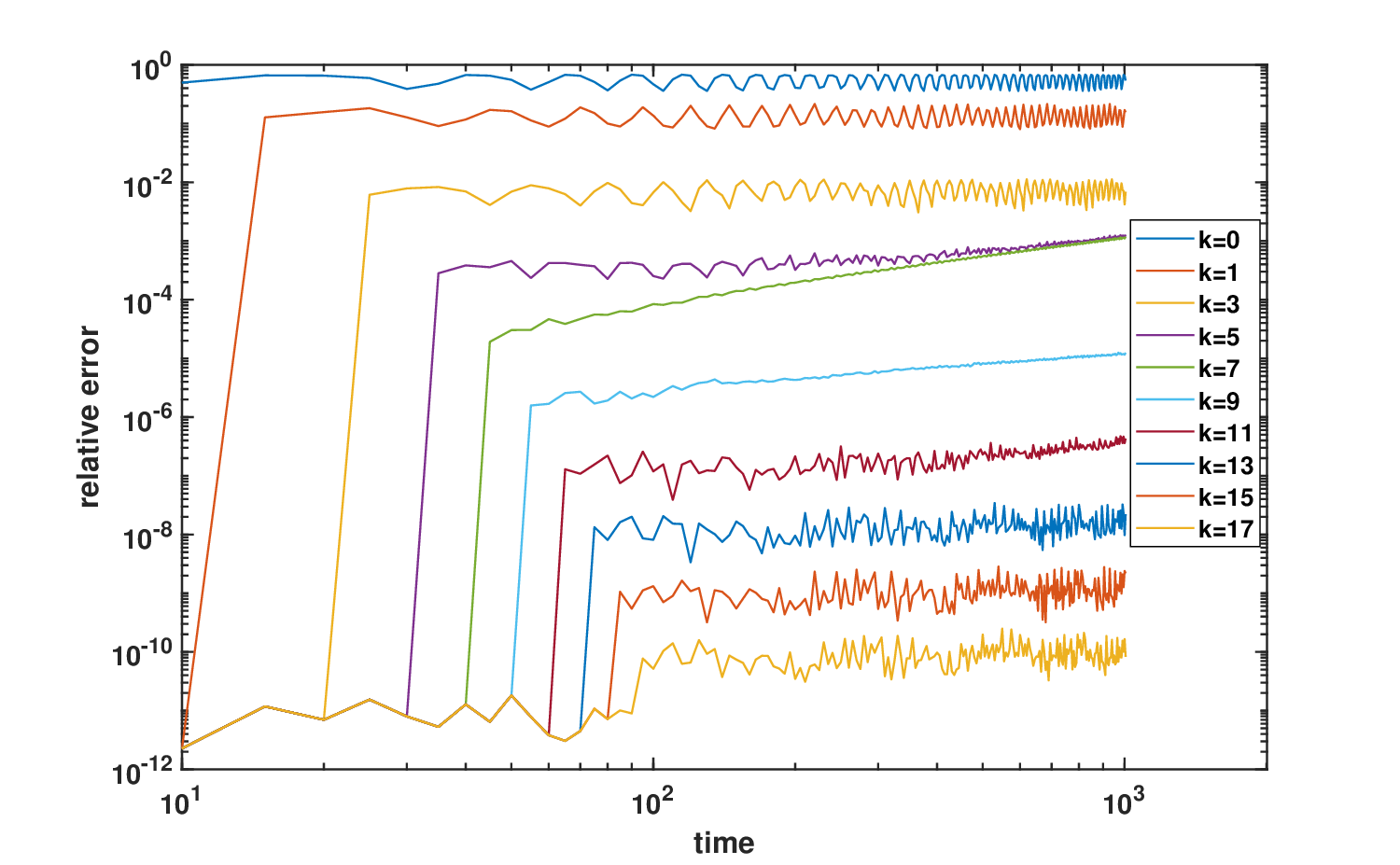}}
		\subfloat[$\epsilon=0.01.$]{\includegraphics[width=0.5\textwidth]{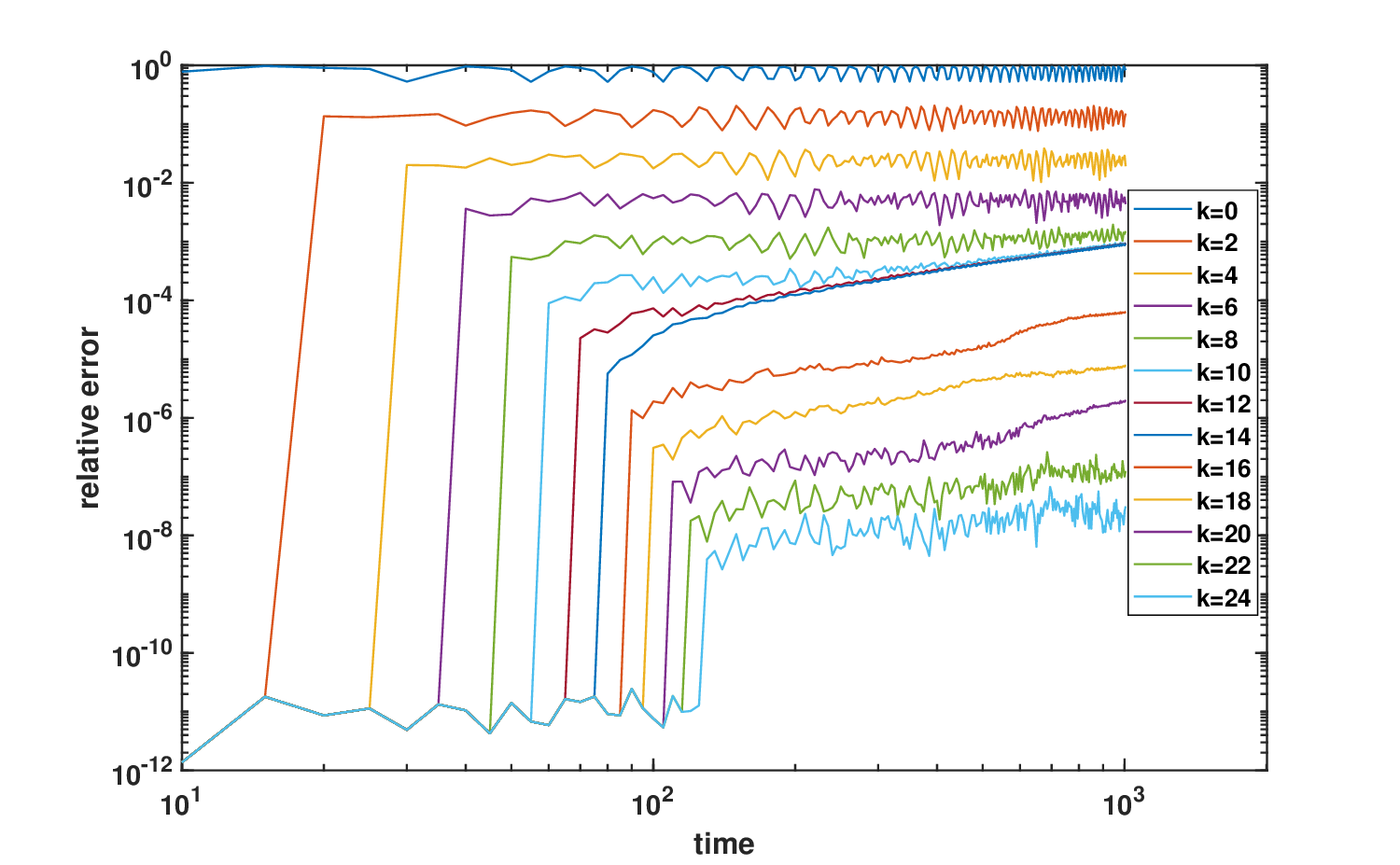}}
		
		\subfloat[$\epsilon=0.005.$]{\includegraphics[width=0.5\textwidth]{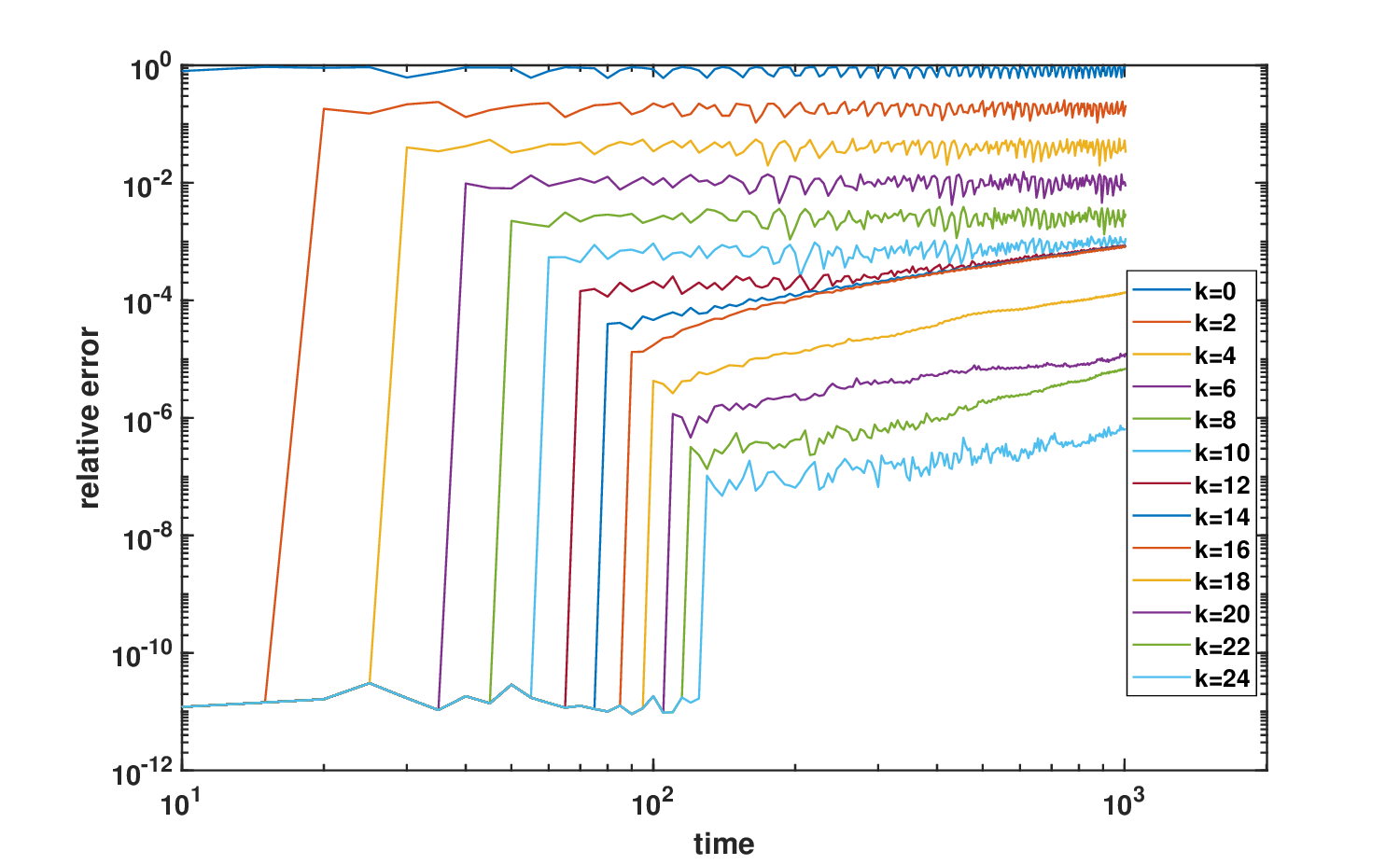}}
		\caption{The evolution curves of the relative error of $U^k_n$ obtained by the method ``AdapParareal'' with $(m_l,p)=(1,1)$  in each parareal iteration for solving the Kolmogorov flow with $\epsilon=0.05,0.01,0.005$}
		\label{KOLFEM-APODTHREE2}
	\end{figure}

	Motivated by the above observations, when turning to the more challenging advection-dominated problems with $\epsilon=0.05,0.01,0.005$, we adopt the setting $(m_l,p)=(1,1)$. The corresponding numerical results are reported in Fig.~\ref{KOLFEM-APODTHREE2}. The dimensions of the POD subspaces for the coarse propagator in the $0-$th iteration for these three cases are $24$, $32$, and $35$, respectively.

From Fig. \ref{KOLFEM-APODTHREE2}, we can first see that, as the parareal iteration continues, the error decreases obviously. However, compared with the cases $\epsilon = 0.5$ and $\epsilon = 0.1$, the convergence of the adaptive parareal method for these cases is slower. The reason lies in the fact that the smaller $\epsilon$, the more difficult it is to simulate. It is well documented in the literature that predictor-corrector parallel-in-time methods, such as the parareal method, typically suffer from slower convergence rates and reduced stability when applied to advection-dominated problems \cite{bal2005convergence, gander2007analysis}. Fortunately, after several iterations, the accuracy obtained by our adaptive parareal method can be very high in long-term evolution. In fact, we can see that even for the case of $\epsilon=0.005$, which is very difficult to simulate well, after $24$ parareal iterations, the error obtained on the whole time interval can be smaller than $2\times10^{-6}$.  
	
	To estimate the speedup achieved by our method, we measure the wall time required by each part per iteration (iteration $k\geq 1$) and summarize the results for solving the Kolmogorov flow in the Table~\ref{tableKolmo_time}  (unit: s). Here
			$$T_{\text{total}}=T_{\mathcal{F}}+T_U+T_C+N\times(T_{\mathcal{G}}+T_A)$$
		with $N$ the total number of time subintervals. 
		
			As shown in Table~\ref{tableKolmo_time}, compared with the time cost for the fine propagator in one time subinterval ($T_\mathcal{F}$), the time cost for the POD basis construction is very cheap.  In contrast to the method "Parareal-II", where the coarse propagator is the bottleneck, our method achieves a more balanced workload because only the propagation of the coarse propagator and the augmentation of the reduced spaces require sequential execution. Furthermore, our coarse propagator involves significantly fewer degrees of freedom than that of the method "Parareal-I". The maximum dimension of the POD subspaces $m_\text{max}$ for different $\epsilon$ is listed in Table~\ref{tableKolmo_PODnum}. 
	
	\begin{table}
		\caption{ The time cost of the method "AdapParareal" with $(m_l,p)=(1,1)$ in each iteration  for solving Kolmogorov flow with different $\epsilon$.}\centering
		\label{tableKolmo_time}
		\begin{tabular}{|ccccccc|}
			\hline
			$\epsilon$ & $T_{\mathcal{F}}$          & $T_{\mathcal{G}}$ & $T_U$  & $T_C$ & $T_A$ & $T_{\text{total}}$\\
			\hline
			 0.5         & 44.16& 0.04         & 0.41& 0.25     & 0.11 &74.82 \\
			\hline
			0.1         & 43.39   & 0.07         &0.97& 0.67    & 0.12 &83.03 \\
			\hline
			0.05          & 39.56 & 0.08      &1.34     &1.06 & 0.13   & 83.96  \\
			\hline
			0.01           & 42.19 & 0.13      &3.09   & 2.40 & 0.15    & 103.68  \\
			\hline
			0.005         & 44.14 &  0.15   &     3.94   &3.14 &  0.16   &  113.22 \\
			\hline  
		\end{tabular}
	\end{table}
	\begin{table}
		\caption{ The maximum dimension of the POD subspaces $m_\text{max}$ for the method "AdapParareal" with $(m_l,p)=(1,1)$ for solving Kolmogorov flow with different $\epsilon$  (unit:s).}\centering
		\label{tableKolmo_PODnum}
		\begin{tabular}{|cccccc|}
			\hline
			$\epsilon$ & $0.5$          & $0.1$ & $0.05$  & $0.01$ & $0.005$\\
			\hline
			$m_\text{max}$         & $45$& $74$         & $93$& $147$     & $172$  \\
			\hline 
		\end{tabular}
	\end{table}
	
	We evaluate the speedup as
		$$\text{speedup}\approx \dfrac{T_{\mathcal{F}}\times N}{T_{\text{total}}\times (k+1)}$$
		since the $0$-th iteration is cheaper than the subsequent $k$-th ($k\geq 1$) iteration.  In our experiments, the number of time subintervals $N=200$.  As shown in the Fig. \ref{KOLFEM-APODTHREELONG}-\ref{KOLFEM-APODTHREE2}, the error obtained by our method is of $O(10^{-6})$ by using less than $10$  iterations for $\epsilon=0.5,0.1,0.05$, and  is of $O(10^{-6})$ within about $20$  iterations for $\epsilon=0.01, 0.005$ for solving Kolmogorov flow.   Therefore, based on our analysis for the speedup in  Section~\ref{sec: Complex},  our method "AdapParareal"  achieves a speedup of about $9.5$ for cases of $\epsilon=0.5,0.1,0.05$, and  achieves a speedup of about $3.7$ for cases of $\epsilon=0.01, 0.005$.
	
	In the introduction, we conjectured that as the iteration proceeds, the accuracy of the spatial discretization for the coarse propagator will approach that of the fine propagator's spatial discretization  progressively. We want to point out that we have also done some tests to demonstrate it. Here, to remove the impact of the temporal discretization, we choose to use the same temporal discretization scheme and time step size used in the fine propagator. That is, the only difference between fine and coarse propagator lies in the spatial discretization, one discretizes the problem  in the finite element space, the other discretizes the problem in the adaptive POD subspace. That is, we set $dT = \delta t = 0. 01$. To see the performance, we show the difference between the approximation obtained by the coarse propagator and the fine propagator starting from the same initial values, and the relative error of $\mathcal{G}_k\left(U_{n-1}^{k}\right)$. That is, we show $\frac{||\mathcal{F}\left(U_{n-1}^{k-1}\right)-\mathcal{G}_k\left(U_{n-1}^{k-1}\right)||_{L^2}}{||\mathcal{F}\left(U_{n-1}^{k-1}\right)||_{L^2}}$ and $\frac{||U_{n}-\mathcal{G}_k\left(U_{n-1}^{k}\right)||_{L^2}}{||U_{n}||_{L^2}}$. 
	
\begin{figure}
		\centering
		\subfloat[$\epsilon=0.05.$]{\includegraphics[width=0.5\textwidth]{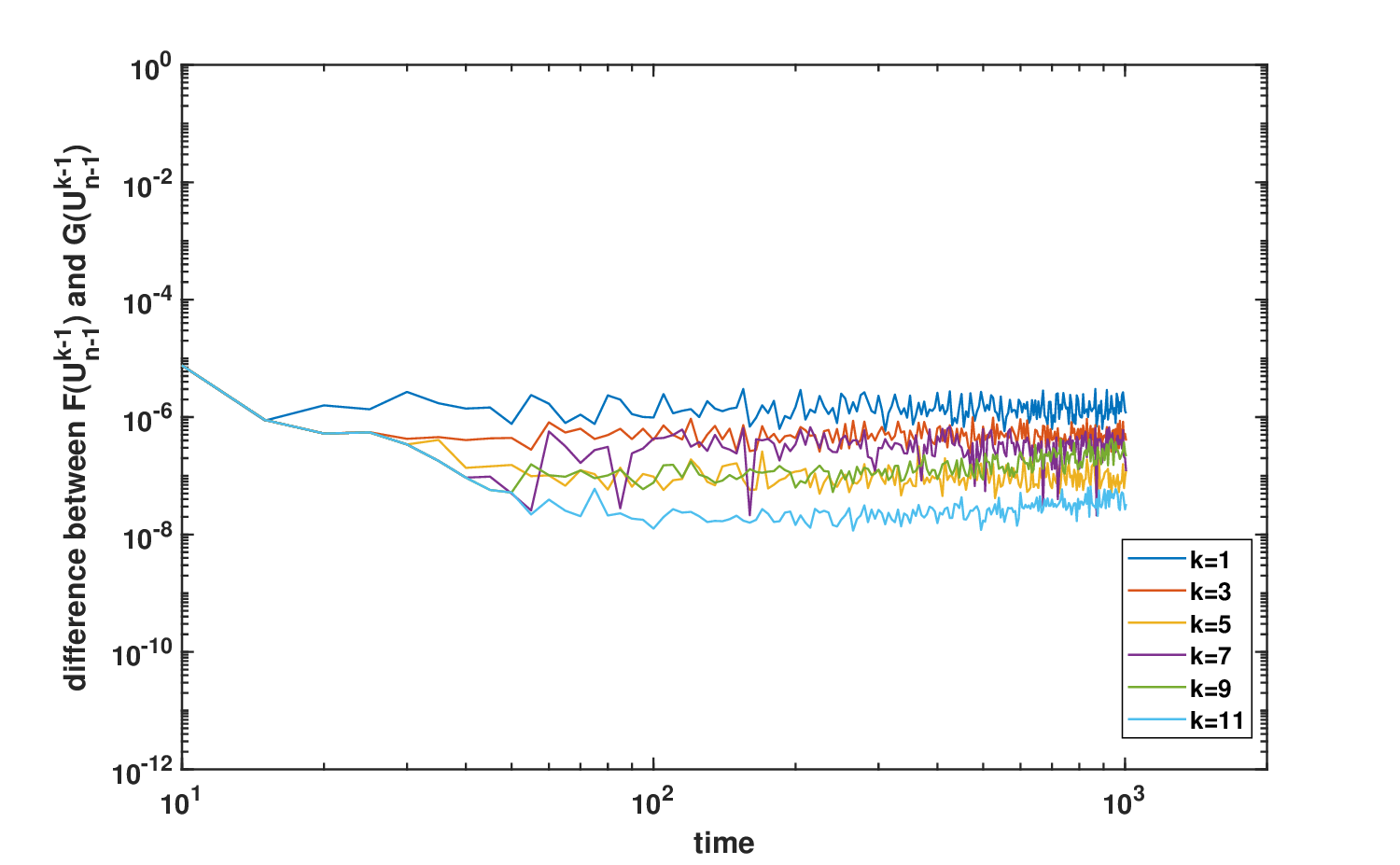}}
		\subfloat[$\epsilon=0.01.$]{\includegraphics[width=0.5\textwidth]{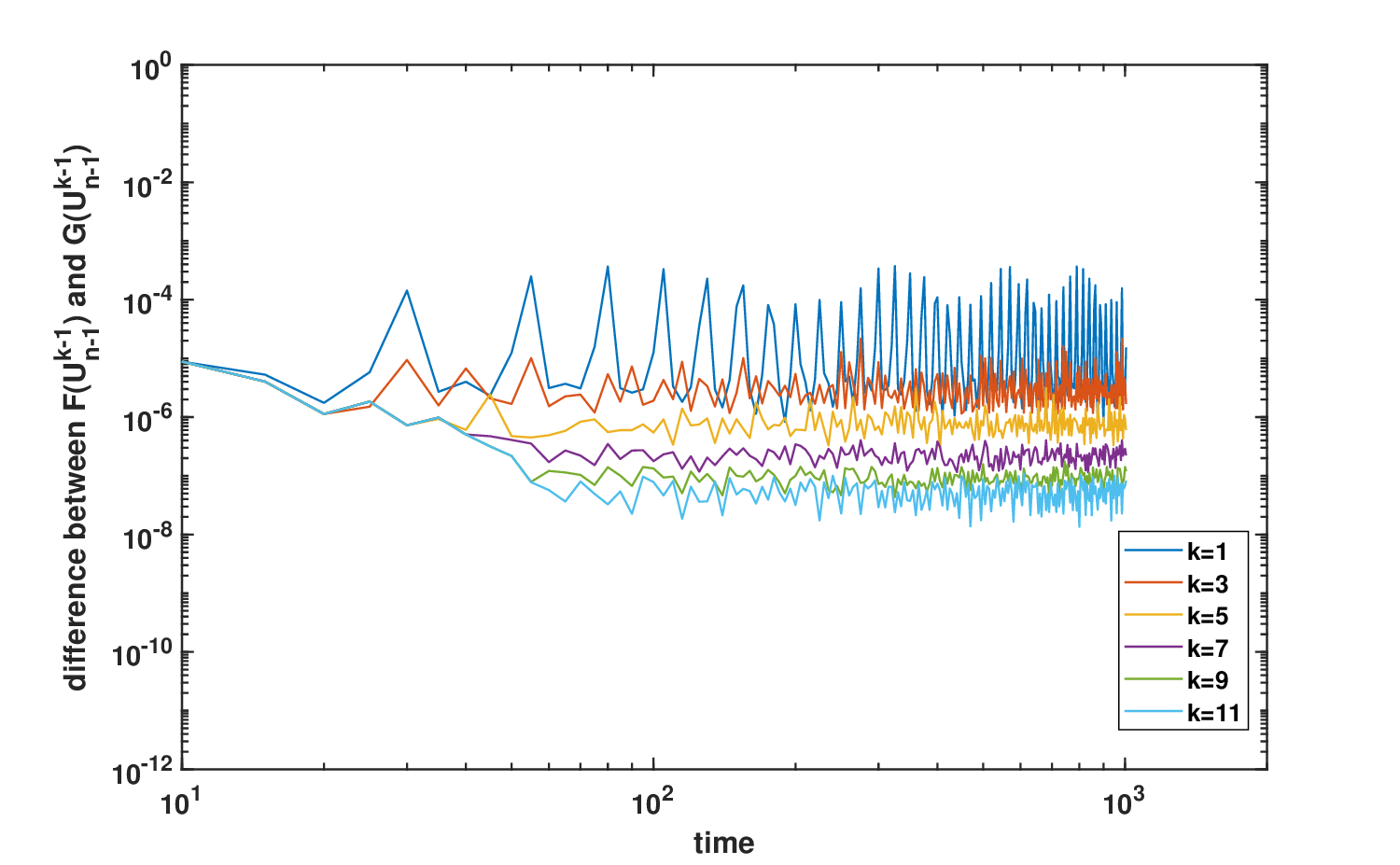}}
		
		\subfloat[$\epsilon=0.005.$]{\includegraphics[width=0.5\textwidth]{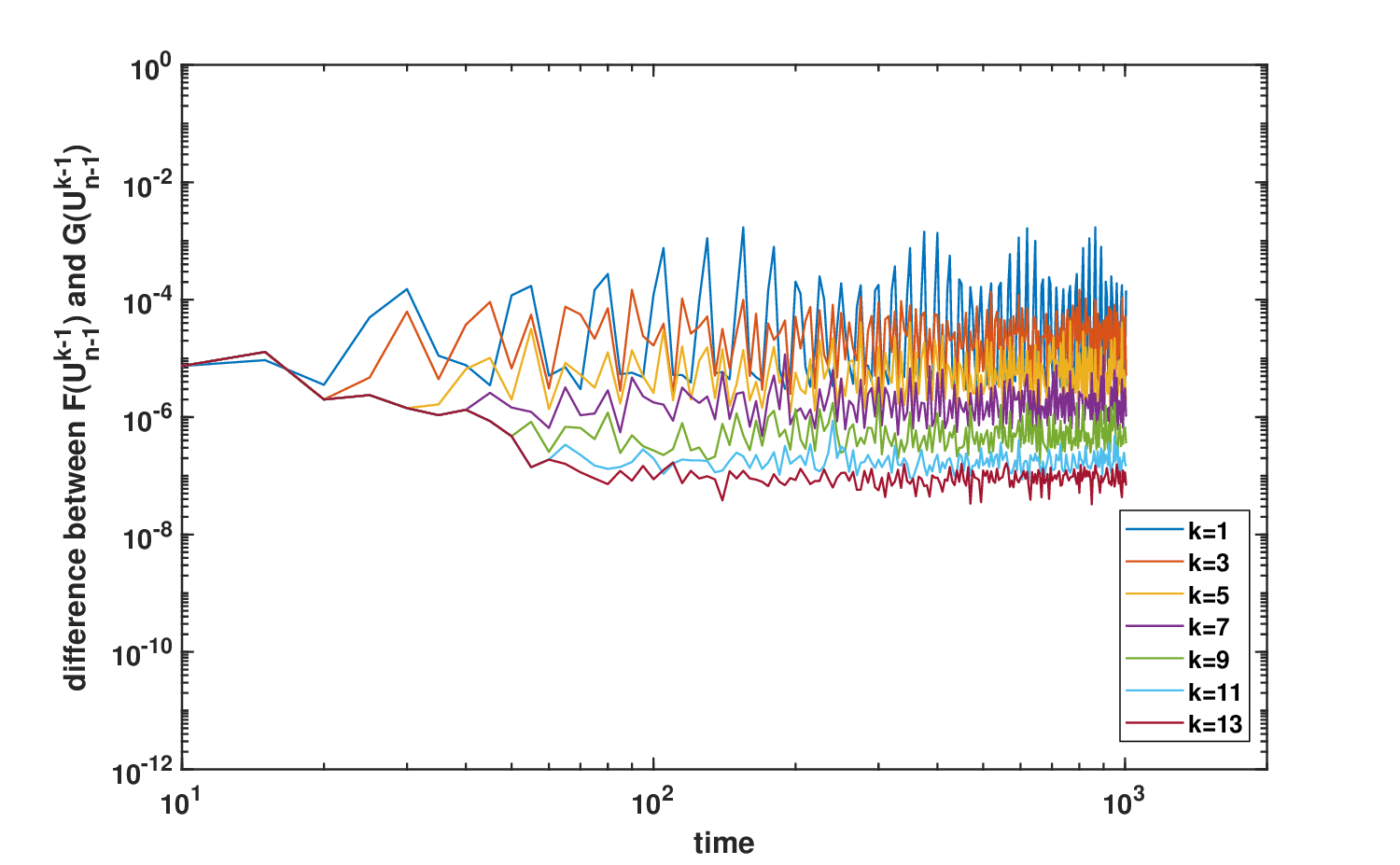}}
		\caption{The evolution curves of $\frac{||\mathcal{F}\left(U_{n-1}^{k-1}\right)-\mathcal{G}_k\left(U_{n-1}^{k-1}\right)||_{L^2}}{||\mathcal{F}\left(U_{n-1}^{k-1}\right)||_{L^2}}$ obtained by the method "AdapParareal" with $(m_l,p)=(1,1)$  in each parareal iteration for solving the Kolmogorov flow with $\epsilon=0.05,0.01,0.005$}
		\label{KOLFEM-APODextra1}
	\end{figure}

\begin{figure}
		\centering
		\subfloat[$\epsilon=0.05.$]{\includegraphics[width=0.5\textwidth]{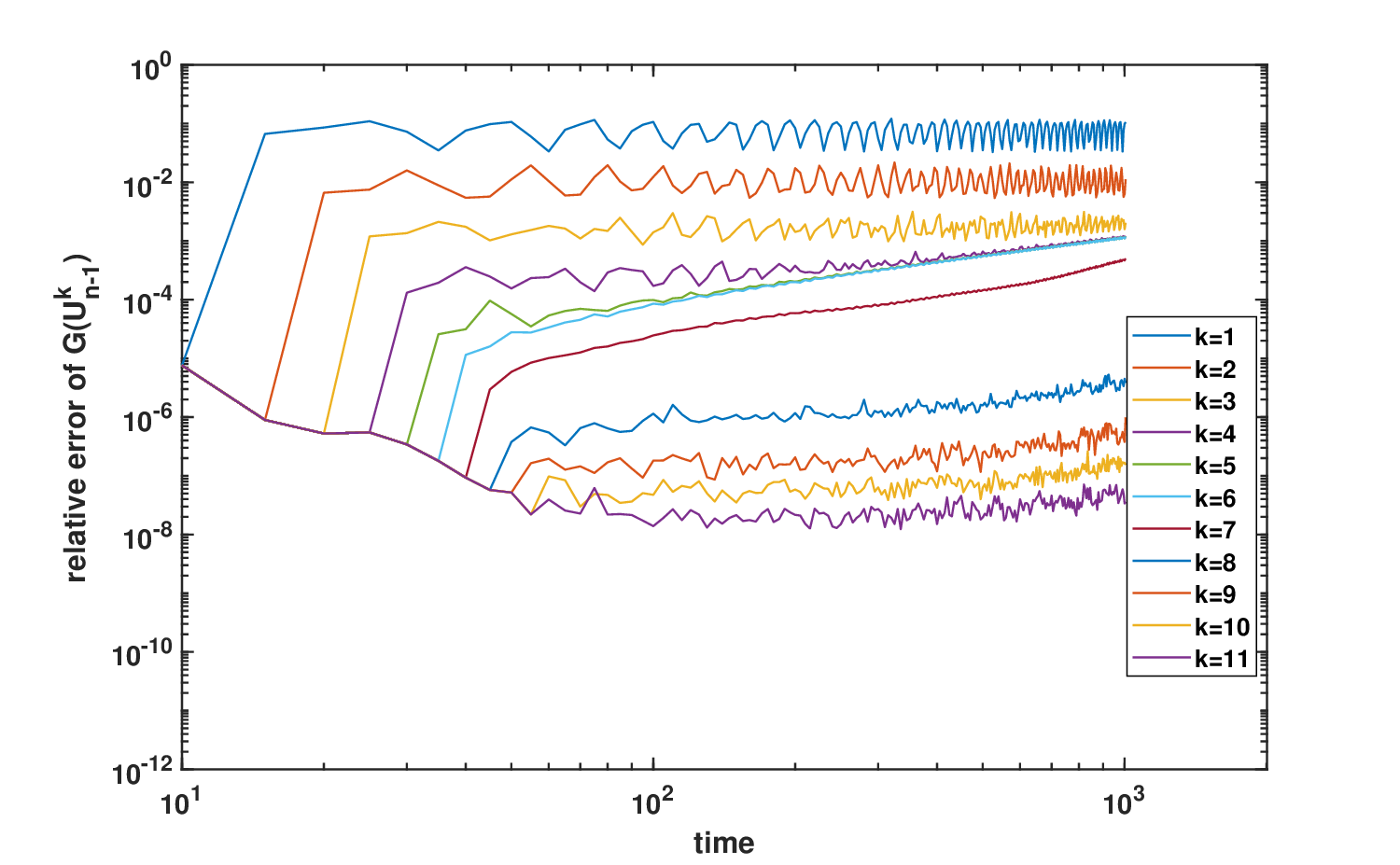}}
		\subfloat[$\epsilon=0.01.$]{\includegraphics[width=0.5\textwidth]{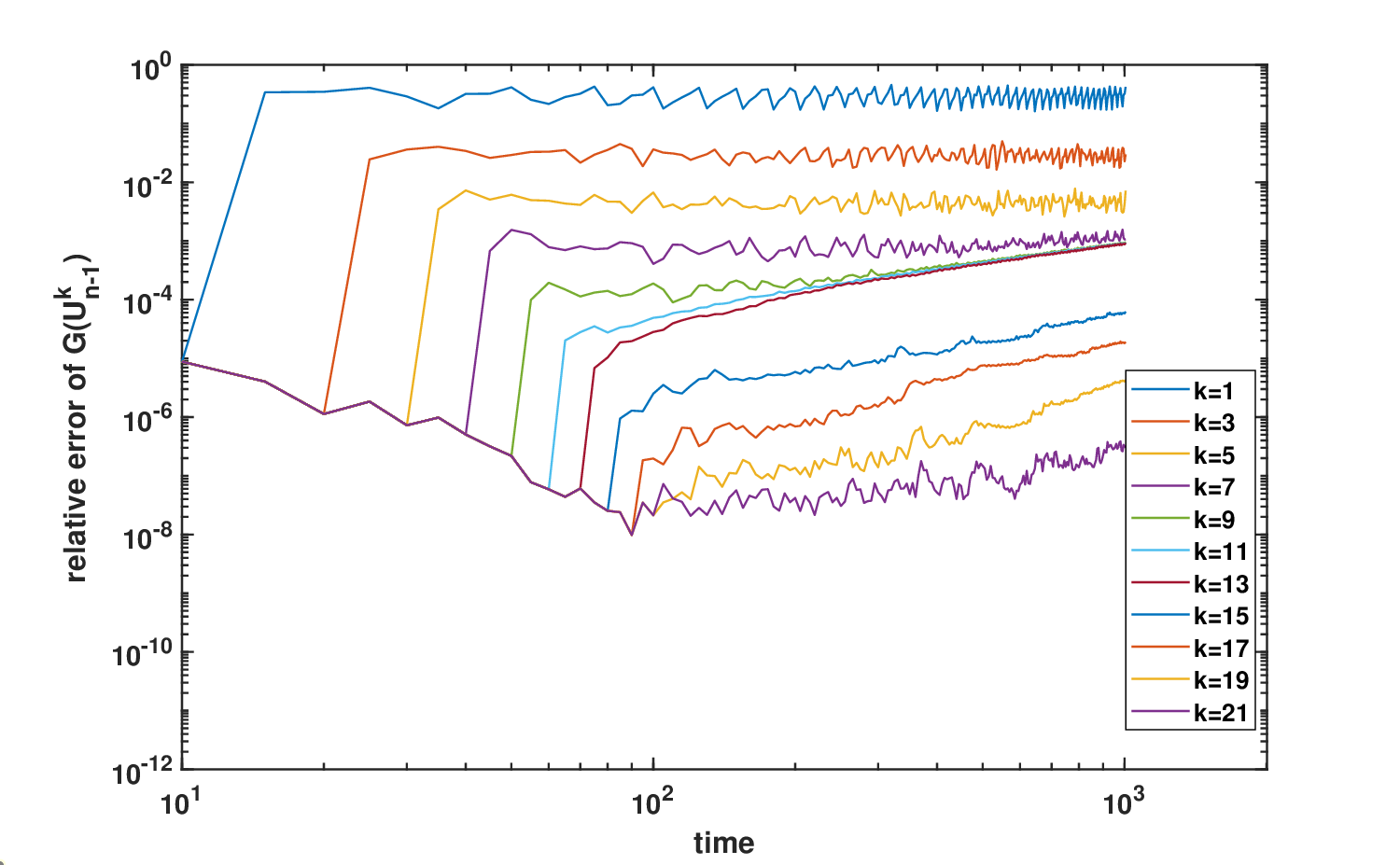}}
		
		\subfloat[$\epsilon=0.005.$]{\includegraphics[width=0.5\textwidth]{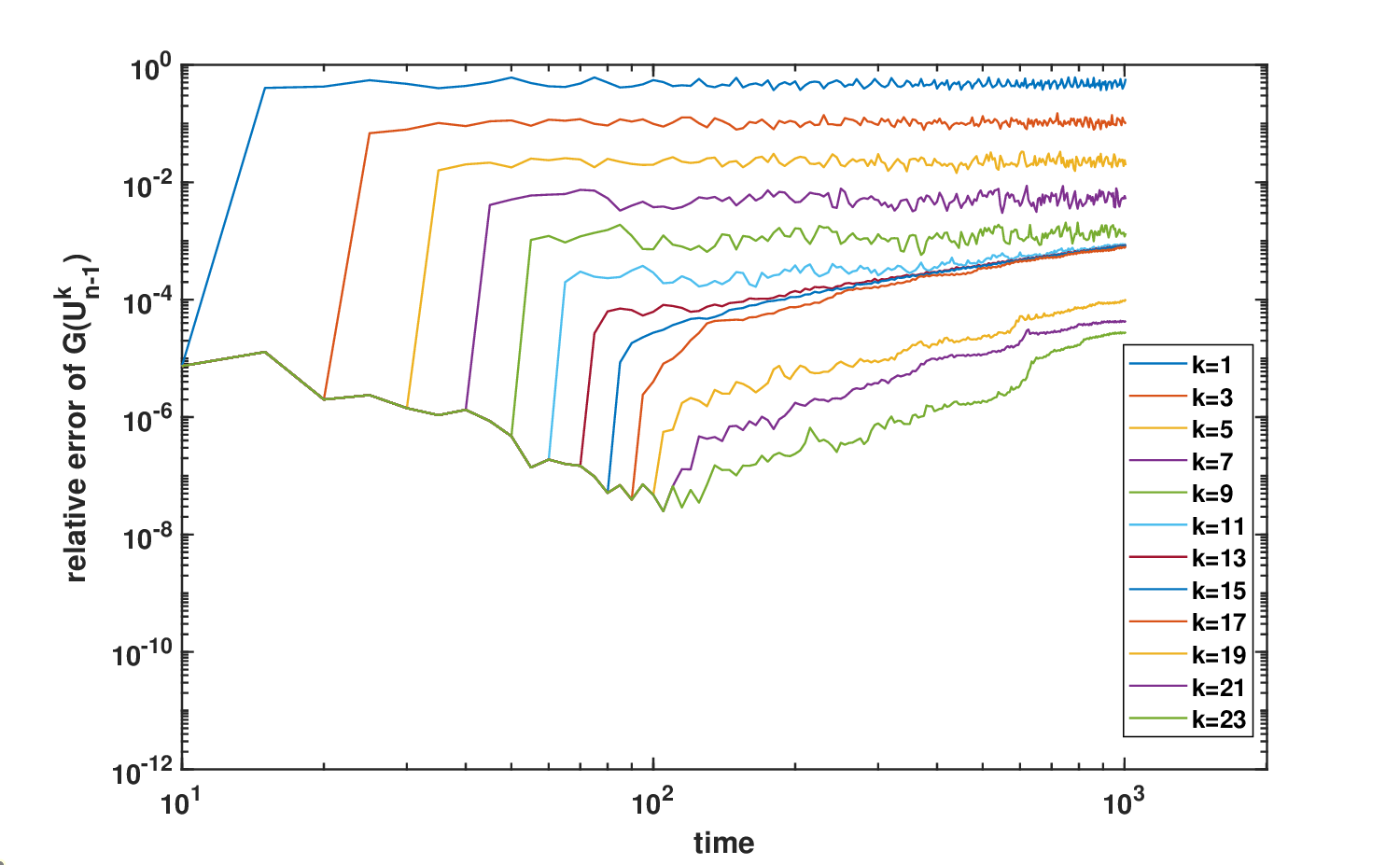}}
		\caption{The evolution curves of $\frac{||U_{n}-\mathcal{G}_k\left(U_{n-1}^{k}\right)||_{L^2}}{||U_{n}||_{L^2}}$ obtained by the method "AdapParareal" with $(m_l,p)=(1,1)$ in each parareal iteration for solving the Kolmogorov flow with $\epsilon=0.05,0.01,0.005$}
		\label{KOLFEM-APODextra}
	\end{figure}

	The numerical results for the Kolmogorov flow with $\epsilon=0.05,0.01,0.005$ respectively are presented in Fig.~\ref{KOLFEM-APODextra1}-\ref{KOLFEM-APODextra}. First, from these figures, we find that $\frac{||\mathcal{F}\left(U_{n-1}^{k-1}\right)-\mathcal{G}_k\left(U_{n-1}^{k-1}\right)||_{L^2}}{||\mathcal{F}\left(U_{n-1}^{k-1}\right)||_{L^2}}$ can reach at least $O(10^{-4})$ for all subintervals and iterations. Furthermore, from these sub-figures, we observe that as the increase of the parareal iteration, the accuracy of the coarse propagator improves significantly. After several iterations, the accuracy of the coarse propagator becomes very high. This shows that as the parareal iteration proceeds,  the accuracy of spatial discretization of the coarse propagator can indeed be improved. From these tests, we can see that our adaptive parareal method can also be viewed as a parareal based adaptive POD method, which uses data obtained by fine propagator in each time subinterval during the parareal iteration to update the POD subspace adaptively. By using the idea of the parareal method, for this adaptive parareal adaptive POD method, the update of the POD subspace for each time subinterval can be carried out essentially in parallel. As the iteration continues, the accuracy of the POD method improves significantly.  
	
	\FloatBarrier
	
	\subsection{Arnold-Beltrami-Childress (ABC) flow}\label{ABC flow}
	We then consider the advection-diffusion equation \eqref{eq1} where the advection field is defined by the Arnold-Beltrami-Childress (ABC) flow \cite{brummell2001linear,  xin2016periodic}. The components of the equation are given as follows: 
	\begin{equation}
	\begin{aligned}
	\mathbf{B}(x,  y,  z,  t) &=(\sin (z+\sin w t)+\cos (y+\sin w t),  \sin (x+\sin w t) \\
	&\quad+\cos (z+\sin w t),  \sin (y+\sin w t)+\cos (x+\sin w t)),  \\
	f(x,  y,  z,  t)&=-\sin (z+\sin w t)-\cos (y+\sin w t),  \\
	c(x,  y,  z,  t) &= 0,  
	~h(x,  y,  z) =0, ~L= 2 \pi, ~w=1.0,~T=1005.0.
	\end{aligned}
	\end{equation}
	
	In this example, $\mathbf{B}$ and $f$ are also separable in time and space. We also use piecewise linear functions as the basis functions for the finite element method. We first divide $\Omega$ into $6$ tetrahedrons to create the initial grid and refine the initial mesh $21$ times uniformly using bisection to obtain the final mesh. The number of degrees of freedom is $2097152$. We fix the parameters $\delta t=1\times 10^{-2}$, $dT=0.5$, $\Delta T=T_0=5.0$, $\delta M=5$, $\gamma_1=\gamma_2=1.0-5.0\times10^{-6}$ and $\gamma_3=1.0-2.0\times10^{-8}$. The number of time subintervals is also set to $200$. 
	
	For this example, similar to the case of the Kolmogorov flow, we will also use the two cases of $\epsilon=0.5$ and $\epsilon=0.1$ to show the performance of the plain parareal methods, including the method "Parareal-I" and "Parareal-II".  
	
\begin{figure}[H]
		\centering
		\subfloat[$\epsilon=0.5.$]{\includegraphics[width=0.5\textwidth]{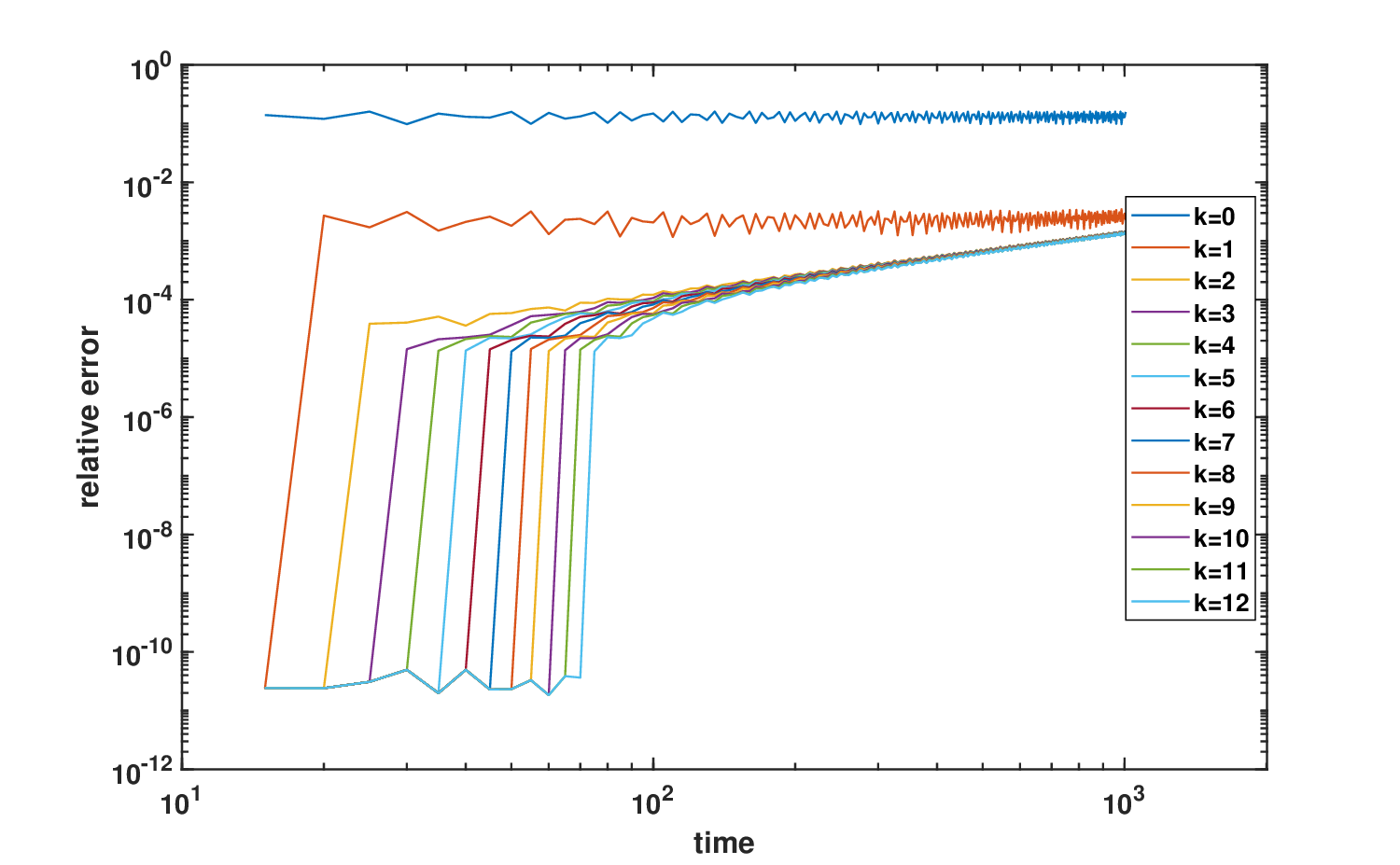}}
		\subfloat[$\epsilon=0.1.$]{\includegraphics[width=0.5\textwidth]{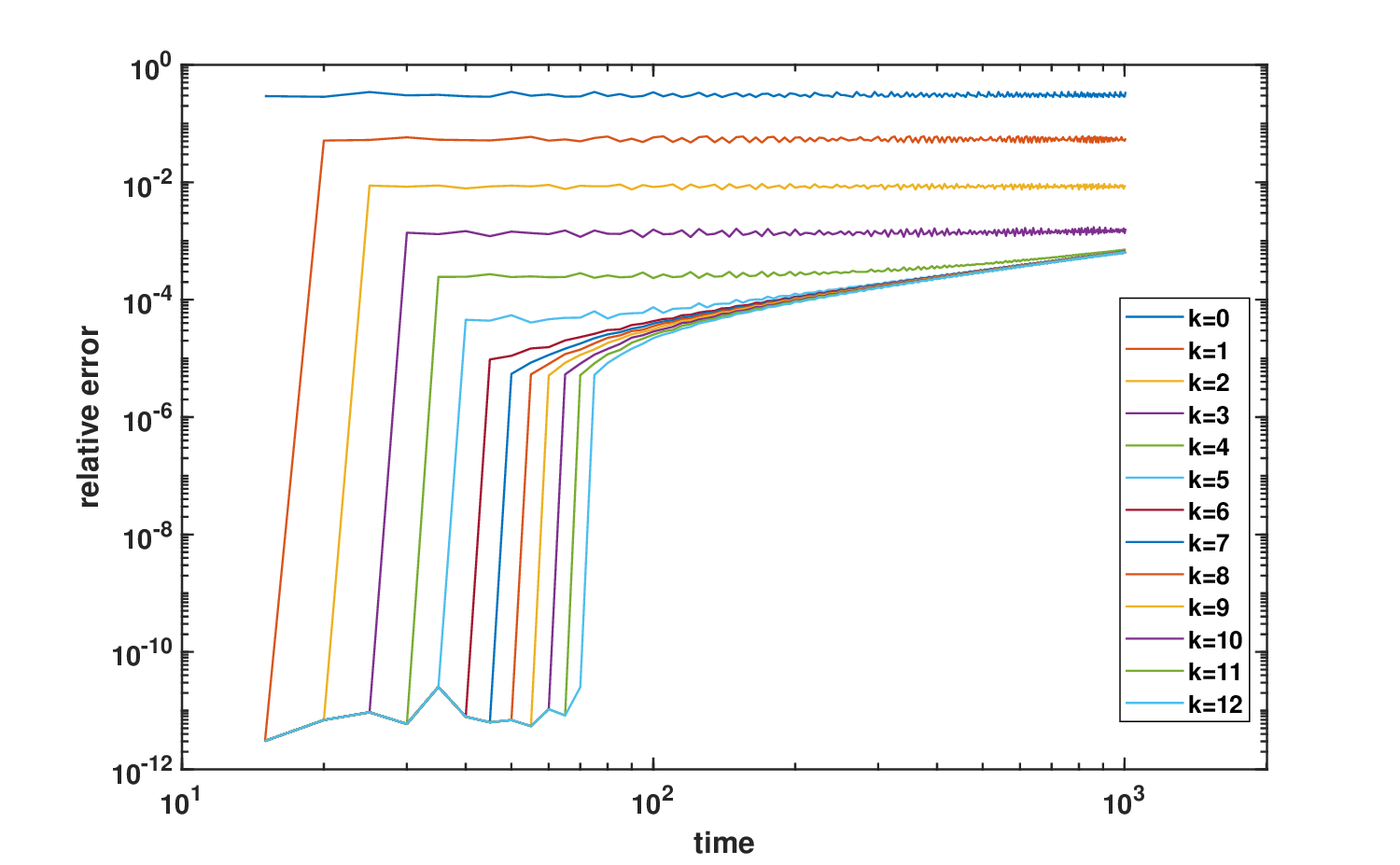}}
		\caption{The evolution curves of the relative error of $U^k_n$ obtained by the method "Parareal-I" in each parareal iteration for solving the ABC flow with $\epsilon=0.5,0.1$}
		\label{ABCFEM-POD}
	\end{figure}

		We first show the results of the method "Parareal-I". Here, we also set $T_0=10.0$. Other parameters are the same as those in the method "AdapParareal". The numerical results are shown in Fig. \ref{ABCFEM-POD}. The dimensions of POD subspaces for the coarse propagator in the $0$-th iteration for these two cases are $30$ and $46$ respectively.  From the two sub-figures in Fig. \ref{ABCFEM-POD}, we can also see that the relative error can reach $O(10^{-3})$ after several iterations, which   indicates the effectiveness of the model order reduction based parareal method. However, similar to the first example, the accuracy cannot be improved any more as the increase of iteration.
	
	\FloatBarrier
			
\begin{figure}[H]
		\centering
		\subfloat[$\epsilon=0.5.$]{\includegraphics[width=0.45\textwidth]{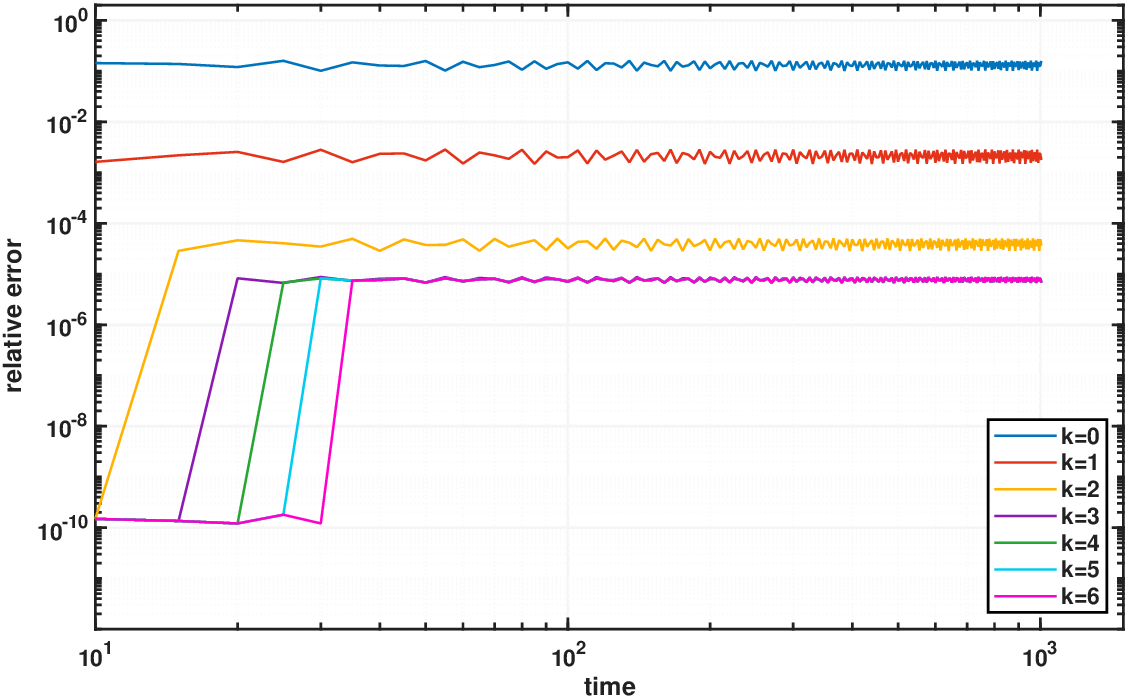}}\hspace{1em}
		\subfloat[$\epsilon=0.1.$]{\includegraphics[width=0.45\textwidth]{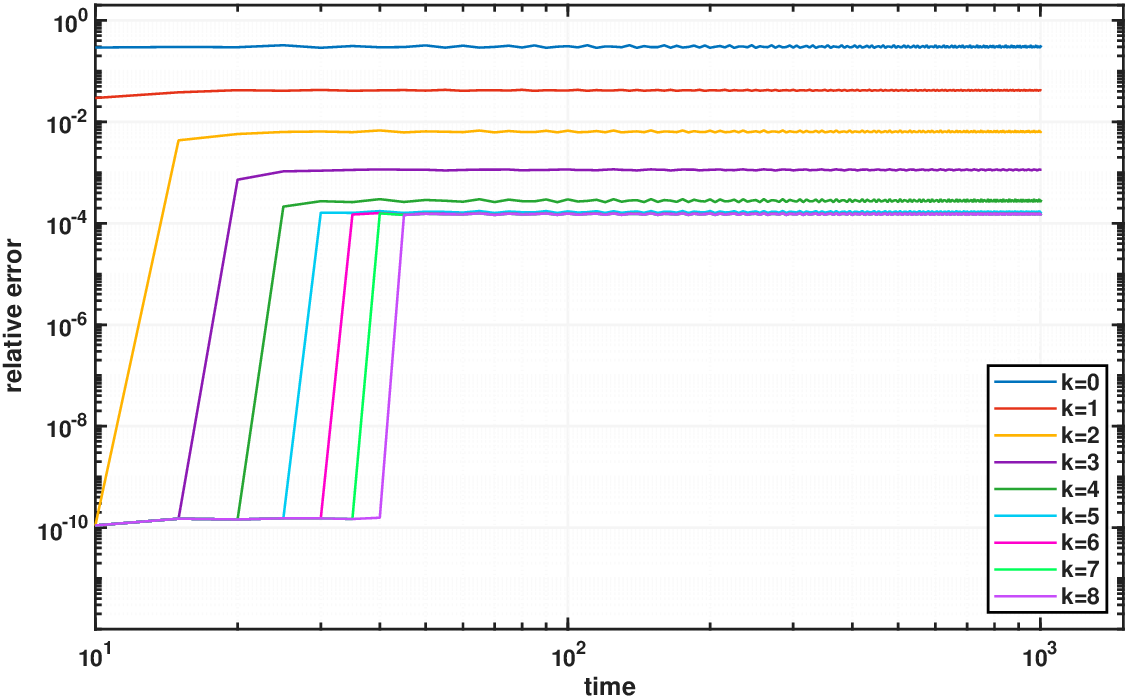}}
		\hspace{1em}
		\caption{The evolution curves of the relative error of $U^k_n$ obtained by the method "Parareal-II" in each parareal iteration for solving the ABC flow with different $\epsilon$.}
		\label{ABCFEM-FEM}
	\end{figure}

	We then show the results obtained by the method "Parareal-II", which uses the finite element with coarse mesh size to do the spatial discretization in the coarse propagator. Similar to the Kolmogorov flow, the coarse mesh is obtained by refining the initial mesh $18$ times, and the number of degrees of freedom  is  $262144$. Other parameters are the same as those used in the method "AdapParareal". The numerical results are shown in Fig.~\ref{ABCFEM-FEM}.  From the sub-figures in Fig.~\ref{ABCFEM-FEM}, we can see that the relative error can reach  $O(10^{-4})$  after several iterations, which  indicates   the method "Parareal-II" can achieve better accuracy than the method "Parareal-I". However, similar to the first example, the accuracy cannot be improved any more even if more parareal iterations are taken.
	
	\FloatBarrier

		We also list the computational cost of the method "Parareal-II" in Table \ref{tableABCFEM_time}(unit:s). The speedup of this method is not satisfactory.
	
	\begin{table}[!htb]
		\caption{ The time cost of the method "Parareal-II" per iteration for solving ABC flow  with different $\epsilon$ (unit:s).}\centering
		\label{tableABCFEM_time}
		\begin{tabular}{|cccc|}
			\hline
			$\epsilon$ & $T_{\mathcal{F}}$   & $T_{\mathcal{G}}$ & $T_{\text{total}}$\\
			\hline
			0.5         & 62.13& 0.92     &246.13 \\
			\hline
			0.1         & 61.89   & 0.88          &237.89 \\
			\hline
		\end{tabular}
	\end{table}

\begin{figure}[H]
		\centering
		\subfloat[$\epsilon=0.5.$]{\includegraphics[width=0.5\textwidth]{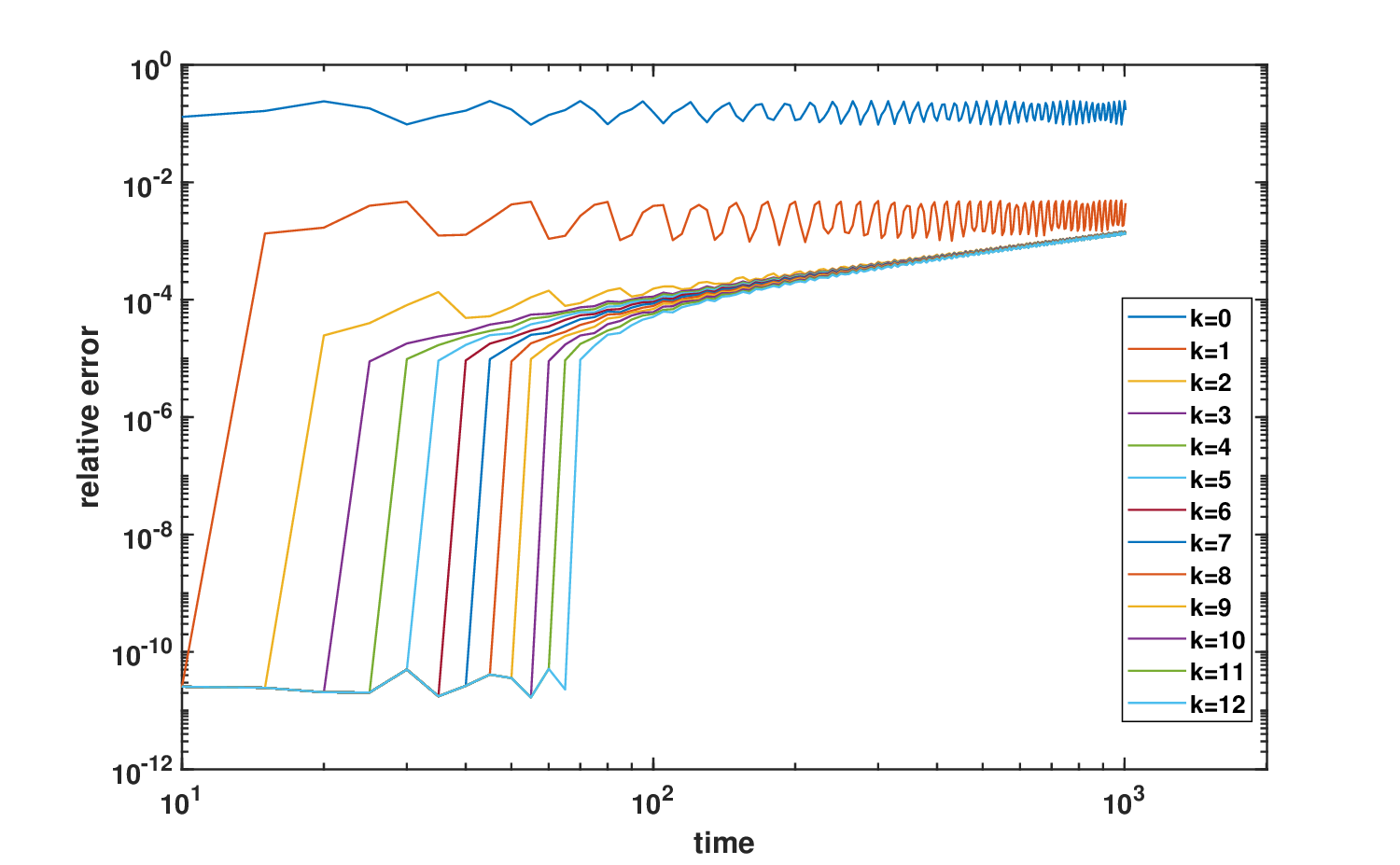}}
		\subfloat[$\epsilon=0.1.$]{\includegraphics[width=0.5\textwidth]{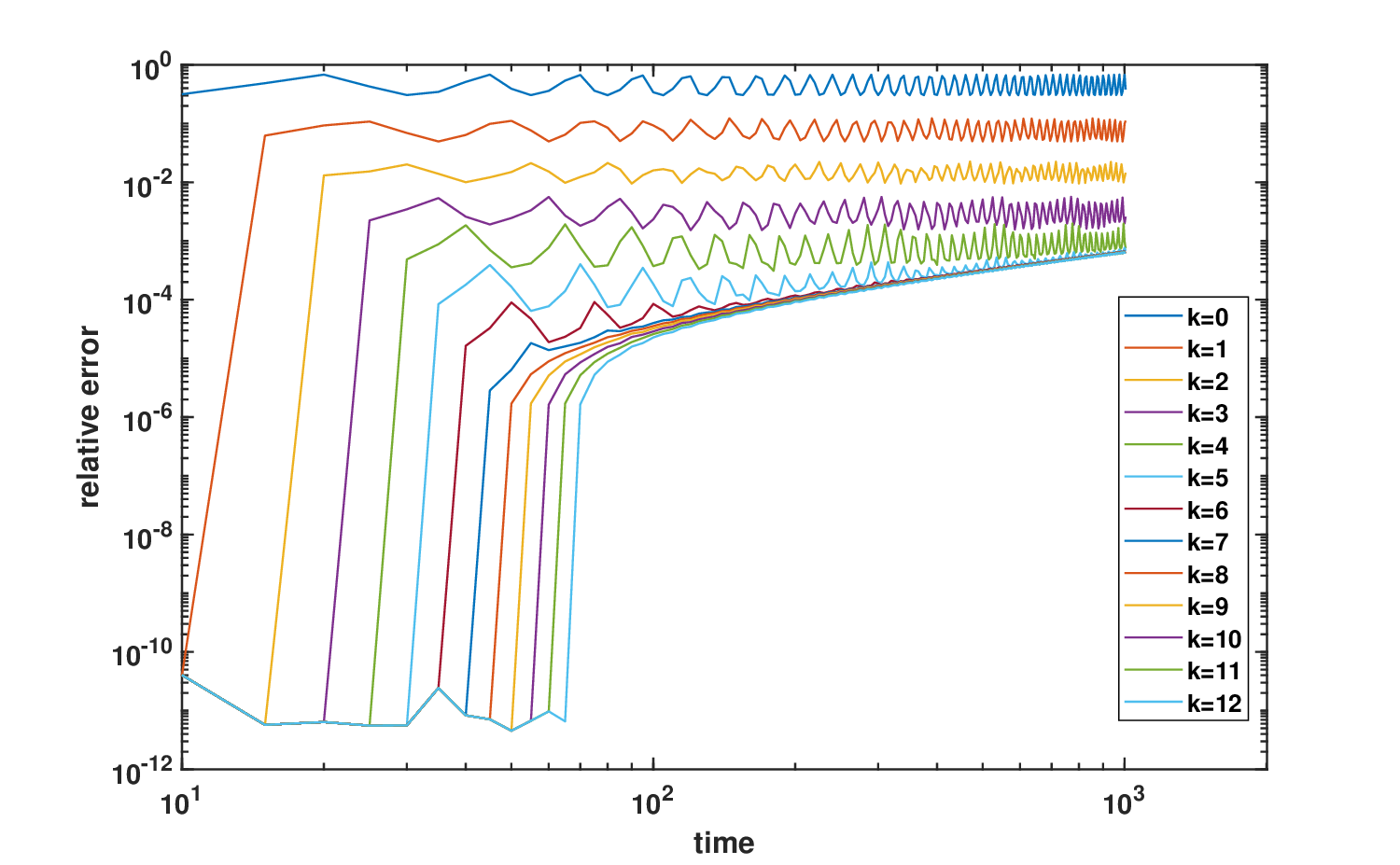}}
		\caption{The evolution curves of the relative error of $U^k_n$ obtained by the method "AdapParareal" with $(m_l,p)=(0,0)$ in each parareal iteration for solving the ABC flow  with $\epsilon=0.5,0.1$}
		\label{ABCFEM-PODONE}
	\end{figure}

\begin{figure}[!tbhp]
		\centering
		\begin{minipage}{\textwidth}
			\centering
			\subfloat[$\epsilon=0.5.$]{\includegraphics[width=0.5\textwidth]{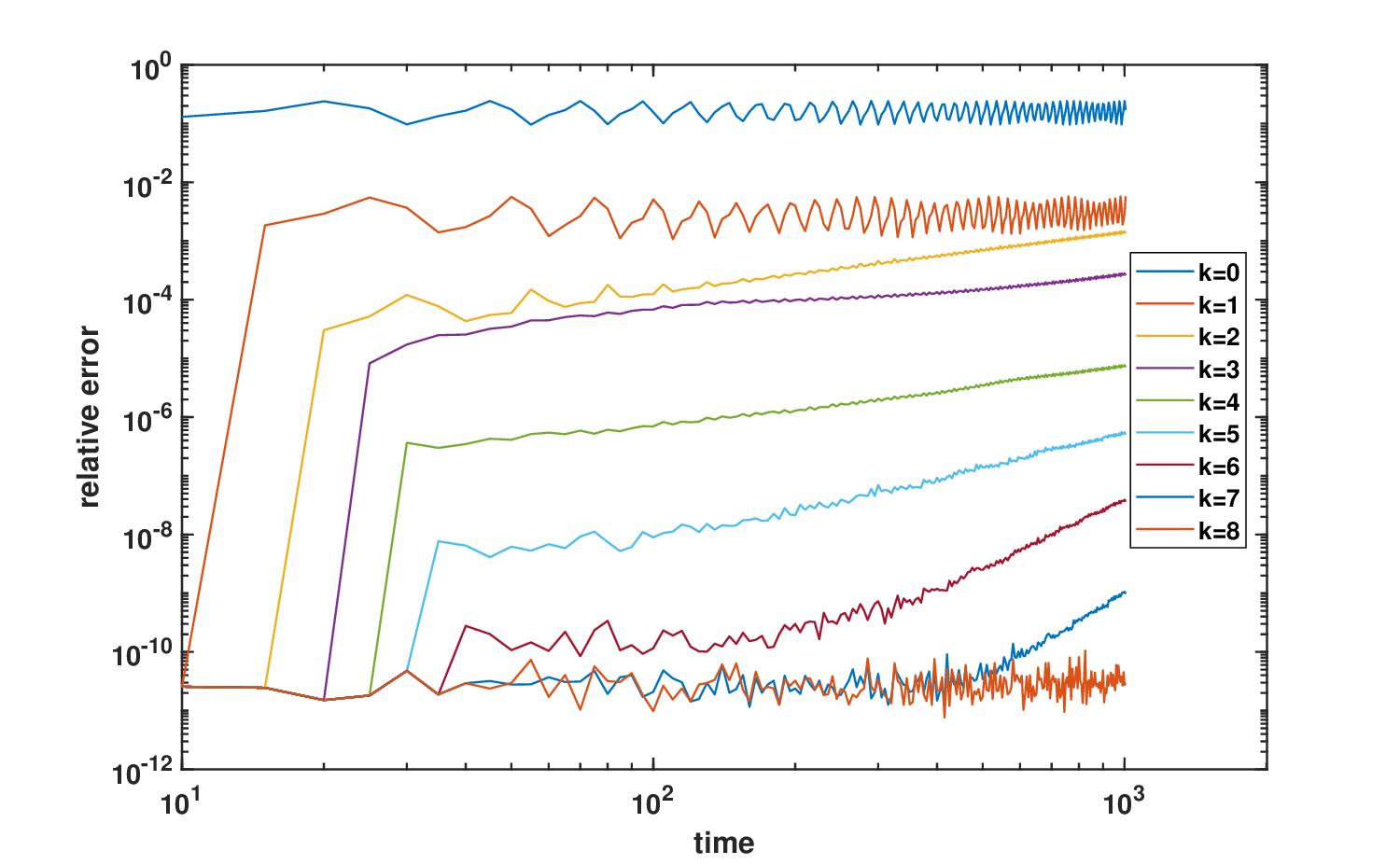}}
			\subfloat[$\epsilon=0.1.$]{\includegraphics[width=0.5\textwidth]{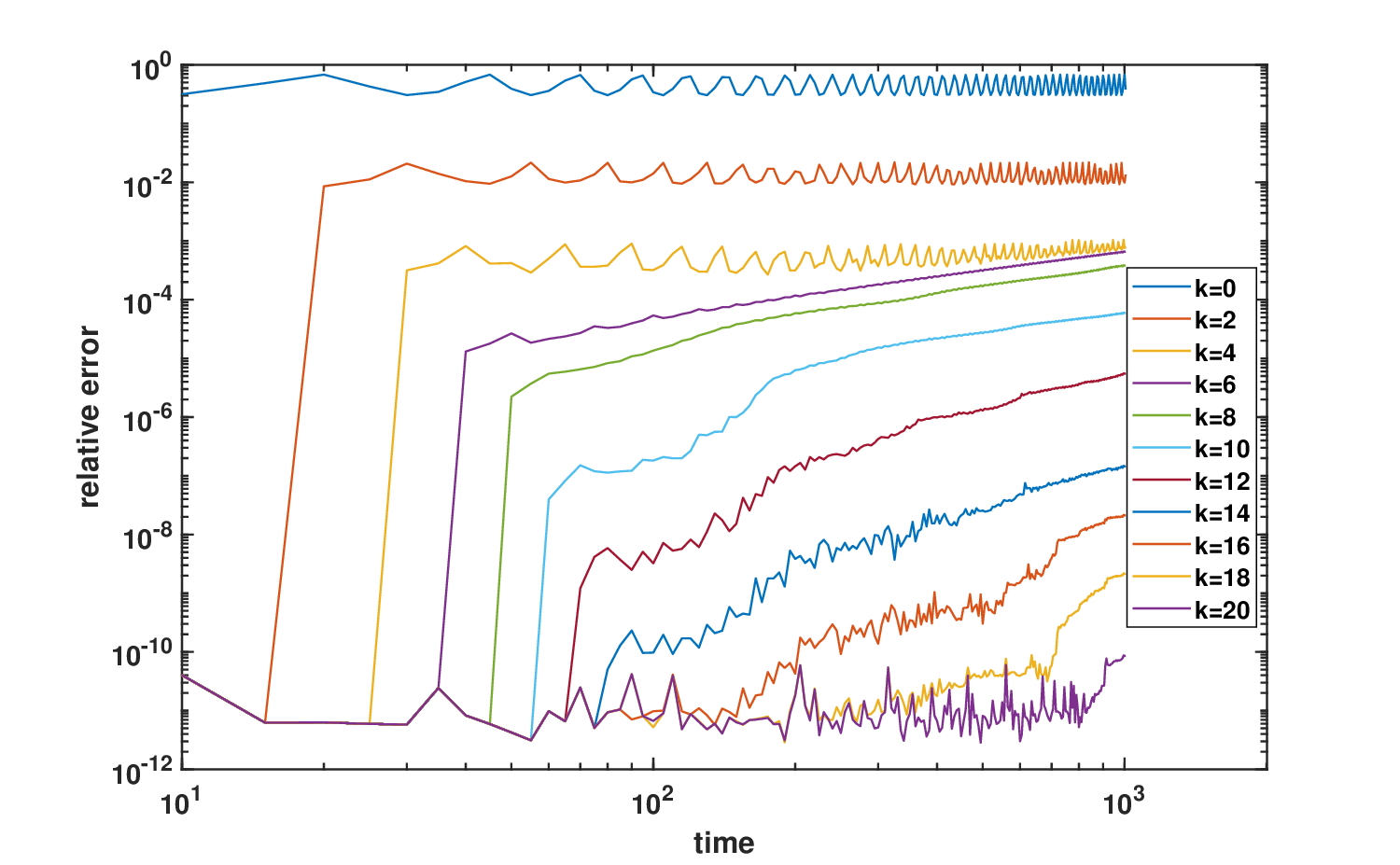}}
		\end{minipage}
		\caption{The evolution curves of the relative error of $U^k_n$ obtained by the method "AdapParareal"  with $(m_l,p)=(1,0)$ in each parareal iteration for solving the ABC flow with $\epsilon=0.5,0.1$}
		\label{ABCFEM-APODlast}
	\end{figure}

\begin{figure}[!tbhp]
		\centering
		\subfloat[$\epsilon=0.5$]{\includegraphics[width=0.5\textwidth]{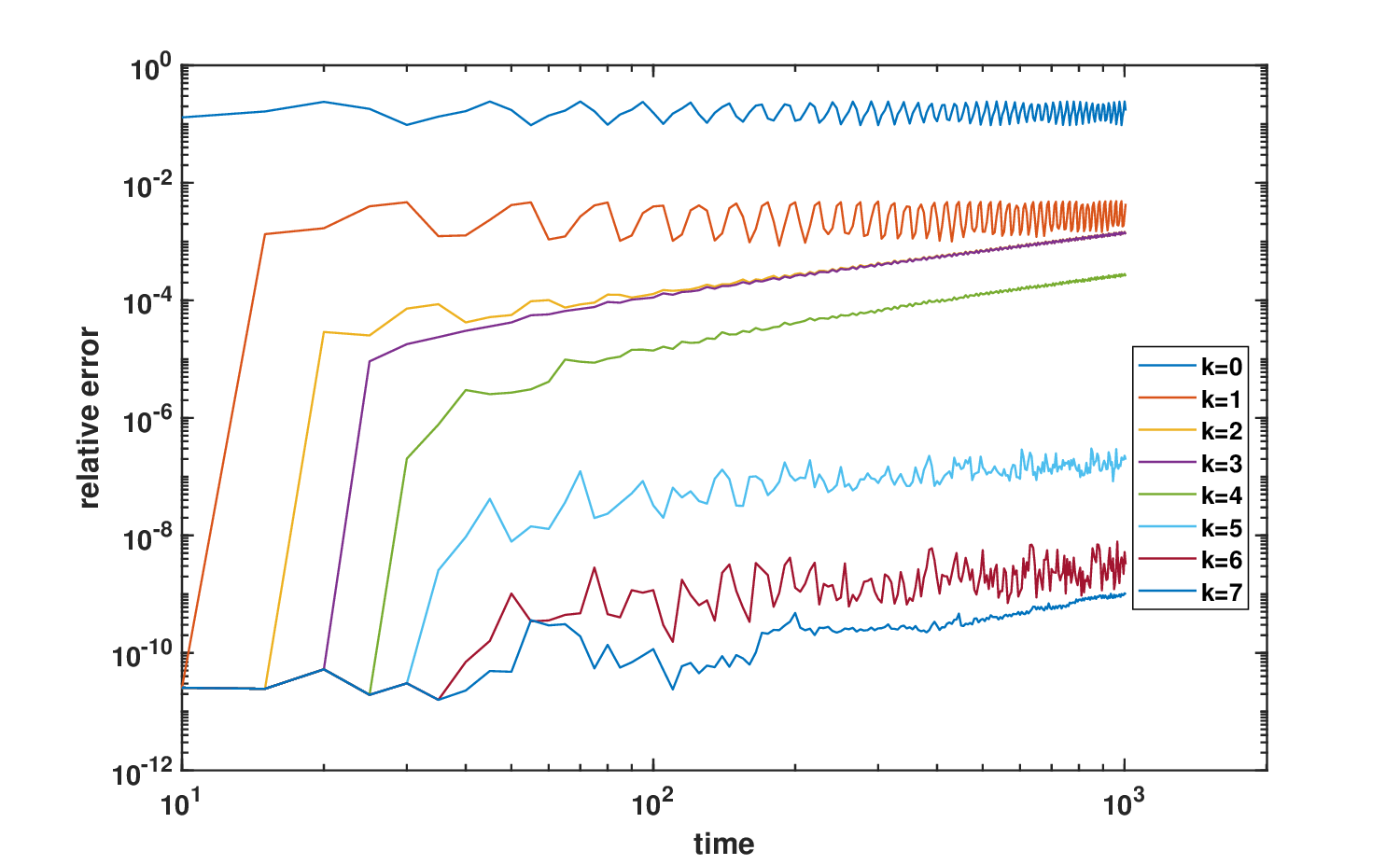}}
		\subfloat[$\epsilon=0.1$]{\includegraphics[width=0.5\textwidth]{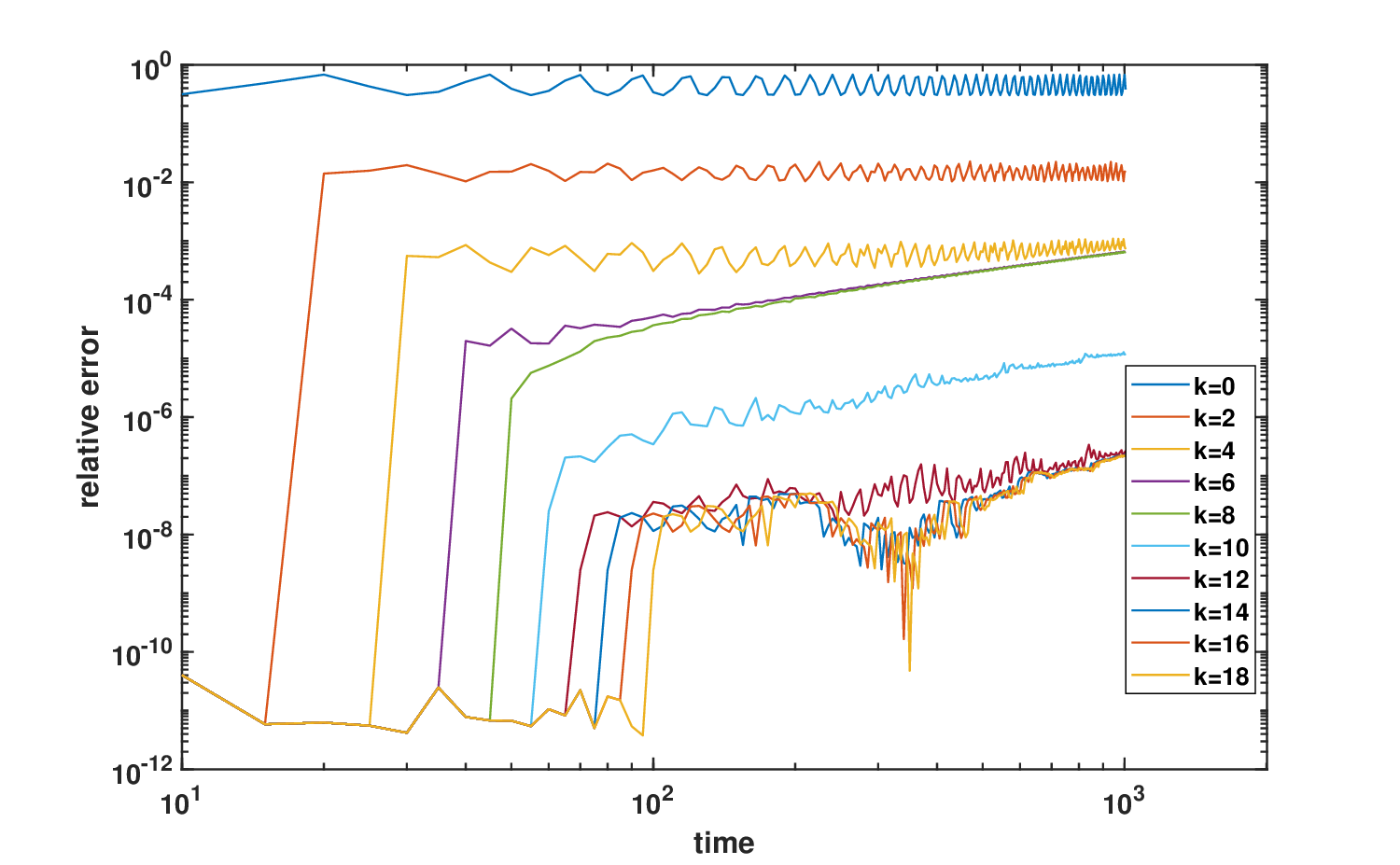}}
		\caption{The evolution curves of the relative error of $U^k_n$ obtained by the method "AdapParareal" with $(m_l,p)=(0,1)$ in each parareal iteration for solving the ABC flow with $\epsilon=0.5,0.1$}
		\label{ABCFEM-APODtwo}
	\end{figure}

\begin{figure}[!tbhp]
		\centering
		\subfloat[$\epsilon=0.5.$]{\includegraphics[width=0.5\textwidth]{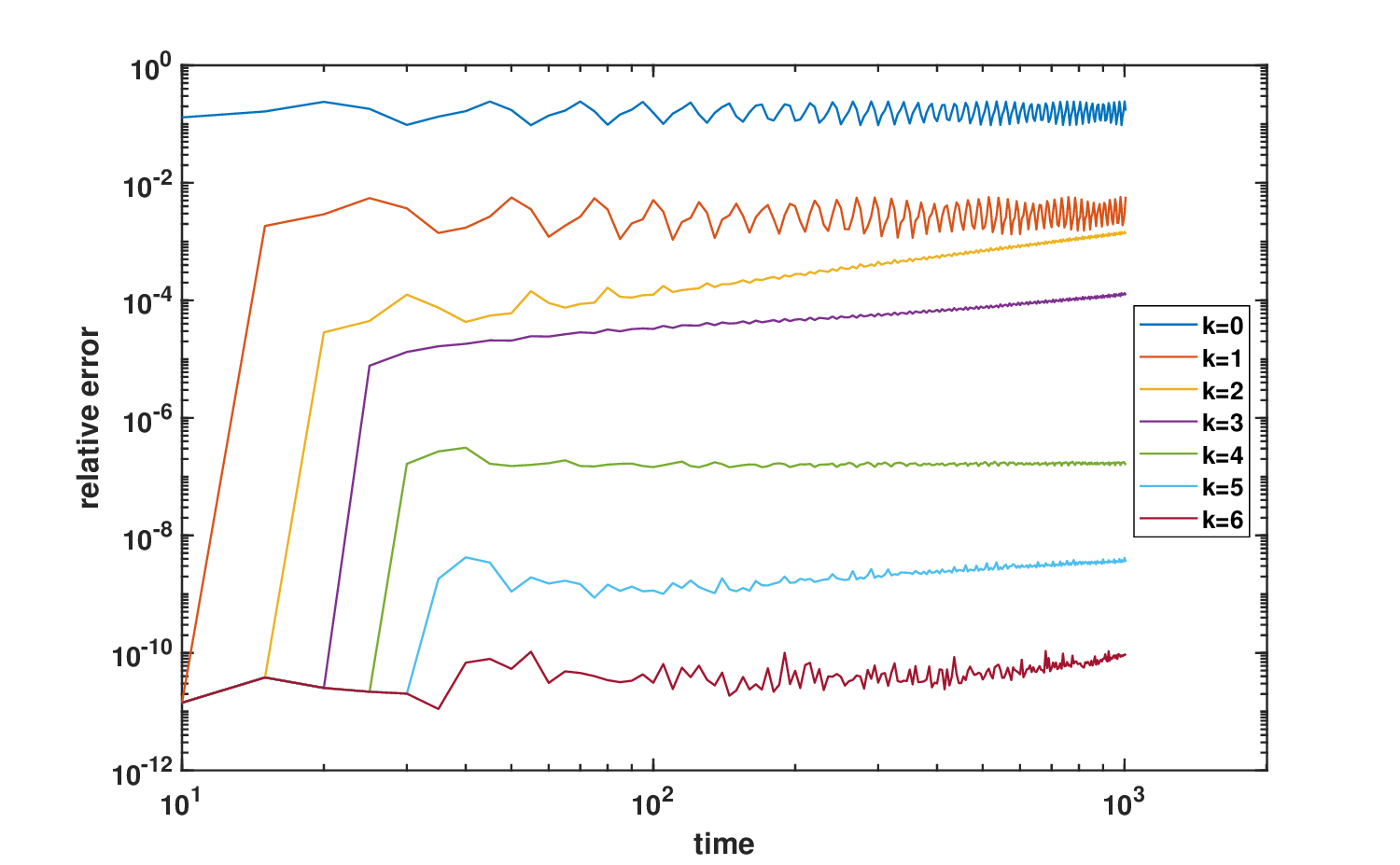}}
		\subfloat[$\epsilon=0.1.$]{\includegraphics[width=0.5\textwidth]{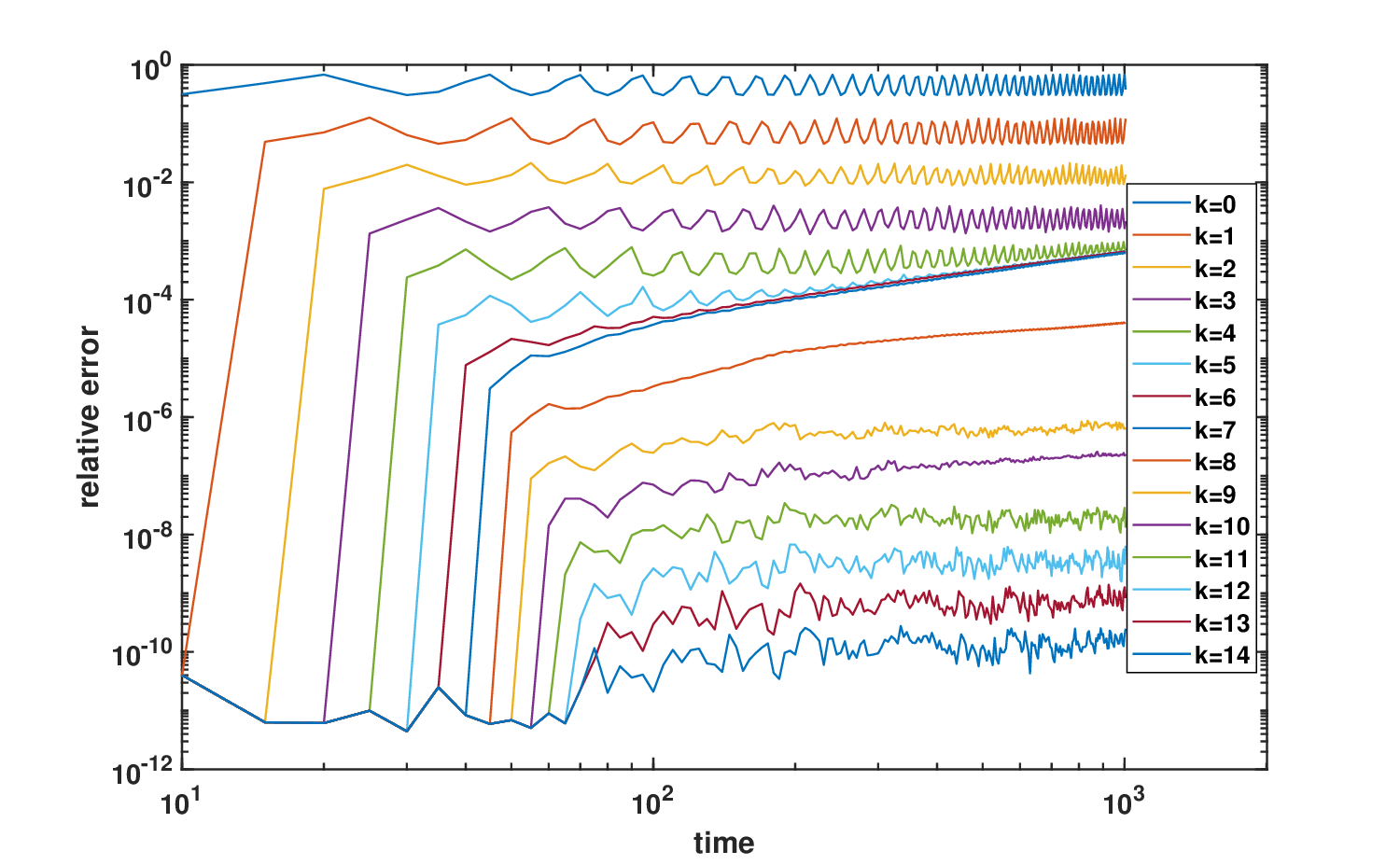}}
		\caption{The evolution curves of the relative error of $U^k_n$ obtained by the method "AdapParareal"  with $(m_l,p)=(1,1)$ in each parareal iteration for solving the ABC flow with $\epsilon=0.5,0.1$}
		\label{ABCFEM-APODTHREELONG}
	\end{figure}

	Then, we also use these two cases with $\epsilon=0.5$ and $0.1$ to show the performance of our new method "AdapParareal" with different $(m_l,p)$. The dimensions of POD subspaces obtained for the coarse propagator in the $0$-th iteration are $18$ and $24$ respectively. The corresponding results are presented in Fig. \ref{ABCFEM-PODONE}-\ref{ABCFEM-APODTHREELONG}. Similar to the example of the Kolmogorov flow, we can see that the relative error can reach at least $O(10^{-3})$ after several iterations from these figures, which also shows that the method "AdapParareal" is effective. Besides, when comparing them with the method "Parareal-I" and "Parareal-II" in Fig. \ref{ABCFEM-POD}-\ref{ABCFEM-FEM}, the numerical results show that higher accuracy can be obtained by the method "AdapParareal". Similarly, the relative error does not increase rapidly over time, which means that our adaptive parareal method has good performance for simulating long-term evolution.

	\FloatBarrier
	
	\FloatBarrier
	
	\FloatBarrier
	
	Similar to the first example, by comparing the numerical results using different $(m_l,p)$ in Fig. \ref{ABCFEM-PODONE}-\ref{ABCFEM-APODTHREELONG}, we also observe that as the number of time subintervals used in constructing the coarse propagator increases, the convergence becomes faster. Meanwhile, the accuracy can also be improved to $O(10^{-10})$ for the case of $(m_l, p) = (1,0)$ and $ (m_l, p) = (1, 1)$. We have done some more tests for using larger $m_l$ or $p$, and are not able to  see obvious improvement for the accuracy. Note that selecting more time subintervals means more computational cost.  Therefore, considering the accuracy and the cost, we recommend to set $(m_l, p) = (1, 1)$.
	
\begin{figure}
		\centering
		\subfloat[$\epsilon=0.05.$]{\includegraphics[width=0.5\textwidth]{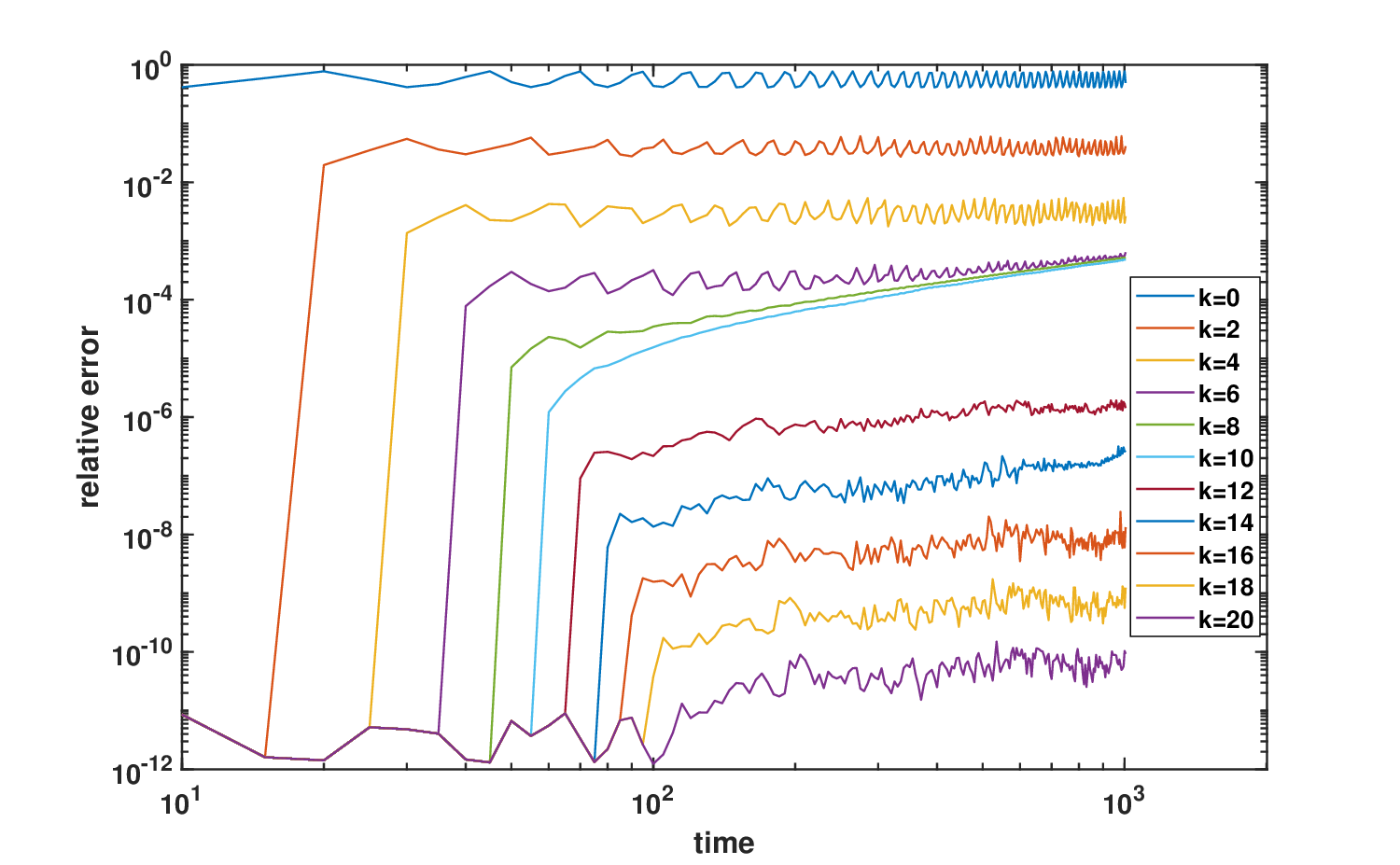}}
		\subfloat[$\epsilon=0.01.$]{\includegraphics[width=0.5\textwidth]{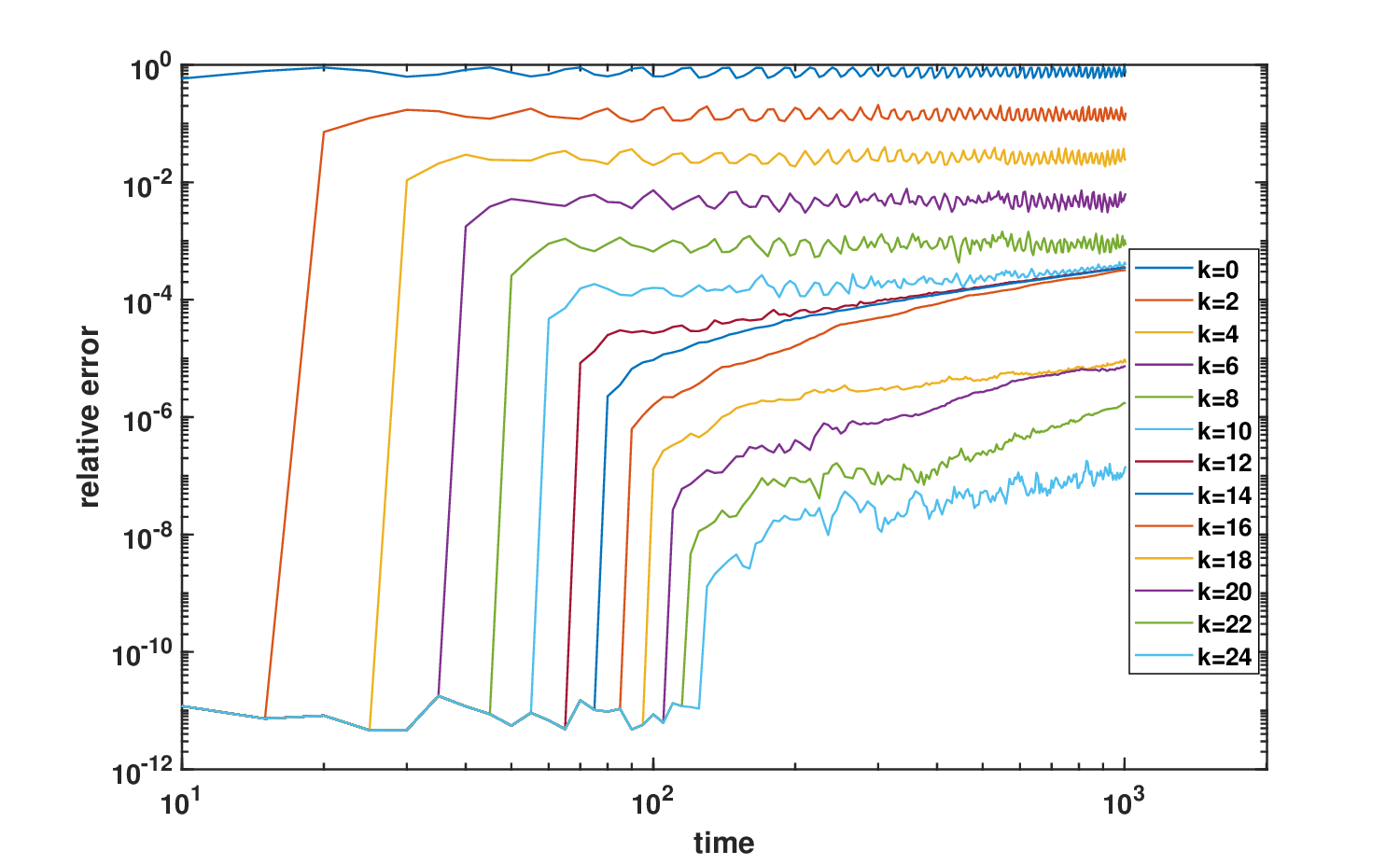}}
		
		\subfloat[$\epsilon=0.005.$]{\includegraphics[width=0.5\textwidth]{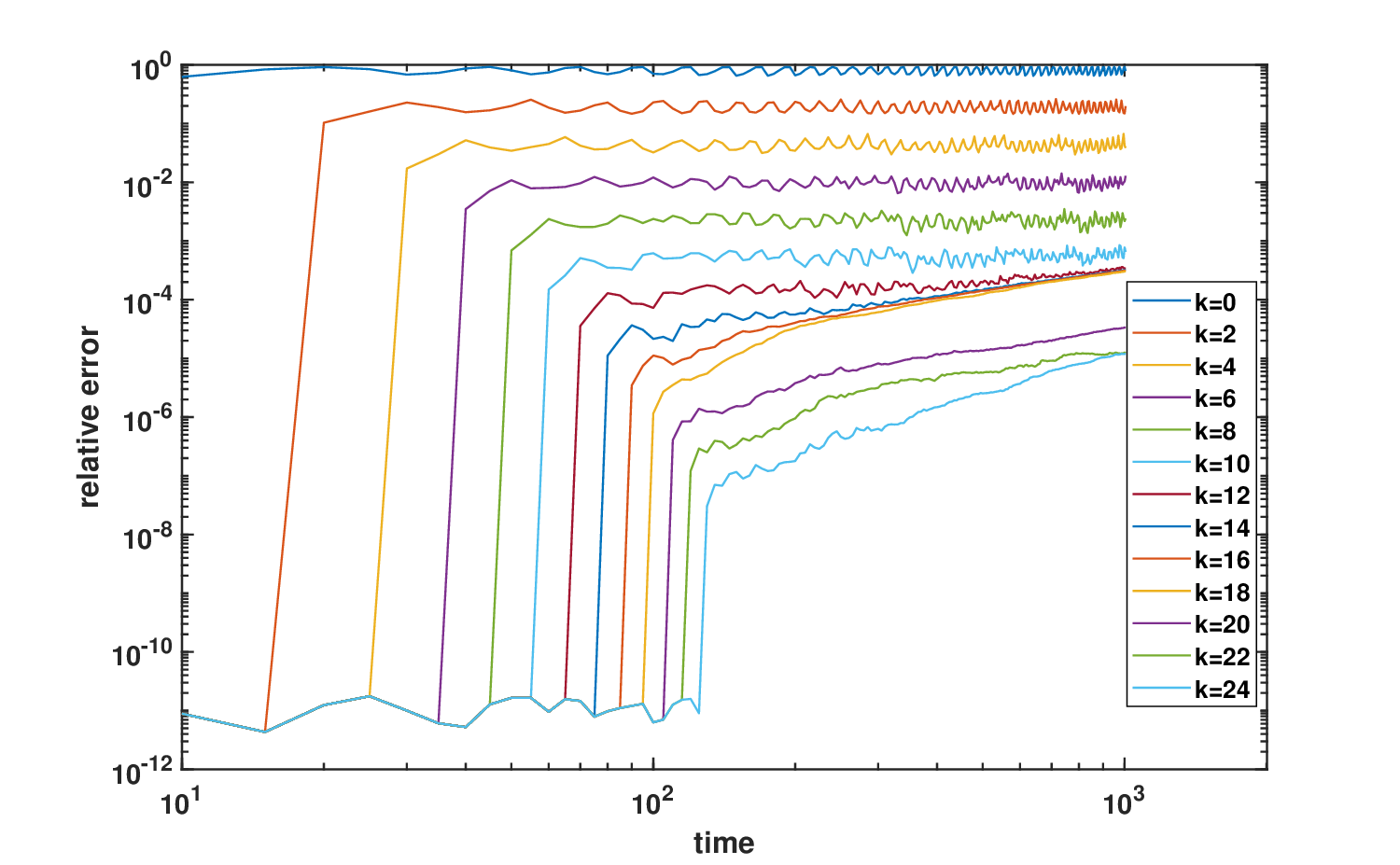}}
		\caption{The evolution curves of the relative error of $U^k_n$ obtained by the method "AdapParareal" with $(m_l,p)=(1,1)$ in each parareal iteration for solving the ABC flow  with $\epsilon=0.05, 0.01, 0.005$}
		\label{ABCFEM-APODTHREE2}
	\end{figure}

	When we simulate the advection-dominated cases with $\epsilon=0.05, 0.01,0.005$ by the method "AdapParareal", we   choose $(m_l, p) = (1, 1)$. The numerical results are presented in Fig. \ref{ABCFEM-APODTHREE2}. The dimensions of POD subspaces for the initial coarse propagator in the $0$-th iteration are $28$, $40$ and $44$ respectively.
	
\FloatBarrier
Similar to the case of the Kolmogorov flow, we can also see from Fig.~\ref{ABCFEM-APODTHREE2} that, the error decreases obviously as the parareal iteration continues. However, compared with the cases $\epsilon = 0.5$ and $\epsilon = 0.1$, the convergence of the adaptive parareal method for these three cases is slower. Fortunately,  the accuracy obtained by our adaptive parareal method can still be very high in long-term evolution after several iterations. We can also see that even for the case of $\epsilon=0.005$, after $24$ parareal iterations, the error obtained on the whole time interval can be smaller than $2\times10^{-5}$.  
	
\begin{figure}[!tbhp]
		\centering
		\subfloat[$\epsilon=0.05.$]{\includegraphics[width=0.5\textwidth]{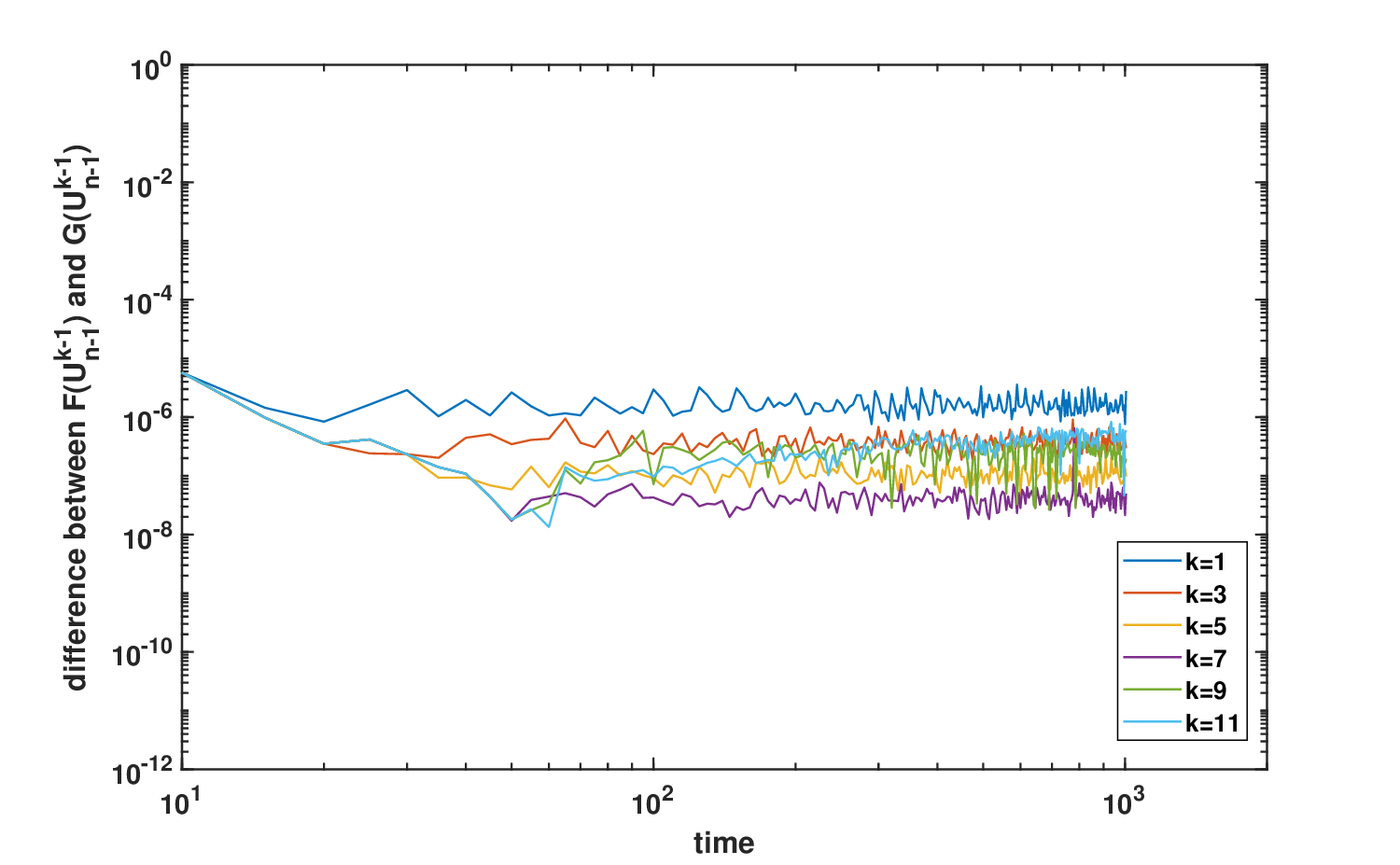}}
		\subfloat[$\epsilon=0.01.$]{\includegraphics[width=0.5\textwidth]{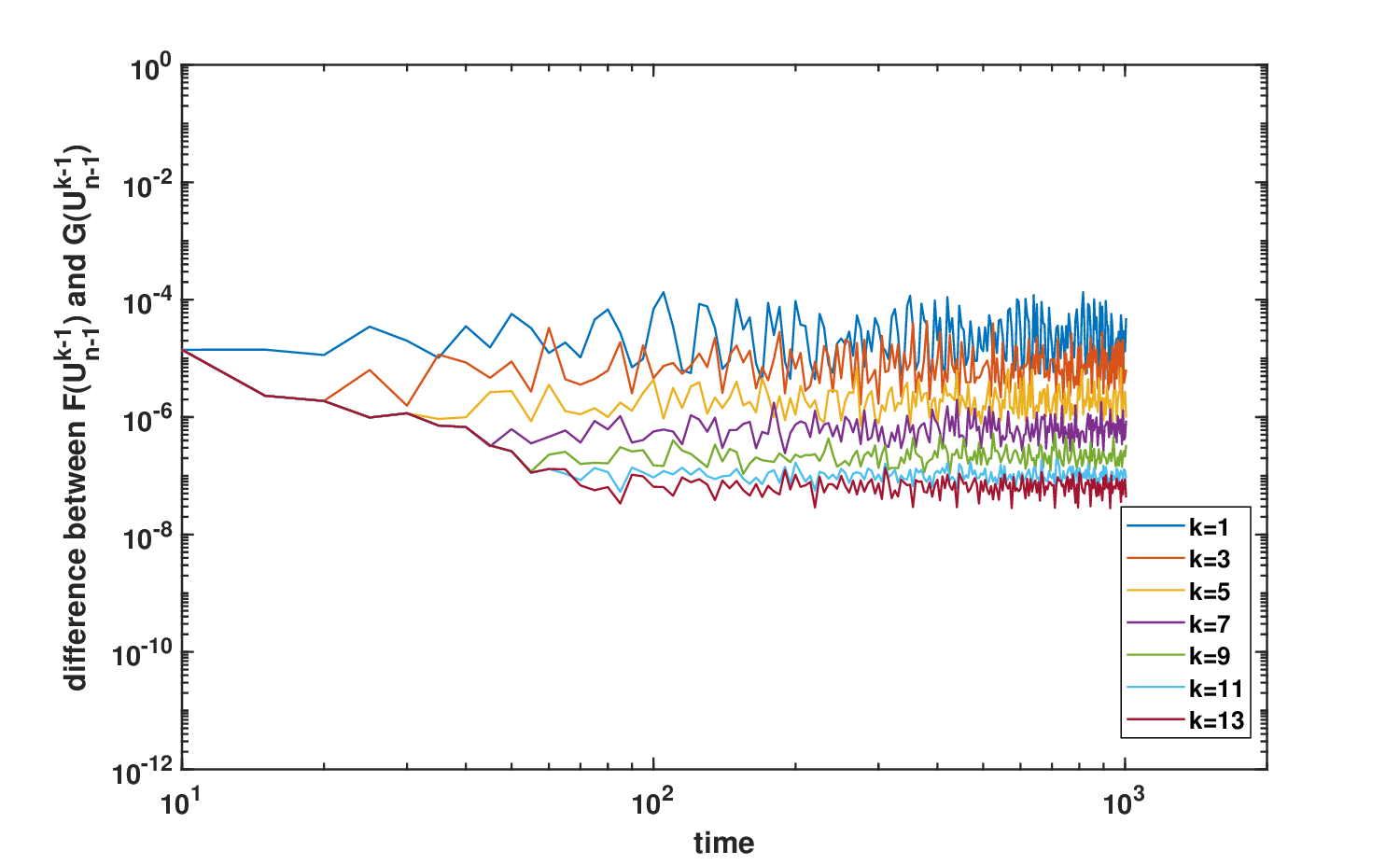}}
		
		\subfloat[$\epsilon=0.005.$]{\includegraphics[width=0.5\textwidth]{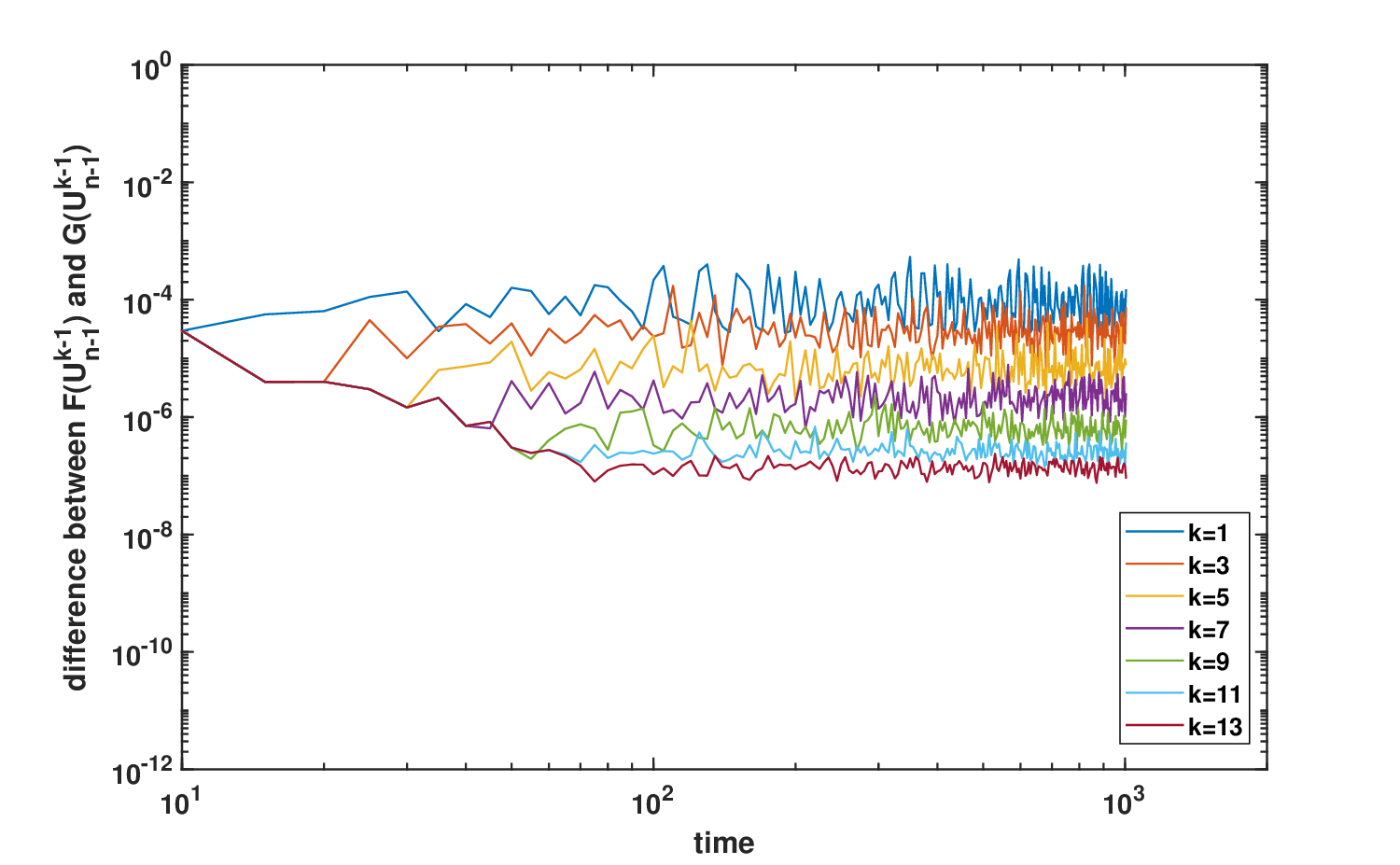}}
		\caption{The evolution curves of $\frac{||\mathcal{F}\left(U_{n-1}^{k-1}\right)-\mathcal{G}_k\left(U_{n-1}^{k-1}\right)||_{L^2}}{||\mathcal{F}\left(U_{n-1}^{k-1}\right)||_{L^2}}$ obtained by the method "AdapParareal" with $(m_l,p)=(1,1)$ in each parareal iteration for solving the ABC flow  with $\epsilon=0.05,0.01,0.005$}
		\label{ABCFEM-APODextra1}
	\end{figure}

\begin{figure}[!tbhp]
		\centering
		\subfloat[$\epsilon=0.05.$]{\includegraphics[width=0.5\textwidth]{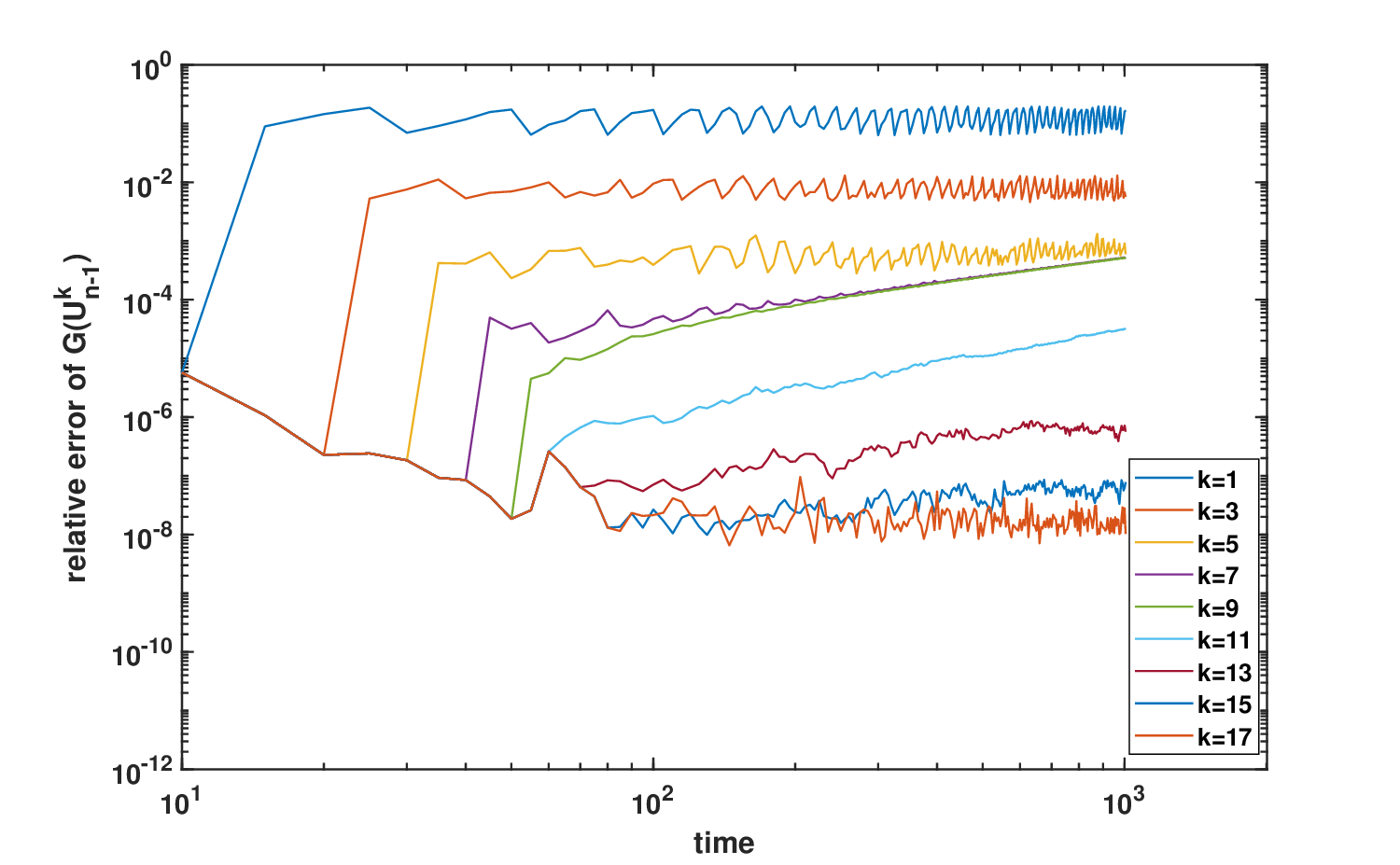}}
		\subfloat[$\epsilon=0.01.$]{\includegraphics[width=0.5\textwidth]{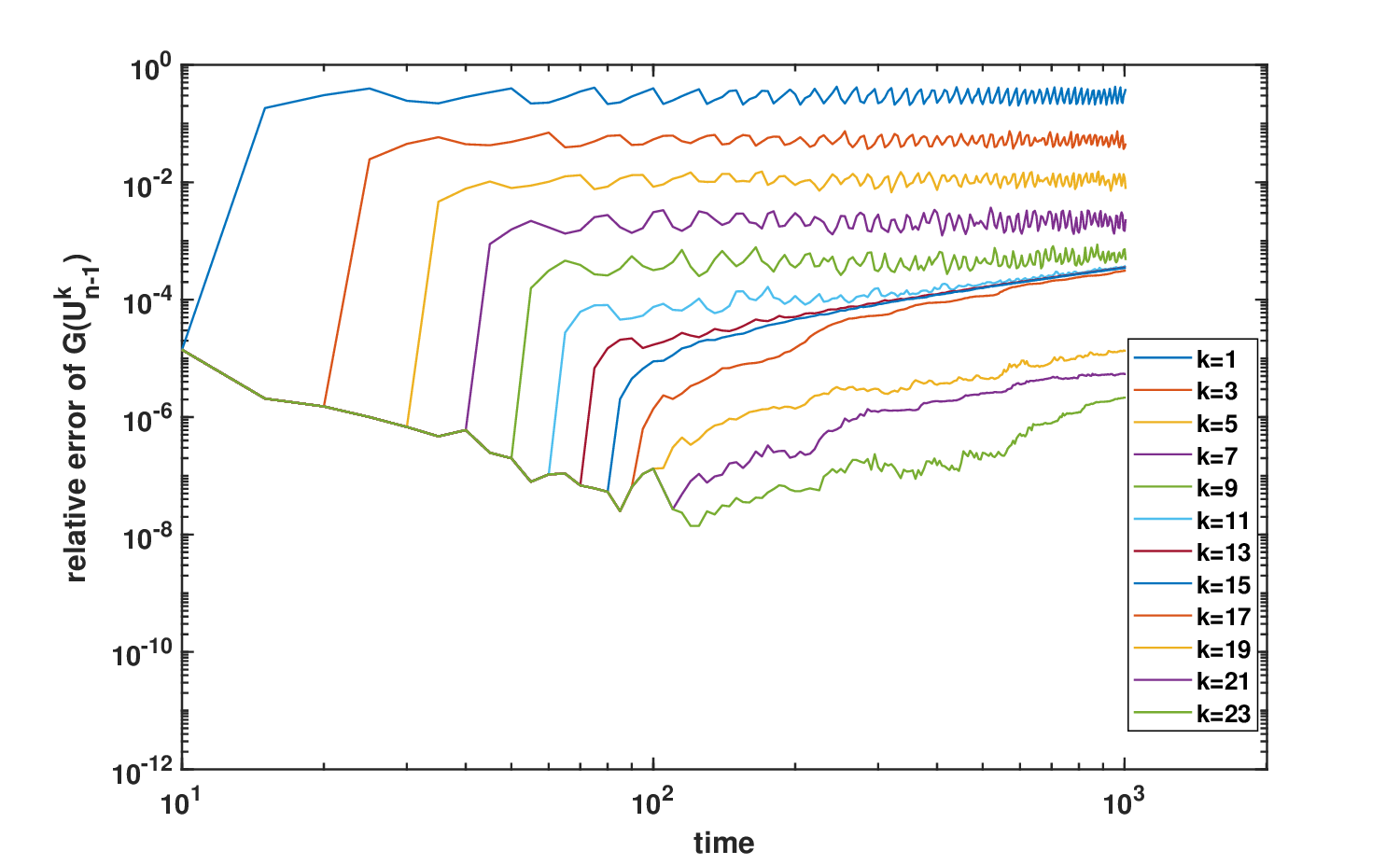}}
		
		\subfloat[$\epsilon=0.005.$]{\includegraphics[width=0.5\textwidth]{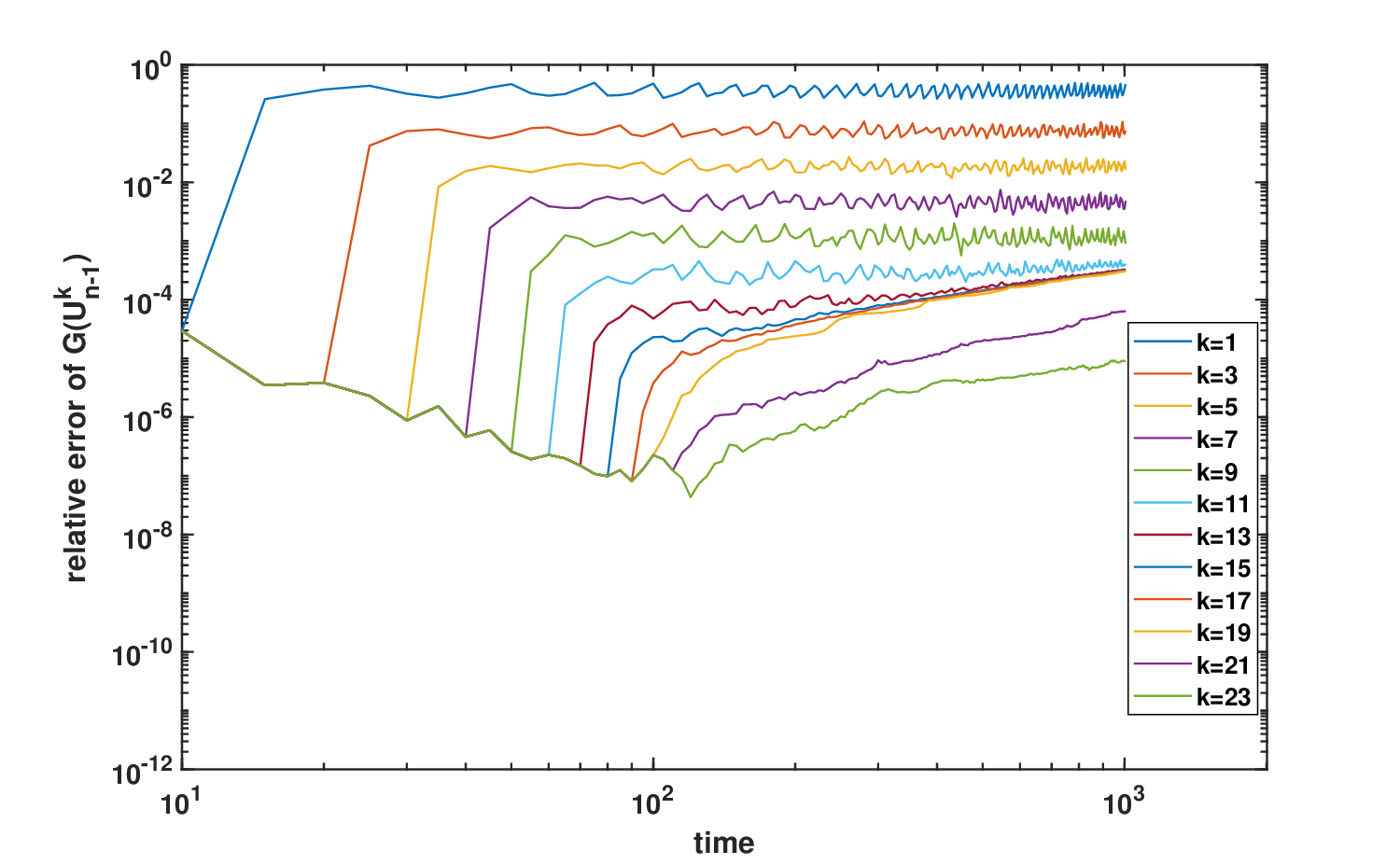}}
		\caption{The evolution curves of $\frac{||U_{n}-\mathcal{G}_k\left(U_{n-1}^{k}\right)||_{L^2}}{||U_{n}||_{L^2}}$ obtained by the method "AdapParareal" with $(m_l,p)=(1,1)$  in each parareal iteration for solving the ABC flow with $\epsilon=0.05,0.01,0.005$}
		\label{ABCFEM-APODextra}
	\end{figure}

	For the ABC flow, we also show the computational cost for each part of our adaptive parareal method, see  Table~\ref{tableABC_time} for the details. Besides, we also show the dimension $m_{\max}$  of the POD subspace for various $\epsilon$  in Table~\ref{tableABC_PODnum}. From  Table~\ref{tableABC_time}, we can see that  the time cost for constructing the  POD basis  is very cheap compared  with the time cost for the fine propagator. As shown in Fig. \ref{ABCFEM-APODTHREELONG}-\ref{ABCFEM-APODTHREE2}, our method achieves an accuracy of $O(10^{-6})$  using  less than $12$  iterations for cases of  $\epsilon=0.5,0.1,0.05$, and achieves an accuracy of  $O(10^{-6})$ within about $22$  iterations for $\epsilon=0.01$,  and achieves $O(10^{-5})$ accuracy within about $22$  iterations for $\epsilon=0.005$. Note that the number of time subintervals is also $N=200$. Therefore, based on our analysis for the speedup in  Section~\ref{sec: Complex}, we can estimate that our method can achieve a speedup of about $8.7$ for cases of $\epsilon=0.5,0.1,0.05$, about $4$ for cases of $\epsilon=0.01$ and $\epsilon=0.005$.
	\begin{table}[H]
		\caption{The time cost of the method "AdapParareal" with $(m_l,p)=(1,1)$ in each iteration for solving ABC flow  with different $\epsilon$  (unit:s).}\centering
		\label{tableABC_time}
		\begin{tabular}{|ccccccc|}
			\hline
			$\epsilon$ & $T_{\mathcal{F}}$          & $T_{\mathcal{G}}$ & $T_U$  & $T_C$ & $T_A$ & $T_\text{total}$  \\
			\hline
			0.5         & 61.36&0.05         & 0.69 & 0.39    & 0.12 & 96.44 \\
			\hline
			0.1         & 62.57  & 0.09        &1.58 & 0.91   & 0.14 & 111.06 \\
			\hline
			0.05          & 56.06 & 0.10      &2.13    & 1.42 & 0.14  & 107.61 \\
			\hline
			0.01           & 61.88 & 0.18     & 5.01   & 3.28  & 0.18  & 142.17  \\
			\hline
			0.005         & 61.52 &  0.18   &    6.24   &4.23 &  0.18   &  143.99 \\
			\hline
		\end{tabular}
	\end{table}
	\begin{table}[H]
		\caption{ The maximum dimension of the POD subspaces $m_\text{max}$ for the method "AdapParareal" with $(m_l,p)=(1,1)$ for solving ABC flow  with different $\epsilon$.}\centering
		\label{tableABC_PODnum}
		\begin{tabular}{|cccccc|}
			\hline
			$\epsilon$ & $0.5$          & $0.1$ & $0.05$  & $0.01$ & $0.005$\\
			\hline
			$m_\text{max}$         & $56$& $87$         & $108$& $185$     & $196$  \\
			\hline 
		\end{tabular}
	\end{table}
	\FloatBarrier
	Similar to the first example, we have also done some tests to demonstrate that the spatial discretization accuracy of the coarse propagator will approach that of the fine propagator's spatial discretization progressively as the iteration proceeds. We also set
	$dT = \delta t = 0. 01$ to remove the  impact of the temporal discretization. We show the difference between the approximation obtained with the coarse propagator and the fine propagator starting from the same initial values and the relative error of $\mathcal{G}_k\left(U_{n-1}^{k}\right)$, that is $\frac{||\mathcal{F}\left(U_{n-1}^{k-1}\right)-\mathcal{G}_k\left(U_{n-1}^{k-1}\right)||_{L^2}}{||\mathcal{F}\left(U_{n-1}^{k-1}\right)||_{L^2}}$ and $\frac{||U_{n}-\mathcal{G}_k\left(U_{n-1}^{k}\right)||_{L^2}}{||U_{n}||_{L^2}}$. The numerical results of the ABC flow with $\epsilon=0.05,0.01,0.005$  respectively are presented in Fig.~\ref{ABCFEM-APODextra1}-\ref{ABCFEM-APODextra}. From these figures, we can first see that $\frac{||\mathcal{F}\left(U_{n-1}^{k-1}\right)-\mathcal{G}_k\left(U_{n-1}^{k-1}\right)||_{L^2}}{||\mathcal{F}\left(U_{n-1}^{k-1}\right)||_{L^2}}$ can reach at least $O(10^{-4})$ for all subintervals and iterations, which indicates that the accuracy of the term $\mathcal{G}_k\left(U_{n-1}^{k-1}\right)$ is very high. Moreover, we also observe that as the increase of the parareal iteration, the accuracy of the coarse propagator improves obviously. After several iterations, the accuracy of the coarse propagator becomes very high. This also shows that as the parareal iteration proceeds, the accuracy of spatial discretization of the coarse propagator can indeed be improved. Therefore, similar to the first example, our adaptive parareal method can also be viewed as a parareal based adaptive POD method.
	
	\FloatBarrier	
	
	\section{Concluding remarks}
	\label{sec: Conclusion} 
	
In this paper, we have proposed a model order reduction based adaptive parareal method for time-dependent partial differential equations, by adaptively updating the spatial discretization of the coarse propagator using data obtained by the fine propagator during the parareal iterations together with the POD method.  The construction of the coarse propagator is local in time. That is, we use the data obtained by the fine propagator from  neighboring time subintervals and neighboring iteration steps to construct the POD subspace  for each time subinterval, which makes  the subspaces used  for different time subintervals  different and can capture the local feature of the solution efficiently. Besides, in order to eliminate the projection error,  our method augments each time local reduced subspace with the current parareal interface state. Therefore, our method can produce a coarse propagator whose spatial approximation progressively approaches that of the fine propagator  while retaining a low  computational cost and low memory cost.

	Furthermore, as stated in the end of subsection \ref{Kolmogorov flow}, this method can also be interpreted as a parareal based adaptive POD method. When using the POD method to solve time-dependent partial differential equations, if the initial snapshots fail to capture the solution's behavior over the entire time interval, the numerical errors will increase rapidly over time. Our method can be viewed as an adaptive POD method which uses the data obtained by fine propagator in each time subinterval during the parareal iteration together with POD method to construct the subspace for spatial discretization. 
	
	We have applied our methods to two typical 3D advection-diffusion equations, the Kolmogorov flow and the ABC flow. Numerical results show both the effectiveness and accuracy of our method in simulating long-term evolution, and its ability to obtain significant speedup. Furthermore, we would like to emphasize that our method performs well for the advection-dominated cases. In our future work, we will focus on designing more efficient coarse propagators, improving updating strategies, and applying this method to other types of time-dependent partial differential equations.
	\FloatBarrier
	\section*{Acknowledgments}
	The authors would like to thank Professor Aihui Zhou for his inspiring discussions during the early stages of this research and his valuable suggestions during the writing of this manuscript, and would also like to thank the anonymous reviewers for their constructive suggestions for the improvement of the manuscript.
	
	\section*{CRediT authorship contribution statement}
	
	All the authors contributed equally to this work.
	
	\bibliography{references}

\begin{thebibliography}{10}

\bibitem{acary2010implicit}
Vincent Acary and Bernard Brogliato.
\newblock Implicit {E}uler numerical scheme and chattering-free implementation
  of sliding mode systems.
\newblock {\em Systems \& Control Letters}, 59(5):284--293, 2010.

\bibitem{baffico2002parallel}
Leonardo Baffico, Stephane Bernard, Yvon Maday, Gabriel Turinici, and Gilles
  Z{\'e}rah.
\newblock Parallel-in-time molecular-dynamics simulations.
\newblock {\em Physical Review E}, 66(5):057701, 2002.

\bibitem{bal2005convergence}
Guillaume Bal.
\newblock On the convergence and the stability of the parareal algorithm to
  solve partial differential equations.
\newblock In {\em Domain decomposition methods in science and engineering},
  pages 425--432. Springer, 2005.

\bibitem{Bal2002}
Guillaume Bal and Yvon Maday.
\newblock A ``parareal'' time discretization for non-linear pde's with
  application to the pricing of an american put.
\newblock {\em Notes in Computational Science and Engineering}, 23:189--202,
  2002.

\bibitem{bal2008symplectic}
Guillaume Bal and Qi~Wu.
\newblock {\em Symplectic parareal}.
\newblock Springer, 2008.

\bibitem{benner2015survey}
Peter Benner, Serkan Gugercin, and Karen Willcox.
\newblock A survey of projection-based model reduction methods for parametric
  dynamical systems.
\newblock {\em SIAM review}, 57(4):483--531, 2015.

\bibitem{blouza2011parallel}
Adel Blouza, Laurent Boudin, and Sidi~M Kaber.
\newblock Parallel in time algorithms with reduction methods for solving
  chemical kinetics.
\newblock {\em Communications in Applied Mathematics and Computational
  Science}, 5(2):241--263, 2011.

\bibitem{borue1996numerical}
Vadim Borue and Steven~A Orszag.
\newblock Numerical study of three-dimensional {K}olmogorov flow at high
  {R}eynolds numbers.
\newblock {\em Journal of Fluid Mechanics}, 306:293--323, 1996.

\bibitem{brenner2008mathematical}
Susanne~C Brenner.
\newblock {\em The mathematical theory of finite element methods}.
\newblock Springer, 2008.

\bibitem{brummell2001linear}
Nicholas~H Brummell, Fausto Cattaneo, and Steven~M Tobias.
\newblock Linear and nonlinear dynamo properties of time-dependent {ABC} flows.
\newblock {\em Fluid Dynamics Research}, 28(4):237, 2001.

\bibitem{caldas2021modified}
Joao~G Caldas~Steinstraesser, Vincent Guinot, and Antoine Rousseau.
\newblock Modified parareal method for solving the two-dimensional nonlinear
  shallow water equations using finite volumes.
\newblock {\em The SMAI journal of computational mathematics}, 7:159--184,
  2021.

\bibitem{carlberg2019data}
Kevin Carlberg, Lukas Brencher, Bernard Haasdonk, and Andrea Barth.
\newblock Data-driven time parallelism via forecasting.
\newblock {\em SIAM Journal on Scientific Computing}, 41(3):B466--B496, 2019.

\bibitem{chen2014use}
Feng Chen, Jan~S Hesthaven, and Xueyu Zhu.
\newblock On the use of reduced basis methods to accelerate and stabilize the
  parareal method.
\newblock {\em Reduced Order Methods for Modeling and Computational Reduction},
  pages 187--214, 2014.

\bibitem{chinesta2017model}
Francisco Chinesta, Antonio Huerta, Gianluigi Rozza, and Karen Willcox.
\newblock Model reduction methods.
\newblock {\em Encyclopedia of Computational Mechanics Second Edition}, pages
  1--36, 2017.

\bibitem{dai2024augmented}
Xiaoying Dai, Miao Hu, Jack Xin, and Aihui Zhou.
\newblock An augmented subspace based adaptive proper orthogonal decomposition
  method for time dependent partial differential equations.
\newblock {\em Journal of Computational Physics}, 514:113231, 2024.

\bibitem{dai2020two}
Xiaoying Dai, Xiong Kuang, Jack Xin, and Aihui Zhou.
\newblock Two-grid based adaptive proper orthogonal decomposition method for
  time dependent partial differential equations.
\newblock {\em Journal of Scientific Computing}, 84:1--27, 2020.

\bibitem{dai2013symmetric}
Xiaoying Dai, Claude Le~Bris, Fr{\'e}d{\'e}ric Legoll, and Yvon Maday.
\newblock Symmetric parareal algorithms for {H}amiltonian systems.
\newblock {\em ESAIM: Mathematical Modelling and Numerical Analysis},
  47(3):717--742, 2013.

\bibitem{dai2013stable}
Xiaoying Dai and Yvon Maday.
\newblock Stable parareal in time method for first-and second-order hyperbolic
  systems.
\newblock {\em SIAM Journal on Scientific Computing}, 35(1):A52--A78, 2013.

\bibitem{farhat2003time}
Charbel Farhat and Marion Chandesris.
\newblock Time-decomposed parallel time-integrators: theory and feasibility
  studies for fluid, structure, and fluid--structure applications.
\newblock {\em International Journal for Numerical Methods in Engineering},
  58(9):1397--1434, 2003.

\bibitem{farhat-2006}
Charbel Farhat, Julien Cortial, Climène Dastillung, and Henri Bavestrello.
\newblock Time-parallel implicit integrators for the near-real-time prediction
  of linear structural dynamic responses.
\newblock {\em International Journal for Numerical Methods in Engineering},
  67(5):697--724, 2006.

\bibitem{fischer2005parareal}
Paul~F Fischer, Fr{\'e}d{\'e}ric Hecht, and Yvon Maday.
\newblock A parareal in time semi-implicit approximation of the
  {N}avier-{S}tokes equations.
\newblock In {\em Domain decomposition methods in science and engineering},
  pages 433--440. Springer, 2005.

\bibitem{gander201550}
Martin~J Gander.
\newblock 50 years of time parallel time integration.
\newblock In {\em Multiple Shooting and Time Domain Decomposition Methods:
  MuS-TDD, Heidelberg, May 6-8, 2013}, pages 69--113. Springer, 2015.

\bibitem{gander2008nonlinear}
Martin~J Gander and Ernst Hairer.
\newblock Nonlinear convergence analysis for the parareal algorithm.
\newblock In {\em Domain decomposition methods in science and engineering
  XVII}, pages 45--56. Springer, 2008.

\bibitem{gander2013analysis}
Martin~J Gander, Yaolin Jiang, Bo~Song, and Hui Zhang.
\newblock Analysis of two parareal algorithms for time-periodic problems.
\newblock {\em SIAM Journal on Scientific Computing}, 35(5):A2393--A2415, 2013.

\bibitem{gander-2026}
Martin~J Gander, Mario Ohlberger, and Stephan Rave.
\newblock A parareal algorithm with low-rank coarse solvers.
\newblock {\em arXiv:2508.08873}, 2025.

\bibitem{gander2008analysis}
Martin~J Gander and Madalina Petcu.
\newblock Analysis of a krylov subspace enhanced parareal algorithm for linear
  problems.
\newblock In {\em ESAIM: Proceedings}, volume~25, pages 114--129. EDP Sciences,
  2008.

\bibitem{gander2007analysis}
Martin~J Gander and Stefan Vandewalle.
\newblock Analysis of the parareal time-parallel time-integration method.
\newblock {\em SIAM Journal on Scientific Computing}, 29(2):556--578, 2007.

\bibitem{gander2025time}
Martin~J Gander, Shulin Wu, and Tao Zhou.
\newblock Time parallelization for hyperbolic and parabolic problems.
\newblock {\em Acta Numerica}, 34:385--489, 2025.

\bibitem{garrido2005convergent}
Izaskun Garrido, Magne~S Espedal, and Gunnar~E Fladmark.
\newblock A convergent algorithm for time parallelization applied to reservoir
  simulation.
\newblock In {\em Domain decomposition methods in science and engineering},
  pages 469--476. Springer, 2005.

\bibitem{golub2013matrix}
Gene~H Golub and Charles~F Van~Loan.
\newblock {\em Matrix computations}.
\newblock JHU press, 2013.

\bibitem{he-2010}
Liping He.
\newblock The reduced basis technique as a coarse solver for parareal in time
  simulations.
\newblock {\em Journal of Computational Mathematics}, 28:676--692, 2010.

\bibitem{ibrahim2023parareal}
Abdul~Q Ibrahim, Sebastian G{\"o}tschel, and Daniel Ruprecht.
\newblock Parareal with a physics-informed neural network as coarse propagator.
\newblock In {\em European Conference on Parallel Processing}, pages 649--663.
  Springer, 2023.

\bibitem{ibrahim2024space}
Abdul~Qadir Ibrahim, Sebastian G{\"o}tschel, and Daniel Ruprecht.
\newblock Space-time parallel scaling of parareal with a physics-informed
  fourier neural operator coarse propagator applied to the black-scholes
  equation.
\newblock {\em PASC'25: Proceedings of the Platform for Advanced Scientific
  Computing Conference}, pages 1--11, 2025.

\bibitem{iizuka2018influence}
Mikio Iizuka and Kenji Ono.
\newblock Influence of the phase accuracy of the coarse solver calculation on
  the convergence of the parareal method iteration for hyperbolic pdes.
\newblock {\em Computing and Visualization in Science}, 19(3):97--108, 2018.

\bibitem{legoll2022adaptive}
Fr{\'e}d{\'e}ric Legoll, Tony Leli{\`e}vre, and Upanshu Sharma.
\newblock An adaptive parareal algorithm: application to the simulation of
  molecular dynamics trajectories.
\newblock {\em SIAM Journal on Scientific Computing}, 44(1):B146--B176, 2022.

\bibitem{lions2001resolution}
Jacques~L Lions, Yvon Maday, and Gabriel Turinici.
\newblock R{\'e}solution d'{EDP} par un sch{\'e}ma en temps parar{\'e}el.
\newblock {\em Comptes Rendus de l'Acad{\'e}mie des Sciences-Series
  I-Mathematics}, 332(7):661--668, 2001.

\bibitem{maday2006reduced}
Yvon Maday.
\newblock Reduced basis method for the rapid and reliable solution of partial
  differential equations.
\newblock {\em Proceedings of International Conference of Mathematicians,
  European Mathematical Society}, III:1255--1270, 2006.

\bibitem{maday-2007}
Yvon Maday.
\newblock Parareal in time algorithm for kinetic systems based on model
  reduction.
\newblock {\em In High-Dimensional Partial Differential Equations in Science
  and Engineering, CRM Proceedings \& Lecture Notes}, 2007.

\bibitem{maday2010parareal}
Yvon. Maday.
\newblock The parareal in time algorithm.
\newblock In F.~Magoul{\`e}s, editor, {\em Substructuring Techniques and Domain
  Decomposition Methods}, chapter~2, pages 19--44. Saxe-Coburg Publications,
  Stirlingshire, UK, 2010.

\bibitem{maday2020adaptive}
Yvon Maday and Olga Mula.
\newblock An adaptive parareal algorithm.
\newblock {\em Journal of Computational and Applied Mathematics}, 377:112915,
  2020.

\bibitem{maday2003parallel}
Yvon Maday and Gabriel Turinici.
\newblock Parallel in time algorithms for quantum control: Parareal time
  discretization scheme.
\newblock {\em International Journal of Quantum Chemistry}, 93(3):223--228,
  2003.

\bibitem{nielsen2012feasibility}
Allan~S Nielsen.
\newblock Feasibility study of the parareal algorithm.
\newblock {\em Technical University of Denmark}, 2012.

\bibitem{obukhov1983kolmogorov}
AM~Obukhov.
\newblock {K}olmogorov flow and laboratory simulation of it.
\newblock {\em Usp. Mat. Nauk}, 38(4):101--111, 1983.

\bibitem{pamela2024neural}
S.J.P. Pamela, N.~Carey, J.~Brandstetter, R.~Akers, L.~Zanisi, J.~Buchanan,
  V.~Gopakumar, M.~Hoelzl, G.~Huijsmans, K.~Pentland, T.~James, and
  G.~Antonucci.
\newblock Neural-parareal: Self-improving acceleration of fusion mhd
  simulations using time-parallelisation and neural operators.
\newblock {\em Computer Physics Communications}, 307:109391, 2025.

\bibitem{pinnau2008model}
Ren{\'e} Pinnau.
\newblock Model reduction via {P}roper {O}rthogonal {D}ecomposition.
\newblock In {\em Model Order Reduction: Theory, Research Aspects and
  Applications}, pages 95--109. Springer, 2008.

\bibitem{ruprecht2012explicit}
D~Ruprecht and R~Krause.
\newblock Explicit parallel-in-time integration of a linear acoustic-advection
  system.
\newblock {\em Computers \& Fluids}, 59:72--83, 2012.

\bibitem{staff2005stability}
Gunnar~A Staff and Einar~M R{\o}nquist.
\newblock Stability of the parareal algorithm.
\newblock In {\em Domain decomposition methods in science and engineering},
  pages 449--456. Springer, 2005.

\bibitem{strikwerda2004finite}
John~C Strikwerda.
\newblock {\em Finite difference schemes and partial differential equations}.
\newblock SIAM, 2004.

\bibitem{verfurth2013posteriori}
R{\"u}diger Verf{\"u}rth.
\newblock {\em A posteriori error estimation techniques for finite element
  methods}.
\newblock OUP Oxford, 2013.

\bibitem{volkwein2011model}
Stefan Volkwein.
\newblock Model reduction using proper orthogonal decomposition.
\newblock {\em Lecture Notes, University of Graz}, 1025, 2011.

\bibitem{xin2016periodic}
Jack Xin, Yifeng Yu, and Andrej Zlatos.
\newblock Periodic orbits of the {ABC} flow with {A} = {B} = {C} =1.
\newblock {\em SIAM Journal on Mathematical Analysis}, 48(6):4087--4093, 2016.

\bibitem{PHG}
The toolbox {P}arallel {H}ierarchical {G}rid ({PHG}).
\newblock \url{http://lsec.cc.ac.cn/phg/}, 2019.

\end{thebibliography}
	
\end{document}